\documentclass[12pt]{article}
\usepackage[utf8]{inputenc}
\usepackage[T1]{fontenc}
\usepackage[english]{babel}

\usepackage[tbtags]{amsmath}
\usepackage{amssymb}
\usepackage{amsfonts}
\usepackage{amsthm}
\usepackage{array}
\usepackage{bbm}
\usepackage{stmaryrd}

\usepackage[svgnames]{xcolor}
\usepackage{hyperref}

\usepackage{amsmath, amsfonts}
\usepackage{enumerate}
\usepackage{color}
\usepackage{hyperref}
\numberwithin{equation}{section}

\def\B1{B_{1/2}}

\def\Bl1{{\lambda_1}}
\def\Box{\hfill\rule{2.5mm}{2.5mm}}

\def\O{{\cal {O}}}
\def\C{{\cal {C}}}

\def\H{{\cal H}}

\def\L{{\cal L}}

\def\R{{\mathbb {R}}}

\def\RR{{\cal R}}

\def\bl1{{\bar \lambda_1}}

\def\build#1_#2^#3{\mathrel{
\mathop{\kern 0pt#1}\limits_{#2}^{#3}}}

\def\h1{\mathop{\rm H^1_{\rm loc,\rm u}}}

\def\l2{\mathop{\rm L^2_{\rm loc,\rm u}}}

\newtheorem{cor}{Corollary}[section]

\newtheorem{cl}[cor]{Claim}
\newtheorem{lem}[cor]{Lemma}
\newtheorem{prop}[cor]{Proposition}
\newtheorem{nb}[cor]{Remark}
\newtheorem{propo}{Proposition}
\newtheorem{thm}[propo]{Theorem}

\newtheorem{proposi}[propo]{Proposition}

\def\C{{\cal {C}}}

\renewcommand{\H}{{\cal H}}

\renewcommand{\l}{\left}

\def\R{{\mathbb {R}}}

\theoremstyle{plain}
\newtheorem*{thm*}{Theorem}
\newtheorem*{propos*}{Proposition}

\theoremstyle{definition}

\theoremstyle{remark}

\makeatletter
\def\blfootnote{\xdef\@thefnmark{}\@footnotetext}
\makeatother

\title{\bf Prescribing the center of mass of a multi-soliton solution for a perturbed semilinear wave equation }

\author{M.A. Hamza\\
{\it \small 
Imam Abdulrahman Bin Faisal University
P.O. Box 1982 Dammam, Saudi Arabia}\\
Omar Saidi \\
{\it \small Institut Préparatoire Aux Etudes D'ingénieurs de Nabeul. Campus Universitaire, Merazka, 8000 Nabeul, Tunisie},\\
{\it \small Laboratoire \'{E}quations Aux Dérivées Partielles }
\\
{\it \small LR03ES04, 2092 Tunis, Tunisie}
}
\date{}

\allowdisplaybreaks

\begin{document}

\title{\bf Refined estimates of the blow-up profile for a strongly perturbed semilinear wave equations in one space dimension}

\maketitle

\begin{abstract}
We consider in this paper a class of strongly perturbed semilinear wave equations with a non-characteristic point in one space dimension, for general initial data. Working in the framework of similarity variables, in \cite {fh} Merle and Zaag constructed an explicit stationary solution of the unperturbed problem and proved an exponential convergence to this family of solutions. If we follow the same strategy under our strongly perturbed equation we just obtain a polynomial convergence which is a rough estimate compared to the one obtained in the unperturbed problem. In order to refine this approximation, we constructed an implicit solution to the perturbed problem which approaches the stationary solutions of the unperturbed problem and we prove the exponential convergence to this prescribed blow-up profile. 
\end{abstract}

\medskip


\medskip


\begin{section}{Introduction }
\end{section}
\subsection{Known results and motivation of the problem}

In the current work, we consider the following one dimensionel semilinear
wave equation:

\begin{equation}\label{y}
\left\{
  \begin{array}{ll}
   \partial_{t}^2 u= \partial_{x}^2 u+|u|^{p-1}u+f(u)\\
  (u(x,0),\partial_{t} u(x,0))=(u_{0}(x),u_{1}(x))\,\,\,\,\,\,\,\,\,\,\,\,\,\,\,\,\,\,\,\,\,\,\,\\
  \end{array}
\right.
\end{equation}
where $u(t):\,\, x\in {\R} \rightarrow u(x,t) \in {\R}$, $u_{0}(x)\in H^{1}_{loc,u}$, $u_{1}(x) \in L^2_{loc,u}$ and 
\begin{equation}\label{fy}
f(v)=\frac{|v|^{p-1}v}{\log^{a} (2+v^2)}\,\,\,\,\,{\rm for }\,\,{\rm all}\,\,\,\,\,\,\,v \in \R\,\,\,\,\,{\rm with} \,\,\,\,\,\,a> 1.
\end{equation}
The space $L^2_{loc,u}$ is the set of all $v\in L^2_{loc}$ such that
$$\|v\|_{L^2_{loc,u}}\equiv \sup_{\alpha\in{\R}}\Big(\int_{|x-\alpha|<1} |v(x)|^2 dx \Big)^{\frac{1}{2}}<+\infty ,$$
the space $ H^{1}_{loc,u}=\{v \mid v,|\partial_{x} v|\in L^{2}_{loc,u}\}.$ 

\bigskip
The Cauchy problem of equation $\eqref{y}$ is wellposed in $H^{1}_{loc,u}\times L^{2}_{loc,u}$. This is followed from the finite speed of propagation and the wellposdness in $H^{1}\times L^{2}$, valid whenever
$1< p <p_{S}=1+\frac{4}{N-2}$. The existence of blow-up solutions $u(t)$ of $\eqref{y}$ follows from
 ODE techniques or the energy-based blow-up criterion of Levine \cite{ha} (see
also Levine and Todorova \cite{HAG} and Todorova \cite{GT}). More blow-up results can be found in Caffarelli and Friedman \cite{LA}, \cite{LA1}, Kichenassamy and Littman \cite{SW}, \cite{SW1}, Killip, Stovall and Visan \cite{RM}.

\bigskip

Note that in this paper, we consider a class of perturbations of the idealized equation (when $f\equiv  0$). This is quite meaningful since, physical models are sometimes damped and hardly come with a pure power source term (see Whitham \cite{W1}). For more application in general relativity, see Donninger, Shlag and Soffer \cite{DSS}.

\bigskip

  If $u(t)$ is a blow-up solution of $\eqref{y}$, we define (see for example Alinhac \cite{sa} and \cite{sa1}) $\Gamma$ as the graph of a function $x \mapsto T(x)$
 such that the domain of definition of $u$ (also called the maximal influence domain)
$$D_{u}=\{ (x,t)| t< T(x)\}.$$
Moreover, from the finite speed of propagation, $T$ is a 1-Lipschitz function.
Let us first introduce the following non degeneracy condition for $\Gamma$. If we introduce
  for all $x \in \R^N$, $t \leq T(x)$ and $\delta > 0$, the cone
  \begin{equation}
   C_{x,t,\delta}=\{(\xi,\tau)\neq (x,t)| 0\leq \tau \leq t- \delta | \xi-x|\},\\
 \end{equation}
 then our non degeneracy condition is the following: $x_{0} $ is a non-characteristic point if
\begin{equation} \label{17}
 \exists \,\,\,\,\delta_{0}=\delta_{0}(x_{0})\in (0,\,1)\,\,{\rm such }\,\,{\rm that} \,\,u\,\, {\rm is}\,\, {\rm defined}\,\, {\rm on}\,\, C_{x_{0},T(x_{0}),\delta_{0}}.
 \end{equation}
If condition $\eqref{17}$ is not true, then we call $x_{0} $ a characteristic point. Already, we know from \cite {C2} and \cite {fh6} that there exist blowup solutions for characteristic points. In our paper we are concerned to the non-characteristic case. We note by $\RR$ the set of non-characteristic points. In \cite {fh} and \cite {fh1} Merle and Zaag have established the following:
\begin{itemize}
\item The set of non-characteristic points $\RR$ is non empty and open.
\item The function $T(x)$ is $\C^1$ on $\RR$ and for all $x_{0}\in \R$, $T'(x_{0})=d(x_{0})\in (-1,1)$.
\end{itemize}

\bigskip

Practically, we define for all $ x_{0} \in \R$, $0< T_{0} \leq T_{0}(x_{0})$, the following self-similar transformation introduced in Antonini and Merle \cite {cf} and used in \cite {MH}, \cite {MH1}, \cite {fh3}, \cite {kl}, \cite {fh4} and \cite {X}:
\begin{equation} \label{2}
y=\frac{x-x_{0}}{T_{0}-t},\,\,\,\,\,\,s=-\log(T_{0}-t)\,\,{\rm and }\,\, w_{x_{0},T_{0}}(y,s)=(T_{0}-t)^\frac{2}{p-1}u(x,t).\,\,\,\,\,\,\,\,\,\,\,\\
\end{equation}
The function $w_{x_{0},T_{0}}$ (we write $w$ for simplicity) satisfies the following equation for all
 $y\in (-1,1)$ and $s\geq -\log(T_{0})$:

\begin{equation}\label{1}
\partial_{s}^2w = \pounds w-\frac{2(p+1)}{(p-1)^2}w+|w|^{p-1}w-\frac{p+3}{p-1}\partial_{s}w-2y \partial^2_{y,s}w+e^{\frac{-2ps}{p-1}}f(e^{\frac{2s}{p-1}}w) 
\end{equation}
where
\begin{equation} \label{12}
 \pounds w=\frac{1}{\rho} \partial_{y}(\rho(1-y^2)\partial_{y} w)\,\,\,{\rm and}\,\,\,\rho(y)=(1-y^2)^{\frac{2}{p-1}}. 
 \end{equation} 
In the new set of variables $(y,s)$, the behavior of $u$ as $t\rightarrow T_{0}$ is equivalent
to the behavior of $w$ as $s\rightarrow +\infty.$ The equation $\eqref{1}$ will be studied in the space $\H$ defined by
\begin{equation}\label{110}
 \H =\{q=(q_{1},q_{2})|\int_{-1}^{1} \Big(q_{2}^2+(\partial_{y} q_{1})^2(1-y^2)+q_{1}^2\Big)\rho dy< +\infty\}.\\
 \end{equation} 
 Note that, in one dimensionel case the energy space is equal to $\H_{0}\times \L^2_{\rho}$, where
 \begin{equation}\label{e91}
\H_{0}=\lbrace r\in H^{1}_{loc}(-1,1)\,\,\,|\,\,\,\Vert r\Vert^2_{\H_{0}}=\int_{-1}^{1}\Big((\partial_{y}r)^2(1-y^2)+r^2\Big)\rho dy<+\infty\rbrace\\
 \end{equation} 
and $\L^r_{\rho}$ is the weighted $\L^r$ space associated with the weight $\rho$ defined in $\eqref{12}$. In the whole paper we denote
 \begin{equation}\label{90}
 F(u)=\int_{0}^{u}f(v)dv.\\
 \end{equation} 

Let us expose now some important properties and identities for our paper proved in some earlier works \cite {fh3}, \cite {kl}, \cite {fh4}.
We start by recalling that 
\begin{equation}\label{cp}
E_{0}(w(s))=\int_{-1}^{1}\Big(\frac{1}{2}(\partial_{s}w)^2+\frac{1}{2}(\partial_{y}w)^2(1-y^2)+\frac{p+1}{(p-1)^2}w^2-\frac{|w|^{p+1}}{p+1}\Big)\rho  {\mathrm{d}}y
\end{equation}
is a Lyapunov functional of equation $\eqref{1}$ when $f\equiv 0$ which is defined in $\H$.
\\Then, we introduce the following functional: 
  \begin{equation} \label{e90}
E(w(s),s)=E_{0}(w(s))-e^{\frac{-2(p+1)s}{p-1}} \int_{-1}^{1} F(e^{\frac{2s}{p-1}}w)\rho {\mathrm{d}}y-\frac{1}{s^{\frac{a+1}{2}}}\int_{-1}^{1} w\partial_{s}w \rho dy.
\end{equation}
In our case and more generally under the assumptions 
$$(H_{f})\,\,\,\,\,\,\,\,\,|f(v)|\leq M\Big(1+\frac{|v|^{p}}{\log^{a} (2+v^2)}\Big),\,\,\,\,\,\,\,\,\,\,\,\,\,\,\,{\rm for }\,\,{\rm all} \,v\in \R\,\,\,\,\,\,\,\,\,\,{\rm with} \,\,\,(M > 0,\,\,\,a> 1)$$
and 
$$(H_{g})\,\,\,\,\,\,\,\,\,|g(x,t,v,z)|\leq  M(1+|z|), \,\,\,\,\,\,\,\,\,\,\,\,\,{\rm for }\,\,{\rm all} \,x,\,\,v\in \R^N\,t,\,\,z\in \R\,\,\,{\rm with }\,\,\,(M > 0),$$
with $f$ and $g$ are perturbed terms added to equation $\eqref{y}$, we have proved in \cite{X}, \cite{XX} and \cite{XXX} that 
$$H(w(s),s)=\exp\Big(\frac{p+3}{(a-1)s^{\frac{a-1}{2}}}\Big) E(w(s),s)+\theta e^{\frac{-(p+1)s}{p-1}},$$
(for a large constant $\theta $) is a Lyapunov functional of equation $\eqref{1}$ which is defined in $\H$. Based on this functional and some energy estimates, we have established when $x_{0} $ is a non-characteristic point (in the sense $\eqref{17}$) that there exists $\widehat{s}_{0}=\widehat{s}_{0}(x_{0},p,a,T_{0}(x_{0}))$ such that for all $s\geq \widehat{s}_{0}$, the following estimate holds:
 \begin{equation}\label{900}
 0 < \varepsilon_{0} \leq \|w_{x_{0},T_{0}(x_{0})}(s)\|_{H^{1}(-1,1)} + \|\partial_{s}w_{x_{0},T_{0}(x_{0})}(s)\|_{L^{2}(-1,1)}\leq M_{1}.
\end{equation}

\begin{nb}
We present now some comments on the parameter $a$ related to the blow-up rate and the blow-up limit. Under the assumptions $(H_{f})$ and $(H_{g})$ and following our earlier work \cite {X}, we have proved the result $\eqref{900}$ when the exponent $p$ is subconformal, i.e. $(1<p<1+\frac{4}{N-1},\,{\rm when }\,N\geq 2)$ for $a>1$ and in \cite {XX} for $a>2$ when the exponent $p$ is conformal, i.e. $(p\equiv 1+\frac{4}{N-1})$. The method used in \cite {X} breaks down when $a\in (0,\,1]$, within some analysis we find an equality of type $\frac{d}{ds}(E(w(s),s))\leq \frac{C}{s^{a}}E(w(s),s)$ with $E(w(s),s)$ is defined in $\eqref{e90}$ and this is a major reason preventing us from deriving the result in the case $a\in (0,\,1]$ and explains the restriction of the parameter $a$ to be in $(1,+\infty )$ in this current paper. Except when the perturbed term $f$ satisfies the condition $\eqref{fy}$, the above estimate will be refined into $\frac{d}{ds}(E(w(s),s))\leq \frac{C}{s^{a+1}}E(w(s),s)$ which implies that we can solve the question of the blow-up rate also when $a\in (0,\,1]$. In this direction, we mention the work of Nguyen and Zaag \cite {rhv} where the authors added the same perturbed term $f$ to the semilinear heat equation and they allowed values of $a$ in $(0,\,1]$ at the expense of taking the particular form $\eqref{fy}$ of the perturbation $f$ via the derivation of a suitable Lyapunov functional. Going back to the problem of the blow-up limit, in this paper we solve the case when $a>1$. It is very interesting to answer the question when $a$ in $(0,\,1]$. For that purpose we may need some refinement to the polynomial decay and to the asymptotic behavior obtained in Claim $\ref{T14}$. We felt that the study of this case can be our next challenge.
\end{nb}
A natural question then is to know if $w_{x_{0}}$ has a limit or not, as $s\longrightarrow +\infty$ (that is as $t\longrightarrow T_{0}$). In order to make a simpler presentation, we start by the case $f\equiv 0$. This case was treated by Merle and Zaag \cite {fh} where the authors have proved the convergence of the solution $w_{x_{0}}$ to the set of stationary solutions $S\equiv \lbrace 0,\,\,\kappa (d,.),\,\,-\kappa (d,.)\,\,|\,\,|d|<1\rbrace$ in one space dimension. In higher dimensions, $N\geq 2$ there is no classification of selfsimilar solutions of equation $\eqref{1}$ when ($f\equiv 0$). In other words, we already know that $\kappa (d,\omega .y)$ defined in $\eqref{400}$ is $\H_{0}$ stationary solution of equation $\eqref{1}$ when ($f\equiv 0$) for any $|d|<1$ and $\omega \in \R^N$ with $|\omega |=1$, but we are unable to say whether there are other stationary solutions or not. Despite that, in higher dimensions Merle and Zaag extended in \cite {fh15} the oppeness of the set of non-characteristic points and regularity of the blow-up curve for the semilinear wave equation in one space dimension to the higher dimensions. In a companion paper \cite {fh10}, Merle and Zaag studied the dynamic of the solution of $\eqref{1}$ when ($f\equiv 0$) near explicit stationary solutions in similarity variables and they extended some properties of the one dimension to the higher dimension. The proof of the convergence in higher dimension is far from being a simple adaptation of the one dimensional case. Indeed, several difficulties arise when $N\geq 2$, among them, we have $(N-1)$ new degenerate directions in the linearized operator of equation $\eqref{1}$ when $f\equiv 0$ around the stationary solution $\kappa (d,y)$ defined in $\eqref{400}$. These new directions naturally come from the derivative of $\kappa (d,y)$ with respect to $(N-1)$ angular directions of $d$.
\\As announced above, we want to briefly recall the work of Merle and Zaag \cite {fh} in the case when $f\equiv 0$ in a clear way. Let us first classify all the $\H_{0}$ stationary solutions of $\eqref{1}$ when $f\equiv 0$ in one dimension. More precisely, we have:
\\If $w\in \H_{0}$ a stationary solution of $\eqref{1}$ when $f\equiv 0$, then, either $w=0$ or there exist $d\in (-1,1)$ and $\omega \in \lbrace -1,\,1\rbrace$ such that $w(y)= \omega\kappa (d,y),$ where  
\begin{equation} \label{400}
\forall\,\,\,(d,y)\in (-1,1)^2\,\,\,\kappa (d,y)=\kappa_{0}\frac{(1-d^2)^{\frac{1}{p-1}}}{(1+dy)^{\frac{2}{p-1}}}\,\,\,{\rm with}\,\,\,\kappa_{0}=\Big(\frac{2(p+1)}{(p-1)^2}\Big)^{\frac{1}{p-1}}.\\
\end{equation}
In a first time, (when $f\equiv 0$) Merle and Zaag in \cite {fh} show that $w_{x_{0}}(y,s)$ defined in $\eqref{2}$ strongly converges in $H^1 \times L^2(-1,\,1)$ to a non-null connected component of the stationary solutions. More precisely, they have proved the following:
\begin{propos*}{\bf (Strong convergence related to the set of approximate stationary solution).}
Consider $w\in \C([s^*,\infty ),\H)$ for some $s^*\in \R$ a solution of equation $\eqref{1}$ (when $f\equiv 0$).
If $x_{0}$ is a non-characteristic point (in the sense $\eqref{17}$), then there exists $\omega^* (x_{0})\in \lbrace -1,1\rbrace$ such that:
$$\inf_{|d|<1}\Vert w_{x_{0}}(.,s)-\omega^*(x_{0})\kappa (d,.)\Vert_{H^1(-1,1)}+\Vert \partial_{s}w_{x_{0}}(.,s)\Vert_{L^2(-1,1)}\longrightarrow 0\,\,\,\,{\rm as}\,\,\,s\longrightarrow \infty .$$
\end{propos*}

\begin{nb}
We would like to emphasize that this proposition holds also when $w_{x_{0}}(y,s)$ defined in $\eqref{2}$ is a solution of $\eqref{1}$ in the perturbed case. In other words, this proposition remains unchanged even under the strong perturbation $f$ defined in $\eqref{fy}$ and can be generalized with the same arguments used in \cite {fh} without any difficulty. Since the perturbed term $f$ is polynomially small, it can be ignored in the proof.
\end{nb}
Then, Merle and Zaag \cite {fh} derive the following result:
\begin{thm*}{\bf (Trapping near the set of non zero stationary solutions of $\eqref{1}$ when $f\equiv 0$).} There exist positive constants $\epsilon_{0}$, $\mu_{0}$ and $C$ such that if 
$w\in \C([s^*,\infty ),\H)$, for some $s^*\in \R$ is a solution of equation $\eqref{1}$ (when $f\equiv 0$) such that:
\begin{equation}\label{1110}
\,\,\forall\,\,s\geq s^*\,\,\,E_{0}(w(s))\geq  E_{0}(\kappa_{0}),
\end{equation}
and 
\begin{equation}\label{1120}
\Big \Vert \left(
  \begin{array}{ccc}
    w(s^*) \\
    \partial_{s} w(s^*)\\
  \end{array}
\right)-\omega^*\left(
  \begin{array}{ccc}
    \kappa (d^*,.) \\
    0\\
  \end{array}
\right)\Big\Vert_{\H}\leq \epsilon^*,
\end{equation}
for some $d^*\in (-1,1)$, $\omega^*\in \lbrace -1,1\rbrace$ and $\epsilon^*\in (0,\epsilon_{0}]$, where $\H$ and its norm are defined in $\eqref{110}$ and $\kappa (d,y)$ is defined in $\eqref{400}$, then there exists $d_{\infty}\in (-1,1)$ such that 
$$|d_{\infty}-d^*|\leq C\epsilon^*(1-d^{*2}),$$
and for all $s\geq s^*$,
 \begin{equation}\label{1150}
\Big\Vert \left(
  \begin{array}{ccc}
    w(s) \\
    \partial_{s} w(s)\\
  \end{array}
\right)-\omega^*\left(
  \begin{array}{ccc}
   \kappa (d_{\infty},.) \\
   0\\
  \end{array}
\right)\Big\Vert_{\H}\leq  C\epsilon^*e^{-\mu_{0}(s-s^*)}.
\end{equation}
\end{thm*}

The proof of $\eqref{1150}$ is performed in the framework of similarity variables defined in $\eqref{2}$.
To prove this result, Merle and Zaag \cite {fh}, among other techniques, linearized $\eqref{1}$ when $f\equiv 0$ around the solution $\kappa (d,y)$ and got an exponential decay. 
Similarly to the problem of the blow-up rate, one can ask whether the result proved by Merle and Zaag in \cite {fh} holds also for the strongly perturbed semilinear wave equation. Before studying our class of perturbation we should recall that, in \cite {MH2}, Hamza and Zaag have proved a similar result as the one obtained in the case when $f\equiv 0$ under the assumption $(H_{g})$ and the more restrictive assumption:
$$|f(v)|\leq M(1+|v|^q),\,\,\,{\rm for }\,\,\,{\rm all} \,\,\,v\in \R\,\,\,{\rm with} \,\,\,(M > 0,\,\,\,q<p).$$
As we said earlier, in this paper we are conserned with the non-characteristic case and the function $f$ is defined in $\eqref{fy}$. Unlike the work of Merle and Zaag \cite {fh} and also Hamza and Zaag \cite {MH2}, we focus here on the assumption $\eqref{fy}$ which is much more complicated, so we need to invent a new idea in order to obtain the exponential decay. In fact, if we do the same as in \cite {fh} and \cite {MH2} we may obtain some terms like $\frac{1}{s^{a}}$ coming from the strong perturbation $f$ defined in $\eqref{fy}$ and we may not able to control these terms. More precisely, we obtain the following
\begin{proposi}\label{tawtaw}
Let $x_{0}\in \RR$, there exist positive constants $\epsilon_{0}$, $\mu$ and $C$ such that if 
$w\in \C([s^*,\infty ),\H)$, for some $s^*\in \R$, is a solution of equation $\eqref{1}$ such that:
\begin{equation}\label{11110}
\,\,\forall\,\,s\geq s^*\,\,\,E(w(s))\geq  E_{0}(\kappa_{0})-\frac{C}{s^{\frac{a+1}{2}}}
\end{equation}
and 
\begin{equation}\label{11120}
\Big\Vert \left(
  \begin{array}{ccc}
    w(s^*) \\
    \partial_{s} w(s^*)\\
  \end{array}
\right)-\omega^*\left(
  \begin{array}{ccc}
    \kappa (d^*,.) \\
    0\\
  \end{array}
\right)\Big\Vert_{\H}\leq \epsilon^*,
\end{equation}
for some $d^*\in (-1,1)$, $\omega^*\in \lbrace -1,1\rbrace$ and $\epsilon^*\in (0,\epsilon_{0}]$, where $\H$ and its norm are defined in $\eqref{110}$ and $\kappa (d,y)$ is defined in $\eqref{400}$, then there exists $d_{\infty}\in (-1,1)$ such that 
$$|d_{\infty}-d^*|\leq C\epsilon^*(1-d^{*2}),$$
and for all $s\geq s^*$,
\begin{equation}\label{nao}
\Big\Vert \left(
  \begin{array}{ccc}
    w_{x_{0}}(s) \\
    \partial_{s} w_{x_{0}}(s)\\
  \end{array}
\right)-\omega^*\left(
  \begin{array}{ccc}
    \kappa (d_{\infty},.) \\
    0\\
  \end{array}
\right)\Big\Vert_{\H}\leq \frac{C}{s^{\mu}},
\end{equation}

\end{proposi}

\begin{nb}
 Proposition $\ref{tawtaw}$ may be regarded as a rough version of our essential goal which
 is written for $f$ defined in $\eqref{fy}$. In order to avoid unnecessary repetition, we kindly refer the reader to \cite {fh} for all the projections of the terms in $\eqref{63}$ not involving $f$ and we will explain briefly in this remark the terms with $f$. Remarking that $\kappa (d,y)$ is not a solution of $\eqref{1}$, we can see that if we linearize $\eqref{1}$ around $\kappa (d,y)$, we get a remaining term of type $|e^{\frac{-2ps}{p-1}}f(e^{\frac{2s}{p-1}}(\kappa (d,y)+q_{1}))|$ (with $q_{1}$ will be defined later in $\eqref{301}$) provided from $f$. If we control the projection of this remaining term we can see that the presence of the additional term $\frac{C}{s^{a}}$ is natural which make the proof of this proposition very easy.
\end{nb}

Our aim in this paper is to obtain a sharp estimate for the prescribed blow-up profile, more precisely an exponential convergence like the one obtained in the case when $f\equiv 0$, which is more advantageous than the polynomial decay. From the above remark the remaining term would act as a forcing term, preventing us from obtaining the exponential decay found in Merle and Zaag \cite {fh} in the case when $f\equiv 0$. In order to overcome this difficulty, instead of linearzing around the explicit function $\kappa (d,y)$, which is not a solution of equation $\eqref{1}$, (it is just an approximate solution) we will linearize around a new implicit profile function which happens to be an exact solution of equation $\eqref{1}$ defined by: 
\begin{equation}\label{87}
\overline{w}_{1}(d,y,s)=\kappa (d,y)\frac{\phi (s-\log(\frac{1+dy}{\sqrt{1-d^2}}))}{\kappa_{0}},
\end{equation}
with $\phi$ is a solution of the associated ODE to the PDE $\eqref{1}$. Moreover, we introduce 
\begin{equation}\label{e505}
\overline{w}(d,y,s)= \left(
  \begin{array}{ccc}
    \overline{w}_{1} (d,y,s)\\
    \overline{w}_{2}(d,y,s)\\
  \end{array}
\right)\,\,\,{\rm where }\,\,\,\overline{w}_{2}(d,y,s)=\partial_{s}\overline{w}_{1}(d,y,s).
\end{equation}
The most important property of our new solution is that
\begin{equation}\label{TAB}
\phi(s)-\kappa_{0} \sim  -\frac{\kappa_{0}}{p-1}\Big(\frac{p-1}{4s}\Big)^{a}\,\,\,{\rm as}\,\,s \rightarrow +\infty.
\end{equation}
This property is crucial in many steps in this paper. We will see later in Appendix $\ref{T160}$ the proof of the equivalence $\eqref{TAB}$ as well as some complementary results and in Appendix $\ref{TS2}$ we will see the details of the construction of $\overline{w}(d,y,s)$ defined in $\eqref{87}$ and $\eqref{e505}$, where we can see this property in a clear way. In addition to that, from the proposition of the case when $f\equiv 0$ and the remark after, we can see also that $w_{x_{0}}(y,s)$, defined in $\eqref{2}$ is a solution of $\eqref{1}$ in the perturbed case approaches $\overline{w}$ defined in $\eqref{e505}$ strongly in $H^1 \times L^2(-1,\,1)$ norm. More precisely, we write, when $x_{0}$ is a non-characteristic point (in the sense $\eqref{17}$) and $\omega^*=\omega^*(x_{0})\in \lbrace -1,1\rbrace $, the following:
\begin{equation}\label{waw1}
\hspace{-0.5cm}\inf_{|d|<1}\Vert w_{x_{0}}(.,s)-\omega^*\overline{w}_{1} (d_{\infty},.,s)\Vert_{H^1}+\Vert \partial_{s}w_{x_{0}}(.,s)-\overline{w}_{2}(d_{\infty},.,s)\Vert_{L^2}\longrightarrow 0\,\,{\rm as}\,\,s\longrightarrow \infty .
\end{equation}
 Our aim in this paper is to prove that $\overline{w}_{1} $ is the blow-up profile such that if $w$ is a solution of equation $\eqref{1}$, we show the exponential convergence of $w$ to this expected profile.

 Before stating our result, we would like to mention that an extensive literature is devoted to the question of the construction of a solution to some PDE. We start by the work of C\^{o}te and Zaag \cite{C2} for the non-perturbed problem, where the authors for any integer $k\geq 2$ and $\zeta_{0} \in \R$ constructed a blow-up solution with a characteristic point $a$, such that the asymptotic behavior of the solution near $(a,T(a))$ shows a decoupled sum of $k$ solitons with alternative signs, whose centers (in the hyperbolic geometry) have $\zeta_{0}$ as a center of mass, for all times. Moreover, Hamza and Zaag \cite {MH3} extended the result of C\^{o}te and Zaag \cite{C2} in the characteristic case where they added a perturbed terms satisfying the hypotheses $(H_{f})$ and $(H_{g})$ and they prescribed the center of mass of a multi-soliton solution for strongly perturbed semilinear wave equation. That question was investigated by Nguyen and Zaag in \cite{VH} where the authors construct an implicit profile function for a strongly perturbed semilinear heat equation and by Bressan in \cite{BR} and \cite{BR1} for the semilinear heat equation with an exponential term source. Also we would like to mention the remarkable result of Ghoul, Nguyen and Zaag in \cite{TVH} where the authors constructed blow-up solutions for non-variational semilinear parabolic system and the constructed solutions are stable under a small perturbation of initial data. Then in \cite{TVH1} the same authors extend their result to a higher order semilinear parabolic equation. Concerning the question of the construction of solutions to some PDE, we note that Mahmoudi, Nouaili and Zaag in \cite{MNZ} constructed a $2\pi$-periodic solution to the nonlinear heat equation with power nonlinearity in one space dimension which blows up in finite time $T$ only at one blow-up point. In the same sense Nouaili and Zaag in \cite{NH1} constructed a solution for the complex nonlinear heat equation which blows up in finite time $T$ only at a one blow-up point, the same type of solution is constructed for the complex Ginzburg-Landau Equation by Zaag in \cite{Z}, by Masmoudi and Zaag in \cite{MAZ} and also by Nouaili and Zaag in \cite{NH} in the critical case. Willing to be as exhaustive as possible in our bibliography about the question of the construction, we would like to mention that this question was solved for (gKdV) by C\^{o}te in \cite{C} and \cite{C1}, for Shr\"{o}dinger maps by Merle, Rapha\"{e}l and Rodnianski in \cite{MRR}, the wave maps by Ghoul, Ibrahim and Nguyen in \cite{GIN} and for the Keller-Segel model by Rapha\"{e}l and Schweyer in \cite{RS} and also Ghoul and Masmoudi in \cite{GHM}.
\begin{nb}
 We believe that the present work can be a good start to understand the dynamics of many PDE's without explicit solution. We can cite for example the PDE treated by Hamza and Zaag \cite{HZN} in one space dimension and in \cite{HZNn} and \cite{HZNnn} for the higher dimensions, where we don't have any stationary solution to the semilinear wave equation and heat equation with nonlinearity including a logarithmic factor which is not scale invariant. In our opinion, we can contruct an implicit solution adapted to the problem studied by the authors and prove a similar result as the one obtained in the present paper.
\end{nb}

\subsection{Main results and strategy of the proof}

Let us give now our main result. After linearizing around $\overline{w}_{1}$, which is the major novelty in our approach we derive, with more works,  the following theorem:
\begin{thm}\label{T4}{\bf (Trapping near the set of the family of the implicit profile $\overline{w}$).} There exist $S_{1}=S_{1}(x_{0},p,a)$, $\epsilon_{0}$, $\mu_{0}$ and $C$ such that if 
$w\in \C([s^*,\infty ),\H)$ for some $s^*\geq S_{1}$ a solution of equation $\eqref{1}$ such that:
\begin{equation}\label{111}
\,\,\forall\,\,s\geq s^*\,\,\,{\rm and }\,\,\,E(w(s),s)\geq  E_{0}(\overline{w}_{1}(d,y,s) )-\frac{C}{s^{\frac{a+1}{2}}}
\end{equation}
and 
\begin{equation}\label{112}
\Big\Vert \left(
  \begin{array}{ccc}
    w(s^*) \\
    \partial_{s} w(s^*)\\
  \end{array}
\right)-\omega^*\left(
  \begin{array}{ccc}
    \overline{w}_{1}(d^*,.,s^*) \\
    \overline{w}_{2}(d^*,.,s^*)\\
  \end{array}
\right)\Big\Vert_{\H}\leq \epsilon^*,
\end{equation}
for some $d^*\in (-1,1)$, $\omega^*\in \lbrace -1,1\rbrace$ and $\epsilon^*\in (0,\epsilon_{0}]$, where $\H$ and its norm are defined in $\eqref{110}$, $\overline{w}_{1}$ and $\overline{w}_{2}$ are defined in $\eqref{87}$ and $\eqref{e505}$, then there exists $d_{\infty}\in (-1,1)$ such that 
$$|d_{\infty}-d^*|\leq C\epsilon^*(1-d^{*2})$$
and for all $s\geq s^*$,
 \begin{equation}\label{115}
\Big\Vert \left(
  \begin{array}{ccc}
    w(s) \\
    \partial_{s} w(s)\\
  \end{array}
\right)-\omega^*\left(
  \begin{array}{ccc}
    \overline{w}_{1}(d_{\infty},.,s) \\
   \overline{w}_{2}(d_{\infty},.,s)\\
  \end{array}
\right)\Big\Vert_{\H}\leq  C\epsilon^* e^{-\mu_{0}(s-s^*)}.
\end{equation}
\end{thm}
\begin{nb}\label{90b}
Thanks to $\eqref{waw1}$, condition $\eqref{112}$ is also guaranteed and the time $s^*$ is completely explicit and characterized by the fact that:
$$s^*=\inf_{s\geq -\log (T(x_{0}))} \inf_{|d|<1}\Big\Vert\left(
  \begin{array}{ccc}
    w(s) \\
    \partial_{s} w(s)\\
  \end{array}
\right)-\omega^*\left(
  \begin{array}{ccc}
    \overline{w}_{1}(d^*,.,s) \\
    \overline{w}_{2}(d^*,.,s)\\
  \end{array}
\right)\Big\Vert_{\H}\leq \epsilon_{0}.$$
(For more detail of this fact we can see Corollary $4$ and the remark after that corollary in Merle and Zaag \cite {fh}). 

\end{nb}

\begin{nb}
 We would like to compare $\eqref{111}$ with $\eqref{1110}$ the condition imposed by Merle and Zaag in the unpertubed case, we can remark that we have an additional polynomially small term. This comes from the fact that the difference between $E(w(s),s)$ and $E_{0}(w(s))$ is polynomially small. In addition to that, we can replace $E_{0}(\overline{w}_{1}(d,y,s) )$ in $\eqref{111}$ by $E_{0}( \kappa_{0})$ thanks to the monotonicity of $E_{0}(w(s))$, the fact that $\overline{w}_{1}(d,y,s)\sim\kappa (d,y)$ as $s\longrightarrow \infty $ and the fact that $E_{0}( \kappa (d,y))=E_{0}( \kappa_{0})$ (for more details about the last equality, the interested reader may consult Merle and Zaag \cite {fh}). 
\end{nb}
\begin{nb}
Our proof under the assumption $(H_{g})$ remains valid with exactly the same ideas and purely technical differences, that we omit to keep this paper in reasonable length. 
\end{nb}
Let us now comment on the method used to prove our results. Noting that the proof of Theorem $\ref{T4}$ is far from being a simple adaptation of the case when $f\equiv 0$ treated by Merle and Zaag \cite {fh}. Indeed, there are additional difficulties arising from the perturbed term $f$ and the linearization around the new solution $\overline{w}_{1}$ defined in $\eqref{87}$ which makes the technical details harder to elaborate. Accordingly, the exponential decay obtained in Theorem $\ref{T4}$ requires more works than the case when $f\equiv 0$ and special arguments where we need to invent new idea to get it. However, in a first time, we follow the same strategy of the case when $f\equiv 0$ to get a rough estimate where the condition $\eqref{111}$ will play an important role to obtain the polynomial decay. Based upon this estimate we can see that $\Vert q(s)\Vert_{\H}\longrightarrow 0$ as $s\,\,\longrightarrow \infty$ with $q$ will be defined later in $\eqref{301}$. Even better, thanks to the information asserted just before with some further refinement of the result obtained in the polynomial decay we conclude the exponential convergence at the end of this paper and we get our main Theorem $\ref{T4}$.
\\Let us mention briefly the structure of the paper: The article is organized around three main results. The modulation theory which is a major step to obtain the polynomial decay and the exponential decay, we present the modulation theory in Section $\ref{TS3}$, the polynomial and the exponential decay in Section $\ref{TS2222}$ and each of them in a separate subsection, where we conclude the proof of Theorem $\ref{T4}$. 
\\We mention that $C$ depends on $p$, $a$ and $x_{0}$ will be used in all this paper to denote a positive constant which may vary from line to line.
\begin{section}{Modulation theory}\label{TS3}
\end{section}
In this part, we work in the space $\H$ defined in $\eqref{110}$, which is a natural choice (the energy space in $w$). We consider for some $s^*$ defined in Remark $\ref{90b}$, the function $w\in \C([s^*,\infty ),\H)$ a solution of equation $\eqref{1}$, where $w$ may be equal to $w_{x_{0}}$ defined in $\eqref{2}$ from $u$ a blow-up solution to equation $\eqref{y}$, with no restriction to $x_{0}$.
\\In this section, we use modulation theory and introduce a parameter $d(s)$ adapted to the dispersive property of the equation $\eqref{1}$ whenever $\eqref{112}$ holds, in order to obtain the polynomial decay in the next section.

\bigskip

In the begining, we are going to call back some tools introduced in Merle and Zaag \cite {fh} in the case when $f\equiv 0$.
We briefly recall the crucial linear operator $L_{d}$ and its spectral properties introduced in Merle and Zaag \cite{fh} defined in the energy space $\H$ by:
\begin{equation}\label{60}
L_{d}\left(
  \begin{array}{ccc}
    q_{1}(y,s) \\
    q_{2}(y,s)\\
  \end{array}
\right)=\left(
  \begin{array}{ccc}
    q_{2} \\
    \pounds q_{1}+\psi(d,y)q_{1}-\frac{p+3}{p-1}q_{2}-2y\partial_{y}q_{2}\\
  \end{array}
\right),
 \end{equation}
where
\begin{equation}\label{38}
 \psi(d,y)=p\kappa(d,y)^{p-1}-\frac{2(p+1)}{(p-1)^2}=\frac{2(p+1)}{(p-1)^2}\Big(\frac{p(1-d^2)}{(1+dy)^2}-1\Big).
 \end{equation}
The linear operator $L_{d}$ will play a fundamental role in our analysis, it is useful to recall some well known results link with the linear operator $L_{d}$. 
\\{\bf Spectral properties related to $L_{d}$.}
\\It is well known that $\lambda =1$ and $\lambda =0$ are the positive eigenvalues of the linear operator $L_{d}$ and the rest of the eigenvalues are negative. It happens that the corresponding eigenfunctions of the positive eigenvalues are:
\begin{equation}\label{61}
F_{1}^{d}=(1-d^2)^{\frac{p}{p-1}}\left(
  \begin{array}{ccc}
   (1+dy)^{-\frac{2}{p-1}-1} \\
    (1+dy)^{-\frac{2}{p-1}-1}\\
  \end{array}
\right)\,{\rm and}\,F_{0}^{d}=(1-d^2)^{\frac{1}{p-1}}\left(
  \begin{array}{ccc}
   \frac{y+d}{(1+dy)^{\frac{2}{p-1}+1}} \\
    0\\
  \end{array}
\right).
\end{equation}
Moreover, it holds for some $C>0$ and any $\lambda \in \lbrace 0,1 \rbrace $ that 
\begin{equation}\label{62}
\forall \,\,|d|<1,\,\, \Vert F_{\lambda}^{d}\Vert_{\H}+(1-d^2)\Vert \partial_{d}F_{\lambda}^{d}\Vert_{\H}\leq C.
\end{equation}
In order to compute the projectors on the eigenfunctions of $L_{d}$, we consider its conjugate with respect to the natural inner product $\Upsilon$ of $\H$ defined by:
\begin{equation}\label{50}
\Upsilon (q,r)=\int_{-1}^{1}(q_{1}(-\pounds r_{1}+r_{1})+q_{2}r_{2})\rho dy.\\
\end{equation}
The computation of $L_{d}^*$ the conjugate operator of $L_{d}$ with respect to $\Upsilon$ is simple but lengthy that we omit, (for more details we kindly address the reader to see Lemma 4.1 page 81 in Merle and Zaag \cite{fh}).
Furthermore, $L_{d}^*$ has two nonnegative eigenvalues with eigenfunctions $W_{\lambda}^{d}$ such that:
\begin{equation}\label{e4}
W_{1,2}^{d}(y)=c_{1}\frac{1-y^2}{(1+dy)^{\frac{2}{p-1}+1}},\,\,\,W_{0,2}^{d}(y)=c_{0}\frac{y+d}{(1+dy)^{\frac{2}{p-1}+1}},
\end{equation}
$W_{\lambda,1}^{d}$ is uniquely determined by 
\begin{equation}\label{e5}
-\pounds r+r=\Big(\lambda -\frac{p+3}{p-1}\Big)r_{2}-2y\partial_{y}r_{2}+\frac{8}{p-1}\frac{r_{2}}{1-y^2}
\end{equation}
with $r_{2}=W_{\lambda,2}^{d}$ and the $\C^1$ function $c_{\lambda}$ fixed by the relation 
\begin{equation}\label{e6}
\Upsilon (W_{\lambda}^{d},F_{\lambda}^{d})=1.
\end{equation}
Finally, we also introduce for $q\in \H$ and for $\lambda =1$ and $\lambda =0$ the following
\begin{equation}\label{35}
\pi_{\lambda}^{d}(q)=\Upsilon (W_{\lambda}^{d},q)\,\,\,{\rm and}\,\,\,q=\pi_{1}^{d}(q)F_{1}^{d}(y)+\pi_{0}^{d}(q)F_{0}^{d}(y)+\pi_{-}^{d}(q).\\
\end{equation}
\begin{nb}
Let us notice $\pi_{\lambda}^{d}(q)$ is the projection of $q$ on the eigenfunction of $L_{d}$ associated to $\lambda $ and that $\pi_{-}^{d}(q)$ is the negative part of $q$ such that:
\begin{equation}\label{18}
\pi_{-}^{d}(q)\in \H_{-}^{d}=\lbrace r\in \H\,\,\,|\pi_{1}^{d}(r)=\pi_{0}^{d}(r)=0\rbrace.\\
\end{equation}
\end{nb} 
Of course, we have the following orthogonality results:
\begin{enumerate}
\item {\bf (Orthogonality)} For all $|d|<1$ and $\lambda \in \lbrace 0,1 \rbrace $, we have $\Upsilon (W_{\lambda}^{d},F_{1-\lambda}^{d})=0$.
\item {\bf (Normalization)} There exists $C>0$ such that for $\lambda \in \lbrace 0,1 \rbrace $ and $|d|<1$,
\begin{equation}\label{e7}
\Vert W_{\lambda}^{d}\Vert_{\H}+(1-d^2)\Vert\partial_{d} W_{\lambda}^{d}\Vert_{\H}\leq C.
\end{equation}
\end{enumerate}

\bigskip

After this reminder, we start the modulation theory and we claim the following:
\begin{prop}\label{T7}{\bf (Modulation of $w$ with respect to $\overline{w}_{1}(d,y,s)$).}
There exist $\epsilon_{1}>0$ and $K>0$ such that if $(w,\partial_{s}w)\in \C([s^*,+\infty ),\H)$ is a solution to equation $\eqref{1}$ which satisfies $\eqref{112}$ for some $|d^*|<1$, $\omega^* \in \lbrace -1,1\rbrace $ and $\epsilon^*\leq \epsilon_{1}$, then the following is true:
 \begin{enumerate}[{\rm(i)}]
\item {\bf (Choice of the modulation parameter)} There exists $d(s)\in \C^1([s^*,+\infty ),(-1,1))$ such that for all $s\in [s^*,+\infty )$,
\begin{equation}\label{101}
\pi_{0}^{d(s)}(q(s))=0,\\
\end{equation}
where $\pi_{0}^{d}$ is defined in $\eqref{35}$, $q=(q_{1},q_{2})$ is defined for all $s\in [s^*,+\infty )$ by 
\begin{equation}\label{301}
\left(
  \begin{array}{ccc}
    w(y,s) \\
    \partial_{s}w(y,s)\\
  \end{array}
\right)=\left(
  \begin{array}{ccc}
    \overline{w}_{1}(d(s),y,s) \\
    \overline{w}_{2}(d(s),y,s)\\
  \end{array}
\right)+\left(
  \begin{array}{ccc}
    q_{1}(y,s) \\
    q_{2}(y,s)\\
  \end{array}
\right).\\
\end{equation}
Moreover,
\begin{equation}\label{NOUV50}
\Big|\log \Big (\frac{1+d(s^*)}{1-d(s^*)}\Big )-\log \Big (\frac{1+d^*}{1-d^*}\Big )\Big|+\Vert q(s^*)\Vert_{\H}\leq K \epsilon^*.
\end{equation}
\item {\bf (Equation on $q$)}
For all $s\in [s^*,+\infty )$
$$
\frac{\partial}{\partial_{s}}\left(
  \begin{array}{ccc}
    q_{1}(y,s) \\
    q_{2}(y,s)\\
  \end{array}
\right)=\overline{L}_{d(s)}\left(
  \begin{array}{ccc}
    q_{1}(y,s) \\
    q_{2}(y,s)\\
  \end{array}
\right)-d' \left(
  \begin{array}{ccc}
   \partial_{d} \overline{w}_{1}(d,y,s) \\
   \partial_{d} \overline{w}_{2}(d,y,s)\\
  \end{array}
\right)\\
 $$

\begin{equation}\label{21}
+\left(
  \begin{array}{ccc}
    0 \\
    h(d,y,q_{1})\\
  \end{array}
\right)+\left(
  \begin{array}{ccc}
    0 \\
    \widehat{f}(\overline{w}_{1},q_{1},s)\\
  \end{array}
\right),
 \end{equation}

where

\begin{equation}
 \overline{L}_{d(s)}\left(
  \begin{array}{ccc}
    q_{1} \\
    q_{2}\\
  \end{array}
\right)\hspace{-0,1cm}=\hspace{-0,1cm}\left(
  \begin{array}{ccc}
    q_{2} \\
   \pounds q_{1}\hspace{-0,2cm}+(\overline{\psi}(d,y,s)\hspace{-0,2cm}+e^{-2s}f'(e^{\frac{2s}{p-1}}\overline{w}_{1}))q_{1}-\frac{p+3}{p-1}q_{2}-2y\partial_{y}q_{2}\\
  \end{array}
\right),
 \end{equation}
  \begin{equation}\label{oamy}
 \hspace{-0.5cm} h(d,y,q_{1})\hspace{-0.1cm}=\hspace{-0.1cm}|\overline{w}_{1}(d,y,s)+q_{1}|^{p-1}(\overline{w}_{1}(d,y,s)+q_{1})
-\overline{w}_{1}(d,y,s)^p-p\overline{w}_{1}(d,y,s)^{p-1}q_{1},
\end{equation}
  \begin{equation}\label{e700}
 \overline{\psi}(d,y,s)=p\overline{w}_{1}^{p-1}-\frac{2(p+1)}{(p-1)^2}
 =\frac{2(p+1)}{(p-1)^2}\Big(\frac{p(1-d^2)\tilde{\phi}^{p-1}(d,y,s)}{(1+dy)^2}-1\Big),
 \end{equation}
 \end{enumerate}
with 
\begin{equation}\label{e7000}
\tilde{\phi}(d,y,s)=\frac{\phi (s-\log(\frac{1+dy}{\sqrt{1-d^2}}))}{\kappa_{0}}.
 \end{equation}
We end with the definition of $\widehat{f}$
 \begin{equation}\label{40}
\widehat{f}(\overline{w}_{1},q_{1},s)=e^{\frac{-2ps}{p-1}}\Big(f(e^{\frac{2s}{p-1}}(\overline{w}_{1}+q_{1}))-f(e^{\frac{2s}{p-1}}\overline{w}_{1})-e^{\frac{2s}{p-1}}q_{1}f'(e^{\frac{2s}{p-1}}\overline{w}_{1})\Big).
\end{equation}
\end{prop}

\begin{nb}\label{Tttaw7}
As we said above, we are going to exploite the linear operator $L_{d(s)}$ defined in $\eqref{60}$. For that reason we express the linear term $\overline{L}_{d(s)}$ differently:
\begin{equation}\label{ena20}
\overline{L}_{d(s)}\left(
  \begin{array}{ccc}
    q_{1}(y,s) \\
    q_{2}(y,s)\\
  \end{array}
\right)=L_{d(s)}\left(
  \begin{array}{ccc}
    q_{1}(y,s) \\
    q_{2}(y,s)\\
  \end{array}
\right)+\left(
  \begin{array}{ccc}
    0 \\
   \overline{V}q_{1}\\ 
  \end{array}
\right),
\end{equation}
with 
\begin{eqnarray}\label{39}
\overline{V}&=&\overline{\psi}(d,y,s)-\psi(d,y)+e^{-2s}f'(e^{\frac{2s}{p-1}}\overline{w}_{1})\nonumber \\
&=&p\kappa^{p-1}(d,y)\Big(\tilde{\phi}^{p-1}(d,y,s)-1\Big)+e^{-2s}f'(e^{\frac{2s}{p-1}}\overline{w}_{1}),
\end{eqnarray} 
where $\overline{\psi}$ and $\psi $ are defined respectively in $\eqref{38}$ and $\eqref{e700}$.
We mention also that:
\begin{equation}\label{ena41}
\left(
  \begin{array}{ccc}
   \partial_{d} \overline{w}_{1}(d,y,s) \\
   \partial_{d} \overline{w}_{2}(d,y,s)\\
  \end{array}
\right)=\left(
  \begin{array}{ccc}
   \partial_{d}\kappa(d,y) \\
    0\\
  \end{array}
\right)\\+ \left(
  \begin{array}{ccc}
   \partial_{d}\Big(\kappa(d,y)(\tilde{\phi}(d,y,s)-1 )\Big ) \\
    \partial_{d}\overline{w}_{2}(d,y,s)\\
  \end{array}
\right).
\end{equation} 
\end{nb}
This remark helps us to exploite the techniques used by Merle and Zaag \cite{fh} in a clear way.
\\{\it Proof}:
The proof of the case when $f\equiv 0$ remains valid for our perturbation with exactly the same ideas and purely technical differences. We present it to convince the reader.
\\Up to replacing $w(y,s)$ by $-w(y,s)$, we can assume that $\omega^*=1$ in $\eqref{112}$.
\\$(i)$ In $\eqref{112}$, we see that there is a parameter $d^*\in (-1,1)$ which makes the distance between the solution $(w(s^*),\partial_{s}w(s^*))$ and a particular element $(\overline{w}_{1}(d,y,s),\overline{w}_{2}(d,y,s))$ small. Now, we would like to sharpen the decomposition and find for all $s\in [s^*,\sigma^*]$ for some $\sigma^*> s^*$ a different parameter $d(s)$ close to $d^*$ which not only makes the difference between $(w(s),\partial_{s}w(s))$ and $(\overline{w}_{1}(d,y,s),\overline{w}_{2}(d,y,s))$ small but also satisfies the orthogonality condition $\eqref{101}$.
\\ From $\eqref{35}$, we see that condition $\eqref{101}$ becomes $\Phi ((w(s),\partial_{s}w(s)),d,s)=0$ where $\Phi\in \C(\H \times (-1,1)\times [s^*,+\infty ),\R)$ is defined by:
\begin{equation}\label{e7t}
 \Phi (v,d,s)=\Upsilon \Big(v-(\overline{w}_{1}(d,y,s),\overline{w}_{2}(d,y,s)),W_{0}^{d}\Big )
 \end{equation}
 with $\Upsilon$ and $W_{0}^{d}$ are given in $\eqref{50}$ and $\eqref{e4}$. The implicit function theorem allows us to conclude. Indeed, 
 \begin{itemize}
\item Note first that we have:
\begin{equation}\label{e8t}
 \Phi ((\overline{w}_{1}(d^*,y,s),\overline{w}_{2}(d^*,y,s)),d^*,s)=0.
 \end{equation}
 \item Then, we compute from $\eqref{e7t}$, the expression of $F_{0}^{d}$ written in $\eqref{61}$ and the orthogonality relation $\eqref{e6}$:
 \begin{eqnarray*}
D_{v}\Phi (v,d,s)(u)&=&\Upsilon (u,W_{0}^{d} )\,\,\,\,{\rm for}\,\,\,\,{\rm all}\,\,\,\,u\in \H,\\
\partial_{d}\Phi (v,d,s)&=&\hspace{-0.3cm}-\Upsilon \Big((\partial_{d}\overline{w}_{1}(d,y,s),\partial_{d}\overline{w}_{2}(d,y,s)),W_{0}^{d}\Big )+\Upsilon \Big(v-(\overline{w}_{1}(d,y,s),\overline{w}_{2}(d,y,s)),\partial_{d}W_{0}^{d}\Big ),\\
&=&\frac{2\kappa_{0}}{(p-1)(1-d^2)}-\Upsilon \Big((\partial_{d}\overline{w}_{1}(d,y,s)-\partial_{d}\kappa (d,y),\partial_{d}\overline{w}_{2}(d,y,s)),W_{0}^{d}\Big )\\
&&+\Upsilon \Big(v-(\overline{w}_{1}(d,y,s),\overline{w}_{2}(d,y,s)),\partial_{d}W_{0}^{d}\Big ),\\
\partial_{s}\Phi (v,d,s)&=&-\Upsilon \Big((\partial_{s}\overline{w}_{1}(d,y,s),\partial_{s}\overline{w}_{2}(d,y,s)),W_{0}^{d}\Big ).
\end{eqnarray*}
 According to the Cauchy-Shwarz inequality, the continuity of $\Upsilon$ in $\H$, the bounds $\eqref{e7}$, $\eqref{ap1}$, $\eqref{e10t}$ and $\eqref{e15t}$, we see that if 
 $$\Big|\log \Big(\frac{1+d}{1-d}\Big )-\log \Big(\frac{1+d^*}{1-d^*}\Big )\Big|+\Vert v-(\overline{w}_{1}(d^*,y,s),\overline{w}_{2}(d^*,y,s))\Vert_{\H}\leq \epsilon_{1},$$
 for some $\epsilon_{1}>0$ small enough independant of $d^*$, then we have 
 \begin{equation}\label{e11t}
\hspace{-0.1cm} \Vert D_{v}\Phi (v,d,s)\Vert +\mid\partial_{s}\Phi (v,d,s)\mid\leq C\,{\rm and}\,\frac{1}{C(1-d^2)}\leq \partial_{d}\Phi (v,d,s)\leq \frac{C}{1-d^2}.
 \end{equation}
 Now, if we introduce $\Psi \in \C(\H\times \R\times [s^*,+\infty ),\R)$ defined by
 $$\Psi (v,\theta ,s )=\Phi (v,d,s)\,\,\,{\rm where}\,\,\,d=\tanh \theta ,$$
 then, since $\theta =\frac{1}{2}\log \Big(\frac{1+d}{1-d}\Big )$ and $\tanh' (\theta)=1-\tanh^2 (\theta)$, we see from $\eqref{e8t}$ and $\eqref{e11t}$ that the implicit function theorem applies to $\Psi $ and we get the existence of $d(s)$ for all $s\in [s^*,\sigma^*)$ for some $\sigma^*\leq \infty$. 
 \\ Now, let's prove that $\sigma^* =+\infty$. We argue by contradiction and assume that $\sigma^* <\infty$, we apply the implicit function theorem around $(v,d)=((w(s_{n}),\partial_{s}w(s_{n})),d(s_{n}))$ where $s_{n}=\sigma^*-\frac{1}{n}$ and the uniform continuity of $(w(s_{n}),\partial_{s}w(s_{n}))$ from $[\sigma^*-\eta_{0},\sigma^*+\eta_{0}]$ to $\H$ for some $\eta_{0}> 0$, we see that for $n$ large enough, we can define $d(s)$ for all $s\in [s_{n},s_{n}+\epsilon_{0}]$ for some $\epsilon_{0}>0$ independant of $n$. Therefore, for $n$ large enough, $d(s)$ exist beyond $\sigma^*$, which is a contradiction. Thus, $\sigma^*=\infty$ which ends the proof of $(i)$ of Proposition $\ref{T7}$.
\end{itemize}
$(ii)$ is a direct consequence of the equation $\eqref{1}$ satisfied by $w$ put in the vectorial form:
\begin{equation}\label{mayous}
\partial_{s}w = v,\,\,\,\,\,\,\,\,\,\,\,\,\,\,\,\,\,\,\,\,\,\,\,\,\,\,\,\,\,\,\,\,\,\,\,\,\,\,\,\,\,\,\,\,\,\,\,\,\,\,\,\,\,\,\,\,\,\,\,\,\,\,\,\,\,\,\,\,\,\,\,\,\,\,\,\,\,\,\,\,\,\,\,\,\,\,\,\,\,\,\,\,\,\,\,\,\,\,\,\,\,\,\,\,\,\,\,\,\,\,\,\,\,\,\,\,\,\,\,\,\,\,\,\,\,\,\,\,\,\,\,\,\,\,\,\,\,\,\,\,\,\,\,\,\,\,\,\,\,\,\,\,\,\,\,\,\,\,
\end{equation}
\begin{equation}\label{so325}
\partial_{s}v = \pounds w-\frac{2(p+1)}{(p-1)^2}w+|w|^{p-1}w-\frac{p+3}{p-1}v
-2y \partial_{y}v+e^{\frac{-2ps}{p-1}}f(e^{\frac{2s}{p-1}}w),
\end{equation}
and the fact that $(\overline{w}_{1}(d,y,s),\overline{w}_{2}(d,y,s))$ is a solution of $\eqref{mayous}$-$\eqref{so325}$, that is $\overline{w}_{1}(d,y,s)$ is a solution of 
\begin{eqnarray}\label{mayous2}
\partial_{s}\overline{w}_{2}&=&\pounds \overline{w}_{1}-\frac{2(p+1)}{(p-1)^2}\overline{w}_{1}+\overline{w}_{1}^{p-1}\overline{w}_{1}-\frac{p+3}{p-1}\overline{w}_{2}\nonumber \\
&&-2y \partial_{y}\overline{w}_{2}+e^{\frac{-2ps}{p-1}}f(e^{\frac{2s}{p-1}}\overline{w}_{1}),
\end{eqnarray}
(see $\eqref{e505}$).
\\Indeed, since we have from $\eqref{301}$, the definition of $\pounds $ written in $\eqref{12}$ and the expression of $h(d,y,q_{1})$ written in $\eqref{oamy}$ 
\begin{eqnarray*}
 w(y,s)&=&\overline{w}_{1}(d,y,s)+q_{1}(y,s),\\
 \partial_{s}w(y,s)&=& \partial_{s}\overline{w}_{1}(d,y,s)+d' \partial_{d}\overline{w}_{1}(d,y,s)+\partial_{s}q_{1}(y,s),\\
 \pounds w(y,s) &=&\pounds\overline{w}_{1}(d,y,s)+\pounds q_{1}(y,s),\\
 |w|^{p-1}w &=& h(d,y,q_{1}(y,s))+\overline{w}_{1}^p(d,y,s)+p\overline{w}_{1}^{p-1}(d,y,s)q_{1}(y,s),\\ 
 -2y\partial^2_{y,s}w(y,s)&=&-2y\partial_{y}\overline{w}_{2}(d,y,s)-2y\partial_{y}q_{2}(y,s),\\
 e^{\frac{-2ps}{p-1}}f(e^{\frac{2s}{p-1}}w(y,s))&=&e^{\frac{-2ps}{p-1}}f(e^{\frac{2s}{p-1}}(\overline{w}_{1}(d,y,s)+q_{1}(y,s))).
 \end{eqnarray*}
 We see that equation $\eqref{21}$ follows immediately from $\eqref{mayous}$-$\eqref{mayous2}$. This conclude the proof of Proposition $\ref{T7}$.
\Box
\subsection{Projection on the eigenspaces of the operator $L_{d}$}
For the proof of the main Theorem $\ref{T4}$, we need to prove in some sense dispersive estimates on $q_{-}=\pi_{-}^{d}(q)$ when $q$ is a solution to $\eqref{21}$. In order to achieve this, we need to manipulate a function of $q_{-}$ (equivalent to the norm $\Vert q_{-}\Vert_{\H}=\Upsilon (q_{-},q_{-})^{\frac{1}{2}}$ in $\H_{-}^{d}$) which will capture the dispersive character of the equation $\eqref{21}$. Such a quantity will be 
\begin{eqnarray}\label{500}
\varphi_{d}(q,r)&=&\int_{-1}^{1}(-\psi(d,y)q_{1}r_{1}+\partial_{y}q_{1}\partial_{y}r_{1}(1-y^2)+q_{2}r_{2})\rho dy\nonumber\\
&=&\int_{-1}^{1}(-q_{1}(\pounds r_{1}+\psi (d,y)r_{1})+q_{2}r_{2})\rho dy,
\end{eqnarray}
with  $\psi (d,y)$ is defined in $\eqref{38}$.
\begin{nb}
 We remark that $\varphi_{d}(q,r)$ is the same bilinear form introduced in Merle and Zaag \cite{fh} in the case when $f\equiv 0$.
\end{nb}
It is worth mentioning that this bilinear form $\varphi_{d}(q,r)$ is in fact the second variation of $E_{0}(w(s))$ defined in $\eqref{cp}$ around $\kappa (d,y)$ defined in $\eqref{400}$, which is a stationary solution of $\eqref{1}$ when $(f\equiv 0)$ and can be seen as the energy norm in $\H_{-}^{d}$ (space where it will be definite positive). We recall briefly from Merle and Zaag \cite{fh} in the following the continuity of $\varphi_{d}(q,r)$:
\\For all $(q,r)\in \H^2$ and $s\in [s^*,+\infty )$, we have:
\begin{equation}\label{oa}
 |\varphi_{d}(q,r)|\leq C\Vert q\Vert_{\H}\Vert r\Vert_{\H}.
 \end{equation}
As a matter of fact, it is reasonable to recall the following a priori estimate:
\begin{equation}\label{200}
\Vert q \Vert_{\H}\leq \epsilon,
\end{equation}
for some $s \geq s^*$ and some $\epsilon > 0.$ In addition to that, we would like to expand $q$ from $\eqref{101}$ according to the linear operator $L_{d}$:
 \begin{equation}\label{202}
q(y,s)=\alpha_{1}(s) F_{1}^{d(s)}(y)+q_{-}(y,s),
\end{equation}
where 
\begin{equation}\label{203}
\alpha_{1}(s) =\pi_{1}^{d(s)}(q),\,\,\,\alpha_{0}(s) =\pi_{0}^{d(s)}(q)=0,\,\,\,\alpha_{-}(s) =\sqrt{\varphi_{d}(q_{-},q_{-})}
\end{equation}
and 
\begin{equation}\label{204}
q_{-}=\left(
  \begin{array}{ccc}
   q_{-,1} \\
    q_{-,2}\\
  \end{array}
\right)=\pi_{-}^{d}(q)=\pi_{-}^{d}\left(
  \begin{array}{ccc}
   q_{1} \\
    q_{2}\\
  \end{array}
\right).
\end{equation}
Beside that, from Proposition $4.7$ Page 90 in Merle and Zaag \cite{fh} and $\eqref{202}$, we see that for all $s\geq  s^*$,
\begin{enumerate}[{\rm(i)}]
\item {\bf (Equivalence of norms in $\H_{-}^{d}$)} For all $q_{-}\in \H_{-}^{d}$,
\begin{equation}\label{raw}
\frac{1}{C}\alpha_{-}(s) \leq \Vert q_{-}(s)\Vert_{\H}\leq C\alpha_{-}(s).
\end{equation}
\item {\bf (Equivalence of norms in $\H$)} For all $q\in \H$,
\begin{equation}\label{raw1}
\frac{1}{C}(|\alpha_{1}(s)|+\alpha_{-}(s)) \leq  \Vert q(s)\Vert_{\H}\leq  C(|\alpha_{1}(s)|+\alpha_{-}(s)).
 \end{equation}
\end{enumerate}
for some $C> 0$. We present now the very heart of our argument. Here, we derive from $\eqref{21}$ a differential inequalities satisfied by $\alpha_{1}(s)$, $\alpha_{-}(s)$ and $d(s)$:
\begin{prop}\label{T8}
There exists $\epsilon_{2} > 0$ such that if $w$ is a solution to equation $\eqref{1}$ satisfying $\eqref{101}$ and $\eqref{200}$ at some time $s$ for some $\epsilon <\epsilon_{2} $, where $q$ is defined in $\eqref{301}$, then 
\begin{enumerate}[{\rm(i)}]
\item {\bf (Control of the modulation parameter)} For all $s\geq s^*$, we have
\begin{equation}\label{305}
\frac{|d'|}{1-d^2}\leq C(\alpha_{1}^2+\alpha_{-}^{2})+\frac{C}{s^{a}}(\alpha_{1}^2+\alpha_{-}^{2})^{\frac{1}{2}}.
\end{equation}
\item {\bf (Projection of equation $\eqref{21}$ on the mode $\lambda =1$ and the negative part)}
For all $s\geq s^*$, we have
\begin{equation}\label{306}
|\alpha'_{1}-\alpha_{1}|\leq C(\alpha_{1}^2+\alpha_{-}^{2})+\frac{C}{s^{a}}(\alpha_{1}^2+\alpha_{-}^{2})^{\frac{1}{2}},
\end{equation}
\begin{equation}\label{310}
\Big(R_{-}+\frac{1}{2}\alpha_{-}^2\Big)'\leq -\frac{4}{p-1}\int_{-1}^{1}q_{-,2}^2\frac{\rho}{1-y^2}dy+C(\alpha_{1}^2+\alpha_{-}^{2})^{\frac{1+\overline{p}}{2}}+\frac{C}{s^{a}}(\alpha_{1}^2+\alpha_{-}^{2})
\end{equation}
for some $R_{-}(s)$ satisfying 
\begin{equation}\label{307}
|R_{-}(s)|\leq C(\alpha_{1}^2+\alpha_{-}^{2})^{\frac{1+\overline{p}}{2}},\,\,\,{\rm where}\,\,\,\overline{p}=\min(p,2)>1.
\end{equation}
\item {\bf (Additional relation)}
For all $s\geq s^*$, we have
\begin{equation}\label{308}
\frac{d}{ds}\int_{-1}^{1}q_{1}q_{2}\rho dy\leq -\frac{1}{5}\alpha_{-}^2+ C\int_{-1}^{1}q_{-,2}^2\frac{\rho}{1-y^2}dy+C\alpha_{1}^2.
\end{equation}
\item {\bf (Energy barrier)} If moreover $\eqref{111}$ holds, then 
\begin{equation}\label{309}
|\alpha_{1}(s)|^2\leq c_{2}\alpha_{-}(s)^2+\frac{C}{s^{\frac{a+1}{2}}}.
\end{equation}
\end{enumerate}
\end{prop}

\begin{nb}
The proof of Proposition $\ref{T8}$ is the major step of this paper. Indeed this proposition allows us to derive in the next section the polynomial and the exponential decay. From the polynomial decay, we show that for $d(s)$ introduced in $\eqref{101}$, we have $\Vert q(s)\Vert_{\H}\longrightarrow 0$ as $s\longrightarrow +\infty$. Thanks to this special information combined with some refinement of some results obtained in the polynomial decay, we get the exponential decay. Here lays a major difference between our approach and the case when $f\equiv 0$.
\end{nb}

{\it{\bf Proof of Proposition $\ref{T8}$}}: Before going into the proof of Proposition $\ref{T8}$, let us first derive  some nonlinear estimates which will be useful for getting Proposition $\ref{T8}$. Following the method used in Merle and Zaag \cite{fh} the proof of  Proposition $\ref{T8}$ requests the following nonlinear estimates.

\begin{lem}\label{T9} {\bf (Nonlinear estimates)} For all $y\in (-1,1)$, we have 
 \begin{equation}\label{31}
|h(d,y,q_{1})|\leq C\delta_{p\geq 2}|\kappa(d(s),y))|^{p-2}|q_{1}(y,s)|^2+C|q_{1}(y,s)|^p,
\end{equation}
\begin{equation}\label{32}
|H(d,y,q_{1})|\leq C\delta_{p\geq 2}|\kappa(d(s),y))|^{p-2}|q_{1}(y,s)|^3+C|q_{1}(y,s)|^{p+1},
\end{equation}
where $\delta_{p\geq 2}$ is $0$ if $1<p<2$ and $1$ otherwise, the function $h$ is defined in $\eqref{oamy}$ and 
\begin{equation}\label{33}
H(d,y,q_{1})=\int_{0}^{q_{1}}h(d,y,q')dq'=\frac{|\overline{w}_{1}+q_{1}|^{p+1}}{p+1}-\frac{\overline{w}_{1}^{p+1}}{p+1}-\overline{w}_{1}^{p}q_{1}-\frac{p}{2}\overline{w}_{1}^{p-1}q_{1}^2.
\end{equation}
\end{lem}
{\it Proof}: The proof of $\eqref{31}$ and $\eqref{32}$ is exactly the same as the one written in Claim 5.3 page $104$ in Merle and Zaag \cite{fh}: just replace $\kappa(d,y)$ by $\overline{w}_{1}(d,y,s)$ and use the fact that $\frac{\overline{w}_{1}(d,y,s)}{\kappa(d,y)}$ is a bounded function from $\eqref{fer}$.

\Box

Let us now introduce the following lemma, where we give some nonlinear estimates related to our perturbation $f$ defined in $\eqref{fy}$ and the new solution $\overline{w}_{1}$ defined in $\eqref{87}$ which will play a central role in our analysis.
\begin{lem}\label{T10}
{\bf (Nonlinear estimates related to $f$ and $\overline{w}_{1}$)} For all $y\in (-1,1)$, we have 
 \begin{equation}\label{100000}
\vert \overline{V}\vert\leq Ce^{-s} +\frac{C}{s^{a}} \kappa^{p-1} (d,y),
\end{equation}
\begin{equation}\label{e11}
|\widehat{f}|\leq C\delta_{p\geq 2}|\kappa(d(s),y)|^{p-2}|q_{1}(y,s)|^2+C|q_{1}(y,s)|^p ,
\end{equation}
\begin{equation}\label{ena50}
|\widehat{F}|\leq C\delta_{p\geq 2}|\kappa(d(s),y)|^{p-2}|q_{1}(y,s)|^3+C|q_{1}(y,s)|^{p+1},
\end{equation}
 where $\delta_{p\geq 2}$ is $0$ if $1<p<2$ and $1$ otherwise and
$$\widehat{F}(\overline{w}_{1},q_{1},s)=\int_{0}^{q_{1}}\widehat{f}(\overline{w}_{1},q',s)dq',$$
with $\widehat{f}(\overline{w}_{1},q_{1},s)$ and $\overline{V}$ defined in $\eqref{40}$.

\end{lem}

{\it Proof}: We start by the proof of $\eqref{100000}$. Inspired by the proof of Lemma 2.1 page 1121 in our paper \cite{X}, the following holds 
\begin{equation}\label{ena3}
e^{-2s}|f'(e^{\frac{2s}{p-1}}\overline{w}_{1})|\leq Ce^{-s}+\frac{C}{s^{a}} |\overline{w}_{1}|^{p-1}.
\end{equation}
From the expression $\eqref{39}$ of $\overline{V}$, the fact that $\frac{\overline{w}_{1}(d,y,s)}{\kappa(d,y)}$ is a bounded function, the first point of Claim $\ref{T14}$ and $\eqref{ena3}$, we can write
\begin{equation}\label{ena4}
\vert \overline{V}\vert\leq Ce^{-s} +\frac{C}{s^{a}} \kappa^{p-1} (d,y).
\end{equation}
We deal now with the proof of $\eqref{e11}$ and $\eqref{ena50}$. From the expression $\eqref{fy}$ of $f$ and $\eqref{oamy}$ of $h$, it holds that
\begin{equation}\label{app100}
|\widehat{f}(\overline{w}_{1},q_{1},s)| \leq C|h(d,y,q_{1})|+|\widehat{f}_{1}(\overline{w}_{1},q_{1},s)+\widehat{f}_{2}(\overline{w}_{1},q_{1},s)|,
\end{equation}
where 
$$\widehat{f}_{1}(\overline{w}_{1},q_{1},s) = \Big(\underbrace{|\overline{w}_{1}+q_{1}|^{p-1}(\overline{w}_{1}+q_{1})-|\overline{w}_{1}|^{p-1}\overline{w}_{1} }_{\widehat{f}_{1,1}(\overline{w}_{1},q_{1})}\Big )\Big[\underbrace{\frac{1}{\log^{a} \Big(2+e^{\frac{4s}{p-1}}(\overline{w}_{1}+q_{1})^2\Big)}-\frac{1}{\log^{a} \Big(2+e^{\frac{4s}{p-1}}\overline{w}_{1}^2\Big)}}_{\widehat{f}_{1,2}(\overline{w}_{1},q_{1},s)}\Big ],$$
\begin{eqnarray}\label{so350}
\widehat{f}_{2}(\overline{w}_{1},q_{1},s)&=&|\overline{w}_{1}|^{p-1}\overline{w}_{1}\Big[\frac{1}{\log^{a} \Big(2+e^{\frac{4s}{p-1}}(\overline{w}_{1}+q_{1})^2\Big)}
-\frac{1}{\log^{a} \Big(2+e^{\frac{4s}{p-1}}\overline{w}_{1}^2\Big)}\nonumber\\
&&+\frac{2aq_{1}}{\log^{a+1} \Big(2+e^{\frac{4s}{p-1}}\overline{w}_{1}^2\Big)}\frac{e^{\frac{4s}{p-1}}\overline{w}_{1}}{\Big(2+e^{\frac{4s}{p-1}}\overline{w}_{1}^2\Big )}\Big ].
\end{eqnarray}
We apply the mean value theorem and the fact that $|a+b|^{p}\leq C(|a|^{p}+|b|^{p})$ for any real numbers $a$ and $b$, to write  
\begin{equation}\label{so38}
|\widehat{f}_{1,1}(\overline{w}_{1},q_{1})|\leq C|q_{1}|(|\overline{w}_{1}|^{p-1}+|q_{1}|^{p-1}).
\end{equation}
Remarking that there exists $C$ such that $|\widehat{f}_{1,2}(\overline{w}_{1},q_{1},s)|\leq C$, we can write from $\eqref{so38}$ that if $\frac{|q_{1}|}{|\overline{w}_{1}|}\geq \frac{1}{2}$, the following  
\begin{equation}\label{so39}
|\widehat{f}_{1}(\overline{w}_{1},q_{1},s)|\leq C|q_{1}|^p.
\end{equation}
In the case when $\frac{|q_{1}|}{|\overline{w}_{1}|}< \frac{1}{2}$, we apply the mean value theorem to derive the following
\begin{equation}\label{so310}
|\widehat{f}_{1,2}(\overline{w}_{1},q_{1},s)|\leq \frac{Ce^{\frac{4s}{p-1}}|\overline{w}_{1}||q_{1}|}{2+e^{\frac{4s}{p-1}}(\overline{w}_{1}+\theta_{1} q_{1})^2}.
\end{equation}
In this case (when $\frac{|q_{1}|}{|\overline{w}_{1}|}< \frac{1}{2}$), we separate the cases $p\geq 2$ and $1<p<2$. In the case when $p\geq 2$, by virtue of $\eqref{so38}$ and $\eqref{so310}$ entails
\begin{equation}\label{so311}
|\widehat{f}_{1}(\overline{w}_{1},q_{1},s)|\leq C|\overline{w}_{1}|^{p-2}|q_{1}|^2.
\end{equation} 
For the case when $1<p<2$, we remark that $|q_{1}|^{p-2}\geq \frac{|\overline{w}_{1}|^{p-2}}{2^{p-2}}$, which implie from inequalities $\eqref{so38}$ and $\eqref{so310}$ that
\begin{equation}\label{so360}
|\widehat{f}_{1}(\overline{w}_{1},q_{1},s)|\leq C|q_{1}|^p.
\end{equation} 
We combine now $\eqref{so39}$, $\eqref{so311}$ and $\eqref{so360}$ to deduce that for all $q_{1}$ and $\overline{w}_{1}$
\begin{equation}\label{so312}
|\widehat{f}_{1}(\overline{w}_{1},q_{1},s)|\leq C\delta_{p\geq 2}|\overline{w}_{1}|^{p-2}|q_{1}|^2+C|q_{1}|^p,
\end{equation} 
where $\delta_{p\geq 2}$ is $0$ if $1<p<2$ and $1$ otherwise. We treat now the term $\widehat{f}_{2}(\overline{w}_{1},q_{1},s)$. We write
\begin{equation}\label{so313}
|\widehat{f}_{2}(\overline{w}_{1},q_{1},s)|\leq C|\overline{w}_{1}|^{p}+C|\overline{w}_{1}|^{p-1}|q_{1}|.
\end{equation} 
Directly, we can write from $\eqref{so313}$ if $\frac{|q_{1}|}{|\overline{w}_{1}|}\geq \frac{1}{2}$ 
\begin{equation}\label{so314}
|\widehat{f}_{2}(\overline{w}_{1},q_{1},s)|\leq C|q_{1}|^{p}.
\end{equation} 
Now, if $\frac{|q_{1}|}{|\overline{w}_{1}|}< \frac{1}{2}$, similarly to $\eqref{so311}$ and $\eqref{so360}$, we apply again the mean value theorem to write the following when $p\geq 2$
\begin{equation}\label{so315}
|\widehat{f}_{2}(\overline{w}_{1},q_{1},s)|\leq C|\overline{w}_{1}|^{p-2}|q_{1}|^2.
\end{equation}
When $1<p<2$, we write 
\begin{equation}\label{so361}
|\widehat{f}_{2}(\overline{w}_{1},q_{1},s)|\leq C|q_{1}|^p.
\end{equation}
We combine $\eqref{so314}$, $\eqref{so315}$ and $\eqref{so361}$ to deduce that for all $q_{1}$ and $\overline{w}_{1}$
\begin{equation}\label{so316}
|\widehat{f}_{2}(\overline{w}_{1},q_{1},s)|\leq C\delta_{p\geq 2}|\overline{w}_{1}|^{p-2}|q_{1}|^2+C|q_{1}|^p.
\end{equation}
Thanks to $\eqref{31}$ of Lemma $\ref{T9}$ combined with $\eqref{app100}$, $\eqref{so312}$ and $\eqref{so316}$, we conclude the proof of $\eqref{e11}$. Then, we exploite the expression of $\widehat{F}$ to conclude the proof $\eqref{ena50}$. This conclude the proof of Lemma $\ref{T10}$.

\Box

We give now the strategy of the proof of $(i)$-$(ii)$ of Proposition $\ref{T8}$. We proceed in two steps:
 \\-In Step 1, we project equation $\eqref{21}$ with the projector $\pi_{\lambda}^{d}$ defined in $\eqref{35}$ for $\lambda =0$ and $\lambda =1$ and derive the smallness condition on $d'$ in $\eqref{305}$ and the equation satisfied by $\alpha_{1}$ in $\eqref{306}$.
 \\-In Step 2, we write an equation satisfied by $(q_{-,1},q_{-,2})$ which is the difficult part in this non self-adjoint framework. We claim that inequality $\eqref{310}$ follows from the existence of the Lyapunov functional $\eqref{e90}$ for equation $\eqref{1}$. Here, the Lyapunov functional structure will be revealed by the quadratic form $\varphi_{d}$ $\eqref{500}$.
 \subparagraph{Step1: Projection of equation $\eqref{21}$ on the modes $\lambda =0$ and $\lambda =1$.}Projecting equation $\eqref{21}$ with the projector $\pi_{\lambda}^{d}$ defined in $\eqref{35}$ for $\lambda =0$ and $\lambda =1$, we write
\begin{eqnarray} \label{63}
\pi_{\lambda}^{d}(\partial_{s}q)&=&\pi_{\lambda}^{d}(L_{d}q)+\pi_{\lambda}^{d}\left(
  \begin{array}{ccc}
    0 \\
    h(d,y,q_{1})\\
  \end{array}
\right)-d'\pi_{\lambda}^{d}\left(
  \begin{array}{ccc}
   \partial_{d}\kappa(d,y) \\
    0\\
  \end{array}
\right) \nonumber\\
&&+ \Sigma_{\lambda ,d}^{1}(s)-d'\Sigma_{\lambda ,d}^{2}(s)+\Sigma_{\lambda ,d}^{3}(s),
\end{eqnarray}
where 
\begin{itemize}
\item $\Sigma_{\lambda ,d}^{1}(s)=\pi_{\lambda}^{d}\left(
  \begin{array}{ccc}
    0 \\
   \overline{V}q_{1}\\
  \end{array}
\right),$
\item $\Sigma_{\lambda ,d}^{2}(s)=\pi_{\lambda}^{d}\left(
  \begin{array}{ccc}
  \partial_{d}\Big( \kappa(d,y)(\tilde{\phi}(d,y,s)-1)\Big ) \\
    \partial_{d} \overline{w}_{2}(d,y,s)\\
  \end{array}
\right),$
\item $\Sigma_{\lambda ,d}^{3}(s)=\pi_{\lambda}^{d}\left(
  \begin{array}{ccc}
    0 \\
    \widehat{f}(\overline{w}_{1},q_{1},s)\\
  \end{array}
\right),$
\end{itemize}
with, $L_{d}$, $\kappa(d,y)$, $\overline{V}$ and $\widehat{f}(\overline{w}_{1},q_{1},s)$ defined respectively in $\eqref{60}$, $\eqref{400}$, $\eqref{39}$ and $\eqref{40}$. 
According to the first step of the proof of Proposition 5.2 page $105$ in Merle and Zaag \cite{fh}, we write directly when $\lambda =0$,
\begin{eqnarray}\label{arb1}
\frac{2\kappa_{0}\mid d'\mid}{(p-1)(1-d^2)} &\leq &\frac{C|d'|(|\alpha_{1}(s)|+\alpha_{-}(s))}{1-d^2}+C\Vert q\Vert^2_{\H} \nonumber\\
&&+\mid\Sigma_{0 ,d}^{1}(s)\mid +\mid d'\Sigma_{0 ,d}^{2}(s)\mid +\mid\Sigma_{0 ,d}^{3}(s)\mid
\end{eqnarray}
and when $\lambda =1$,
\begin{eqnarray}\label{arb2}
\mid \alpha'_{1}(s)-\alpha_{1}(s)\mid &\leq & \frac{C|d'|(|\alpha_{1}(s)|+\alpha_{-}(s))}{1-d^2}+C\Vert q\Vert^2_{\H}\nonumber\\
&&+\mid\Sigma_{1 ,d}^{1}(s)\mid +\mid d'\Sigma_{1 ,d}^{2}(s)\mid + \mid \Sigma_{1 ,d}^{3}(s)\mid .
\end{eqnarray}
Our focal interst now is to treat the new terms $\Sigma_{\lambda ,d}^{1}(s)$, $\Sigma_{\lambda ,d}^{2}(s)$ and $\Sigma_{\lambda ,d}^{3}(s)$ with $\lambda \in \lbrace 0,\,1\rbrace$. With Claim $\ref{T14}$, we are in position to give an estimation to $\Sigma_{\lambda ,d}^{1}(s)$.
We use the definition $\eqref{35}$ of $\pi_{\lambda}^{d}$, the expression of $\Upsilon$ given in $\eqref{50}$ and the inequality $\eqref{100000}$, we see that 
 \begin{equation} \label{103}
\Big|\Sigma_{\lambda ,d}^{1}(s)\Big|\\ \leq C e^{-s}\int_{-1}^{1}|q_{1}|\,|W_{\lambda ,2}^{d}|\rho dy +\frac{C}{s^{a}}\int_{-1}^{1}|\kappa^{p-1}(d,y)|\,|q_{1}|\,|W_{\lambda ,2}^{d}|\rho dy .
\end{equation} 
We are going now to estimate one by one the terms of the right-hand side of inequality $\eqref{103}$.
 According to $\eqref{47}$, the H\"{o}lder inequality and the Hardy-Sobolev's inequality in Lemma $\ref{T21}$, we can see the following estimate:
\begin{equation} \label{102} 
\int_{-1}^{1}|\kappa^{p-1}(d,y)|\,|q_{1}||W_{\lambda ,2}^{d}|\rho dy\leq C\hspace{-0.2cm}\int_{-1}^{1}|\kappa^{p}(d,y)|\,|q_{1}|\rho dy 
\leq C\Vert \kappa \Vert^{p}_{\H_{0}}\Vert q \Vert_{\H}\leq C\Vert q \Vert_{\H}.
\end{equation}
Again via $\eqref{47}$, the H\"{o}lder inequality and the Hardy-Sobolev's inequality in Lemma $\ref{T21}$ to obtain
\begin{equation} \label{1003}
\int_{-1}^{1}|q_{1}|\,|W_{\lambda ,2}^{d}|\rho dy \leq   C\Vert q \Vert_{\H}.
\end{equation} 
Collecting $\eqref{103}$, $\eqref{102}$ and $\eqref{1003}$ together, we deduce
\begin{equation} \label{10003}
\Big|\Sigma_{\lambda ,d}^{1}(s)\Big| \leq \frac{C}{s^{a}}\Vert q \Vert_{\H}.
\end{equation} 
We would like now to estimate $\Sigma_{\lambda ,d}^{2}(s)$. From the definition $\eqref{35}$ of $\pi_{\lambda}^{d}$, the expression of $\Upsilon$ given in $\eqref{50}$ and inequality $\eqref{47}$, it holds that:
\begin{eqnarray}\label{e22}
\Big|\Sigma_{\lambda ,d}^{2}(s)\Big|  &\leq &\int_{-1}^{1}|\partial_{d}\Big( \kappa(d,y)(\tilde{\phi}(d,y,s)-1 )\Big)|\,|-\pounds W_{\lambda ,1}^{d}+W_{\lambda ,1}^{d}|\rho dy \nonumber\\
&&+ C\int_{-1}^{1}|\kappa(d,y)|\, |\partial_{d} \overline{w}_{2}(d,y,s)|\rho dy.
\end{eqnarray}
By using equation $\eqref{e5}$ satisfied by $W_{\lambda ,1}^{d}$ and $\eqref{47}$, we can see that 
\begin{equation}\label{e50}
|-\pounds W_{\lambda ,1}^{d}+W_{\lambda ,1}^{d}| \leq C\frac{\kappa(d,y)}{1-y^2}.
\end{equation}
Using Claim $\ref{T14}$, inequalities $\eqref{e54}$, $\eqref{zet}$ and $(ii)$ of Claim $\ref{T21}$ to deduce that 
\begin{equation}\label{e53}
\int_{-1}^{1}|\partial_{d}\Big( \kappa(d,y)(\tilde{\phi}(d,y,s)-1 )\Big)|\,|-\pounds W_{\lambda ,1}^{d}+W_{\lambda ,1}^{d}|\rho dy  \leq \frac{C}{s^{a}(1-d^2)}.
\end{equation}
We get directly from inequality $\eqref{e51}$ and $(ii)$ of Claim $\ref{T21}$ the following
\begin{equation}\label{e57}
\int_{-1}^{1}\kappa(d,y) |\partial_{d} \overline{w}_{2}(d,y,s)|\rho dy   \leq \frac{C}{s^{a}(1-d^2)}.
 \end{equation}
Plugging $\eqref{e53}$ and $\eqref{e57}$ into $\eqref{e22}$, we observe that
\begin{equation}\label{e60}
\Big|\Sigma_{\lambda ,d}^{2}(s) \Big|\\   \leq \frac{C}{s^{a}(1-d^2)}.
 \end{equation}
We treat now the term $\Sigma_{\lambda ,d}^{3}(s)$ which is provided from the perturbation $f$ defined in $\eqref{fy}$. We use the definition $\eqref{35}$ of $\pi_{\lambda}^{d}$, the expression of $\Upsilon$ given in $\eqref{50}$, inequality $\eqref{e11}$ of Lemma $\ref{T10}$ and $\eqref{47}$, we write
\begin{equation}\label{e15}
\Big|\Sigma_{\lambda ,d}^{3}(s)\Big|\leq C\int_{-1}^{1}|\kappa(d,y)|^{p-1}|q_{1}|^2\rho dy+C\delta_{p\geq 2}\int_{-1}^{1}|\kappa(d,y)|\,|q_{1}|^p\rho dy ,
\end{equation}
where $\delta_{p\geq 2}$ is $0$ if $1<p<2$ and $1$ otherwise. From the H\"{o}lder inequality and the Hardy-Sobolev's inequality in Lemma $\ref{T21}$, we write 
\begin{equation} \label{ena200}
\int_{-1}^{1}|\kappa (d,y)|^{p-1} |q_{1}|^2\rho dy \leq \Vert \kappa \Vert^{p-1}_{\H_{0}}\Vert q \Vert^2_{\H}\leq C\Vert q \Vert^2_{\H},
\end{equation}
for the same reason, we write also the following 
\begin{equation}\label{jar}
\int_{-1}^{1}|\kappa(d,y)|\,|q_{1}|^p\rho dy \leq \Vert \kappa \Vert_{\H_{0}}\Vert q \Vert^{p}_{\H}\leq C\Vert q \Vert^{p}_{\H}.
\end{equation}
Observing inequalities $\eqref{e15}$, $\eqref{ena200}$, $\eqref{jar}$ along with the a priori estimate $\eqref{200}$ yield
\begin{equation}\label{jar1}
\Big|\Sigma_{\lambda ,d}^{3}(s)\Big|\leq C\Vert q \Vert^{2}_{\H} .
\end{equation}
According to $\eqref{63}$, $\eqref{arb1}$, $\eqref{arb2}$, $\eqref{10003}$, $\eqref{e60}$ and $\eqref{jar1}$ to obtain when $\lambda =0$
\begin{equation}\label{arb10}
\frac{2\kappa_{0}\mid d'\mid}{(p-1)(1-d^2)} \leq \frac{C|d'|}{1-d^2}(|\alpha_{1}(s)|+\alpha_{-}(s))+C\Vert q\Vert^2_{\H}+\frac{C}{s^{a}}\Vert q\Vert_{\H}+\frac{C|d'|}{s^{a}(1-d^2)},
\end{equation}
and when $\lambda =1$
\begin{equation}\label{arb11}
\mid \alpha'_{1}(s)-\alpha_{1}(s)\mid\leq \frac{C|d'|}{1-d^2}(|\alpha_{1}(s)|+\alpha_{-}(s))+C\Vert q\Vert^2_{\H}+\frac{C}{s^{a}}\Vert q\Vert_{\H}+\frac{C|d'|}{s^{a}(1-d^2)}.
\end{equation}
Using the smallness condition $\eqref{200}$, the equivalence of the norms in $\eqref{raw1}$, the fact that 
\\$C\frac{|d'|}{1-d^2}\Vert q\Vert_{\H}\leq \frac{\kappa_{0}}{4(p-1)}\frac{|d'|}{1-d^2}$ and the fact that for some $s^*$ large enough and for all $s\geq s^*$, we have $C\frac{|d'|}{s^{a}(1-d^2)}\leq \frac{\kappa_{0}}{4(p-1)}\frac{|d'|}{1-d^2}$, we get $\eqref{305}$ and $\eqref{306}$ for $\epsilon$ small enough.
\subparagraph{Step 2: Differential inequality on $\alpha_{-}$.}
In the following, we project equation $\eqref{21}$ on the negative modes, which gives a partial differential inequality satisfied by $q_{-}$.
In order to simplify the presentation, let us introduce the following terms which will be useful in many steps of this part.
\begin{equation}\label{sta}
\Sigma_{- ,d}^{1}(s) = \pi_{-}^{d}\left(
  \begin{array}{ccc}
    0 \\
   \overline{V}q_{1}\\
  \end{array}
\right),\,\,\,\,\,\,\,\,\,\,\,\,\,\,\,\,\,\,\,\,\,\,\,\,\,\,\,\,\,\,\,\,\,\,\,\,\,\,\,\,\,\,\,\,\,\,\,\,\,\,\,\,\,\,\,\,\,\,\,\,\,\,\,\,
\end{equation}
\begin{equation}\label{sta1}
\Sigma_{- ,d}^{2}(s) = d'\pi_{-}^{d}\left(
  \begin{array}{ccc}
  \partial_{d}\Big( \kappa(d,y)(\tilde{\phi}(d,y,s)-1)\Big ) \\
    \partial_{d} \overline{w}_{2}(d,y,s)\\
  \end{array}
\right),
\end{equation}
\begin{equation}\label{sta2}
\Sigma_{- ,d}^{3}(s) = \pi_{-}^{d}\left(
  \begin{array}{ccc}
    0 \\
    h(d,y,q_{1})+\widehat{f}(\overline{w}_{1},q_{1},s)\\
  \end{array}
\right),
\end{equation}

\begin{eqnarray*}
\Sigma_{- ,d}^{3,1}(s)&=&\pi_{-}^{d}\left(
  \begin{array}{ccc}
    0 \\
    h(d,y,q_{1})\\
  \end{array}
\right),\\
\Sigma_{- ,d}^{3,2}(s)&=&\pi_{-}^{d}\left(
  \begin{array}{ccc}
    0 \\
    \widehat{f}(\overline{w}_{1},q_{1},s)\\
  \end{array}
\right).
\end{eqnarray*}
We now claim the following:
\begin{cl}\label{T11}{\bf (Preliminary estimates)}
There exists $\epsilon_{3}> 0$ such that if $\epsilon < \epsilon_{3}$ in the hypotheses of Proposition $\ref{T8}$, then 
\begin{equation}\label{e200}
\Vert \partial_{s}q_{-}-L_{d}q_{-}-\Sigma_{- ,d}^{1}(s)-\Sigma_{- ,d}^{2}(s)
-\Sigma_{- ,d}^{3}(s) \Vert_{\H} \leq  C(\alpha_{1}^2+\alpha_{-}^{2})^{\frac{3}{2}},
\end{equation}
\begin{equation}\label{2030}
\Big|\varphi_{d}\Big(q_{-}, \Sigma_{- ,d}^{1}(s)\Big)\Big|\leq \frac{C}{s^{a}}(\alpha_{1}^2+\alpha_{-}^{2}),
\end{equation} 
\begin{equation}\label{rama}
\Big|\varphi_{d}\Big(q_{-}, \Sigma_{- ,d}^{2}(s)\Big)\Big|\leq C(\alpha_{1}^2+\alpha_{-}^{2})^{\frac{3}{2}}+\frac{C}{s^{a}}(\alpha_{1}^2+\alpha_{-}^{2}),
\end{equation}
\begin{equation}\label{e201}
\Big|\varphi_{d}\Big(q_{-}, \Sigma_{- ,d}^{3,1}(s)\Big )-\int_{-1}^{1}q_{2}h(d,y,q_{1})\rho dy
  \Big|\leq C(\alpha_{1}^2+\alpha_{-}^{2})^{\frac{3}{2}},
\end{equation}
\begin{equation}\label{hs60}
\Big|\varphi_{d}\Big(q_{-},\Sigma_{- ,d}^{3,2}(s)\Big)-\int_{-1}^{1}q_{2}\widehat{f}(\overline{w}_{1},q_{1},s)\rho dy \Big|\leq C(\alpha_{1}^2+\alpha_{-}^{2})^{\frac{3}{2}}
,
\end{equation}
\begin{equation}\label{e202}
\Big|\int_{-1}^{1}q_{2}h(d,y,q_{1})\rho dy
  -\frac{d}{ds}\int_{-1}^{1}H(d,y,q_{1})\rho dy\Big| \leq C(\alpha_{1}^2+\alpha_{-}^{2})^{2}+\frac{C}{s^{a}}(\alpha_{1}^2+\alpha_{-}^{2}),
\end{equation}
\begin{equation}\label{hs45}
 \Big |\int_{-1}^{1}q_{2}\widehat{f}(\overline{w}_{1},q_{1},s)\rho dy -\frac{d}{ds}\int_{-1}^{1}\widehat{F}(\overline{w}_{1},q_{1},s)\rho dy \Big |\leq  C (\alpha_{1}^2+\alpha_{-}^{2})^{\frac{1+\overline{p}}{2}}+\frac{C}{s^{a}}(\alpha_{1}^2+\alpha_{-}^{2}),
 \end{equation}
where $h(d,y,q_{1})$, $H(d,y,q_{1})$ and $\widehat{F}(\overline{w}_{1},q_{1},s)$ are defined respectively in $\eqref{oamy}$, $\eqref{33}$ and $\eqref{hs44}$.
\end{cl}

\begin{nb}
Note that the terms in $\eqref{e202}$ and $\eqref{hs45}$ cannot be controlled directly and have to be seen as time derivatives.
\end{nb}
Let us now use Claim $\ref{T11}$ to derive the proof of the differential inequality $\eqref{310}$ satisfied by $\alpha_{-}$, then we will prove it later.
\\{\it Proof of $\eqref{310}$ admitting Claim $\ref{T11}$}: 
 In order to avoid unecessary repetition, we kindly refer the interested reader to Merle and Zaag \cite{fh} (the proof of inequality (186) page $107$). Only we focus here on the new terms coming from $f$ defined in $\eqref{fy}$ and the new solution $\overline{w}$ defined in $\eqref{87}$ and $\eqref{e505}$.
\\In fact, the whole proof of $\eqref{310}$ is based on the fact that the derivative of $\alpha_{-}^2$ is related to the quadratic form $\varphi_{d}(q_{-},L_{d}(q_{-}))$, defined in $\eqref{500}$, which inherits the properties of the Lyapunov functional defined in $\eqref{e90}$ (and give an almost self-adjoint behavior). 
\\We use here the same bilinear form $\varphi_{d}$ introduced in Merle and Zaag \cite{fh}, therefore using the bound $\eqref{305}$ on $|d'|$, we get 
\begin{equation}\label{e205}
|\alpha_{-}\alpha_{-}'-\varphi_{d}(q_{-}, \partial_{s}q_{-})|\leq C|d'|\frac{\alpha_{-}^2}{1-d^2}\leq C\Vert q\Vert^4_{\H}+\frac{C}{s^{a}}\Vert q\Vert^3_{\H}.
 \end{equation}
For the reader's interest, we mention that the proof of inequality $\eqref{e205}$ is written in inequality $(206)$ page $107$ in Merle and Zaag \cite{fh}.
\\From Claim $\ref{T11}$, the continuity of $\varphi_{d}$ $\eqref{oa}$, the equivalence norms $\eqref{raw}$ and inequality $\eqref{e205}$, we write 
$$\Big|\alpha_{-}\alpha_{-}'-\varphi_{d}(q_{-},L_{d}(q_{-}))-\frac{d}{ds}\int_{-1}^{1}H(d,y,q_{1})\rho dy-\frac{d}{ds}\int_{-1}^{1}\widehat{F}(\overline{w}_{1},q_{1},s)\rho dy\Big| $$
 \begin{eqnarray}\label{e210}
 &\leq & C\Vert q\Vert^2_{\H} +\frac{C}{s^{a}}\Vert q\Vert^{2}_{\H}+\Big|\varphi_{d}\Big(q_{-},\partial_{s}q_{-}-L_{d}(q_{-})-\Sigma_{- ,d}^{1}(s)-\Sigma_{- ,d}^{2}(s)-\Sigma_{- ,d}^{3}(s)\Big)\Big|\nonumber\\
&\leq &C\Vert q\Vert^2_{\H} +\frac{C}{s^{a}}\Vert q\Vert^2_{\H}+C\Vert q_{-}\Vert_{\H}\Vert q\Vert^3_{\H}\leq C\Vert q\Vert^{1+\overline{p}}_{\H} +\frac{C}{s^{a}}\Vert q\Vert^{2}_{\H}.
 \end{eqnarray} 
We use now the expression $\eqref{500}$ of $\varphi_{d}$ and $\eqref{60}$ of $L_{d}$ to write
\begin{equation}\label{e206}
 \varphi_{d}(q_{-},L_{d}(q_{-}))=-\frac{4}{p-1}\int_{-1}^{1}q_{-,2}^2\frac{\rho}{1-y^2}dy.
 \end{equation}
 Let us recall that $\varphi_{d}(q,r)$ is the same bilinear form used by Merle and Zaag \cite{fh}. For that reason, we refer the reader to page $107$ and the beginning of page $108$ in Merle and Zaag \cite{fh} to see the details of the proof of equality $\eqref{e206}$. 
\\Using $\eqref{e205}$, $\eqref{e210}$ and $\eqref{e206}$, we see that the estimate $\eqref{310}$ holds with 
 \begin{equation}\label{e207}
 R_{-}(s)=-\int_{-1}^{1}H(d,y,q_{1})\rho dy-\int_{-1}^{1}\widehat{F}(\overline{w}_{1},q_{1},s)\rho dy.
 \end{equation}
Using $\eqref{32}$ of Lemma $\ref{T9}$, $\eqref{ena50}$ of Lemma $\ref{T10}$, condition $\eqref{200}$ (considering first the case $p\geq 2$ and then the case $1<p<2$), we see that $\eqref{307}$ holds. It remains to prove Claim $\ref{T11}$ in order to conclude the proof of $(i)-(ii)$ of Proposition $\ref{T8}$.

\bigskip

{\it Proof of Claim $\ref{T11}$}: 
\\{\bf Proof of $\eqref{e200}$}: We first project equation $\eqref{21}$ and using the negative projector $\pi_{-}^{d}$ introduced in $\eqref{18}$:
\begin{equation}\label{e211}
\pi_{-}^{d}(\partial_{s}q)=\pi_{-}^{d}(\overline{L}_{d}q)-d'\pi_{-}^{d}\left(
  \begin{array}{ccc}
   \partial_{d}\kappa(d,y) \\
    0\\
  \end{array}
\right)+\Sigma_{- ,d}^{2}(s)+\Sigma_{- ,d}^{3}(s),
 \end{equation}
where $\Sigma_{- ,d}^{2}$ and $\Sigma_{- ,d}^{3}$ are defined respectively in $\eqref{sta1}$ and $\eqref{sta2}$. According to the proof of Claim $5.4$ page $108$ in \cite{fh}, the Remark $\ref{Tttaw7}$, equations $\eqref{sta}$ and $\eqref{e211}$, we can directly write the following
\begin{equation}\label{arb100}
 \Vert \partial_{s}q_{-}-L_{d}q_{-}-\Sigma_{- ,d}^{1}(s)-\Sigma_{- ,d}^{2}(s)
-\Sigma_{- ,d}^{3}(s) \Vert_{\H} 
\leq  C(\alpha_{1}^2+\alpha_{-}^{2})^{\frac{3}{2}},
 \end{equation}
which end the proof of $\eqref{e200}$. 
\\{\bf Proof of $\eqref{2030}$:} This inequality is a direct consequence of the continuity of $\varphi_{d}$ $\eqref{oa}$, from the equivalence norms in $\eqref{raw}$ and $\eqref{raw1}$ and inequality $\eqref{100000}$ in Lemma $\ref{T10}$. 
\\{\bf Proof of $\eqref{rama}$:} This inequality is a direct consequence of the continuity of $\varphi_{d}$, $\eqref{oa}$, inequality $\eqref{305}$ on $d'$ and the a priori estimate $\eqref{200}$. 
\\{\bf Proof of $\eqref{e201}$ and $\eqref{hs60}$:} Recall from $\eqref{202}$ that we have 
\begin{equation}\label{e500}
 q(y,s)=\alpha_{1}F_{1}^{d}(y)+q_{-}(y,s),\nonumber
 \end{equation}
 \begin{equation}\label{e501}
 \left(
  \begin{array}{ccc}
    0 \\
    h(d,y,q_{1})\\
  \end{array}
\right)=\beta_{1}(s)F_{1}^{d}(y)+\beta_{0}(s)F_{0}^{d}(y)+\pi_{-}^{d}\left(
  \begin{array}{ccc}
    0 \\
    h(d,y,q_{1})\\
  \end{array}
\right),\nonumber
 \end{equation}
 \begin{equation}\label{hs61}
 \left(
  \begin{array}{ccc}
    0 \\
    \widehat{f}(\overline{w}_{1},q_{1},s)\\
  \end{array}
\right)=\widehat{\beta}_{1}(s)F_{1}^{d}(y)+\widehat{\beta}_{0}(s)F_{0}^{d}(y)+\pi_{-}^{d}\left(
  \begin{array}{ccc}
    0 \\
    \widehat{f}(\overline{w}_{1},q_{1},s)\\
  \end{array}
\right),\nonumber
 \end{equation}
 where $\beta_{\lambda}(s)=\pi_{\lambda}^{d}\left(
  \begin{array}{ccc}
    0 \\
    h(d,y,q_{1})\\
  \end{array}
\right)$ and $\widehat{\beta}_{\lambda}(s)=\pi_{\lambda}^{d}\left(
  \begin{array}{ccc}
    0 \\
    \widehat{f}(\overline{w}_{1},q_{1},s)\\
  \end{array}
\right)$. 
Note from the definition $\eqref{500}$ of $\varphi_{d}$, we have
\begin{equation}\label{et&àà}
\int_{-1}^{1}q_{2}h(d,y,q_{1})\rho dy =\varphi_{d}\Big(q,\left(
  \begin{array}{ccc}
    0 \\
    h(d,y,q_{1})\\
  \end{array}
\right)\Big).
\end{equation}
In addition by virtue of $\eqref{et&àà}$, the bilinearity of $\varphi_{d}$, the bound $\eqref{62}$ on the norm of $F_{\lambda}^{d}$ and the equivalence of norms in $\eqref{raw}$, we write 
 \begin{eqnarray}\label{binthabi}
 \Big|\varphi_{d}\Big(q_{-},\pi_{-}^{d}\left(
  \begin{array}{ccc}
    0 \\
    h(d,y,q_{1})\\
  \end{array}
\right)\Big)-\int_{-1}^{1}q_{2}h(d,y,q_{1})\rho dy \Big|&\leq & C (|\alpha_{1}|+|\alpha_{-}|)(|\beta_{1}|+|\beta_{0}|)\nonumber\\
&&+|\alpha_{1}|\Big|\varphi_{d}(F_{1}^{d},\left(
  \begin{array}{ccc}
    0 \\
    h(d,y,q_{1})\\
  \end{array}
\right))\Big|.\,\,\,\,\,\,\,\,\,\,\,\,\,\,\,\,\,\,
\end{eqnarray}
Since, we have from the expression $\eqref{500}$ of $\varphi_{d}$, the fact that $|F_{1,2}^{d}(y)|\leq C  \kappa(d,y) $ and Lemma $\ref{T9}$,      
 \begin{equation}\label{e502}
 \Big|\varphi_{d}(F_{1}^{d},\left(
  \begin{array}{ccc}
    0 \\
    h(d,y,q_{1})\\
  \end{array}
\right))\Big|=\Big|\int_{-1}^{1}F_{1,2}^{d}(y)h(d,y,q_{1})\rho dy \Big|\leq C(\alpha_{1}^2+\alpha_{-}^{2}),
 \end{equation}
 \begin{equation}\label{e503}
 |\beta_{1}|+|\beta_{0}|\leq C\int_{-1}^{1}\kappa(d,y) h(d,y,q_{1})\rho dy\leq C(\alpha_{1}^2+\alpha_{-}^{2}),
 \end{equation}
 we combine $\eqref{binthabi}$, $\eqref{e502}$ and $\eqref{e503}$ to conclude $\eqref{e201}$.
 \\Note from the definition $\eqref{500}$, the bilinearity of $\varphi_{d}$, the bound $\eqref{62}$ on the norm of $F_{\lambda}^{d}$ and the equivalence norm $\eqref{raw}$, we derive that  
 $$\Big|\varphi_{d}\Big(q_{-},\pi_{-}^{d}\left(
  \begin{array}{ccc}
    0 \\
    \widehat{f}(\overline{w}_{1},q_{1},s)\\
  \end{array}
\right)\Big)-\int_{-1}^{1}q_{2}\widehat{f}(\overline{w}_{1},q_{1},s)\rho dy \Big|$$
 $$=\Big|\varphi_{d}\Big(q_{-},\pi_{-}^{d}\left(
  \begin{array}{ccc}
    0 \\
    \widehat{f}(\overline{w}_{1},q_{1},s)\\
  \end{array}
\right)\Big)-\varphi_{d}\Big(q,\left(
  \begin{array}{ccc}
    0 \\
    \widehat{f}(\overline{w}_{1},q_{1},s)\\
  \end{array}
\right)\Big)\Big|$$
 $$\leq C (|\alpha_{1}|+|\alpha_{-}|)(|\widehat{\beta}_{1}|+|\widehat{\beta}_{0}|)+|\alpha_{1}|\Big|\varphi_{d}(F_{1}^{d},\left(
  \begin{array}{ccc}
    0 \\
    \widehat{f}(\overline{w}_{1},q_{1},s)\\
  \end{array}
\right))\Big|.$$
Since, we have from the expression $\eqref{500}$ of $\varphi_{d}$, the fact that $|F_{1,2}^{d}(y)|\leq C  \kappa(d,y) $, the a priori estimate $\eqref{200}$, $\eqref{e15}$, $\eqref{ena200}$ and $\eqref{jar}$,
 \begin{equation}\label{hs65}
 \Big|\varphi_{d}(F_{1}^{d},\left(
  \begin{array}{ccc}
    0 \\
    \widehat{f}(\overline{w}_{1},q_{1},s)\\
  \end{array}
\right))\Big| =\Big|\int_{-1}^{1}F_{1,2}^{d}(y)\widehat{f}(\overline{w}_{1},q_{1},s)\rho dy \Big|
\leq  C(\alpha_{1}^2+\alpha_{-}^{2}) 
 \end{equation}
 and from $\eqref{jar1}$, we obtain
 \begin{equation}\label{hs66}
 |\widehat{\beta}_{1}|+|\widehat{\beta}_{0}|\leq C\int_{-1}^{1}\kappa(d,y) \widehat{f}(\overline{w}_{1},q_{1},s)\rho dy\leq C(\alpha_{1}^2+\alpha_{-}^{2}),
 \end{equation}
 this gives $\eqref{hs60}$.
 \\{\bf Proof of $\eqref{e202}$:} Since $q_{2}=\partial_{s}q_{1}+d'\partial_{d}\overline{w}_{1}$ by $\eqref{301}$ and $\eqref{e505}$, we use the expression $\eqref{33}$ of $H$ to write 
 $$$$
 \begin{eqnarray}\label{e504}
\int_{-1}^{1}q_{2}h(d,y,q_{1})\rho dy&=&\frac{d}{ds}\int_{-1}^{1}H(d,y,q_{1})\rho dy+d'\int_{-1}^{1} \partial_{d}\overline{w}_{1}\Big(h(d,y,q_{1})-\partial_{d}H(d,y,q_{1})\Big)\rho dy\nonumber\\
&&-\int_{-1}^{1}\frac{\partial H(d,y,q_{1})}{\partial\overline{w}_{1} }\frac{\partial\overline{w}_{1}}{\partial s }\rho dy\nonumber\\
&=&\frac{d}{ds}\int_{-1}^{1}H(d,y,q_{1})\rho dy+ d'\frac{p(p-1)}{2}\int_{-1}^{1}\partial_{d}\overline{w}_{1}\overline{w}_{1}^{p-2}q^2_{1}\rho dy\nonumber\\
 &&-\int_{-1}^{1}\Big( h(d,y,q_{1})+\frac{p(p-1)}{2}\overline{w}_{1}^{p-2}q^2_{1} \Big)\frac{\partial\overline{w}_{1}}{\partial s }\rho dy .
 \end{eqnarray}
Since, we see from the Hardy-Sobolev inequality of Lemma $\ref{T21}$, inequality $\eqref{e54}$ and $\eqref{fek}$ that the function $\phi$ defined in $(iii)$ of Lemma $\ref{T5}$ is bounded, it follows that $\Vert \partial_{d}\overline{w}_{1}\overline{w}_{1}^{p-2}\Vert_{L^{\frac{p+1}{p-1}}_{\rho}} \leq \frac{C}{1-d^2}$. From this fact, we derive the following:
 \begin{equation}\label{hs41}
 \int_{-1}^{1}\partial_{d}\overline{w}_{1}\overline{w}_{1}^{p-2}q^2_{1}\rho dy\leq \frac{C}{1-d^2}\Vert q_{1}\Vert^2_{L^{p+1}_{\rho}}\leq \frac{C(\alpha_{1}^2+\alpha_{-}^{2})}{1-d^2}.
 \end{equation}
 From Lemma $\ref{T9}$, Corollary $\ref{T15}$ and the H\"{o}lder inequality, it follows that 
\begin{eqnarray}\label{na3}
 \int_{-1}^{1}\Big( h(d,y,q_{1})+\frac{p(p-1)}{2}\overline{w}_{1}^{p-2}q^2_{1} \Big)\frac{\partial\overline{w}_{1}}{\partial s }\rho dy &\leq &\frac{C}{s^{a}}\Vert q_{1}\Vert^{2}_{L^{p+1}_{\rho}}\Vert \kappa \Vert^{p-1}_{L^{p+1}_{\rho}}\nonumber\\
&& +\frac{C\delta_{p\geq 2}}{s^{a}}\Vert q_{1}\Vert^{p}_{L^{p+1}_{\rho}}\Vert \kappa \Vert_{L^{p+1}_{\rho}},
 \end{eqnarray}
 where $\delta_{p\geq 2}$ is $0$ if $1<p<2$ and $1$ otherwise. Therefore, using $(ii)$ of Lemma $\ref{T21}$, the a priori estimate $\eqref{200}$, the equivalence of norms $\eqref{raw1}$ and $\eqref{na3}$, we write
 \begin{equation}\label{hs1000}
\int_{-1}^{1}\Big( h(d,y,q_{1})+\frac{p(p-1)}{2}\overline{w}_{1}^{p-2}q^2_{1} \Big)\frac{\partial\overline{w}_{1}}{\partial s }\rho dy \leq \frac{C}{s^{a}}(\alpha_{1}^2+\alpha_{-}^{2}).
 \end{equation}
Finally from the bound $\eqref{305}$ on $d'$, the a priori estimate $\eqref{200}$, identity $\eqref{e504}$, inequalities $\eqref{hs41}$ and $\eqref{hs1000}$, we deduce that
 \begin{equation}\label{hs42}
 \Big|\int_{-1}^{1}q_{2}h(d,y,q_{1})\rho dy -\frac{d}{ds}\int_{-1}^{1}H(d,y,q_{1})\rho dy\Big| \leq C(\alpha_{1}^2+\alpha_{-}^{2})^2+\frac{C}{s^{a}}(\alpha_{1}^2+\alpha_{-}^{2}).
  \end{equation}
{\bf Proof of $\eqref{hs45}$:} We write first
\begin{eqnarray}\label{hs44}
\widehat{F}(\overline{w}_{1},q_{1},s)&=&\int_{0}^{q_{1}}\widehat{f}(\overline{w}_{1},q',s)dq'=e^{-\frac{2(p+1)}{p-1}s}\Big(F(e^{\frac{2s}{p-1}}(\overline{w}_{1}(d,y,s)+q_{1}))\\
&&-F(e^{\frac{2s}{p-1}}\overline{w}_{1}(d,y,s))-e^{\frac{2s}{p-1}}q_{1}f(e^{\frac{2s}{p-1}}\overline{w}_{1}(d,y,s))-e^{\frac{4s}{p-1}}\frac{q^2_{1}}{2}f'(e^{\frac{2s}{p-1}}\overline{w}_{1}(d,y,s))\Big),\nonumber
  \end{eqnarray}
with $F$ and $\widehat{f}$ defined respectively in $\eqref{90}$ and $\eqref{40}$. Since $q_{2}=\partial_{s}q_{1}+d'\partial_{d}\overline{w}_{1}$ by $\eqref{301}$ and $\eqref{e505}$, thanks to the expression $\eqref{hs44}$ of $\widehat{F}$, we write
 \begin{eqnarray}\label{hs46}
\int_{-1}^{1}q_{2}\widehat{f}(\overline{w}_{1},q_{1},s)\rho dy &=&\int_{-1}^{1}\partial_{s}q_{1}\widehat{f}(\overline{w}_{1},q_{1},s)\rho dy + d'\int_{-1}^{1}\partial_{d}\overline{w}_{1}\widehat{f}(\overline{w}_{1},q_{1},s)\rho dy\nonumber\\
&=&\frac{d}{ds}\int_{-1}^{1}\widehat{F}(\overline{w}_{1},q_{1},s)\rho dy+d'\int_{-1}^{1} \partial_{d}\overline{w}_{1}\Big(\widehat{f}(\overline{w}_{1},q_{1},s)-\frac{\partial \widehat{F}(\overline{w}_{1},q_{1},s)}{\partial\overline{w}_{1} }\Big)\rho dy\nonumber\\
&&-  \int_{-1}^{1}\frac{\partial \widehat{F}(\overline{w}_{1},q_{1},s)}{\partial\overline{w}_{1} }\frac{\partial\overline{w}_{1}}{\partial s }\rho dy-\int_{-1}^{1}\frac{\partial \widehat{F}(\overline{w}_{1},q_{1},s)}{\partial s }\rho dy.
 \end{eqnarray}
 By a careful calculation, we can see that 
 \begin{equation}\label{hs48}
 \partial_{d}\overline{w}_{1}\Big(\widehat{f}(\overline{w}_{1},q_{1},s)-\frac{\partial \widehat{F}(\overline{w}_{1},q_{1},s)}{\partial\overline{w}_{1} }\Big)=e^{-\frac{2(p-2)}{p-1}s}\frac{q_{1}^2}{2}\partial_{d}\overline{w}_{1}f''(e^{\frac{2s}{p-1}}\overline{w}_{1}).
 \end{equation}
We derive from $\eqref{hs44}$
\begin{eqnarray}\label{hs101}
 \frac{\partial \widehat{F}(\overline{w}_{1},q_{1},s)}{\partial s }&=&-\frac{2(p+1)}{p-1}\widehat{F}(\overline{w}_{1},q_{1},s)+g_{p}(A+X)-g_{p}(A)\nonumber\\
&&-Xg_{p}'(A )-\frac{X^2}{2}g''_{p}(A),
 \end{eqnarray}
where $A(y,s)=e^{\frac{2s}{p-1}}\overline{w}_{1}$, $X(y,s)=e^{\frac{2s}{p-1}}q_{1}$ and the function $g_{p}(X)=\frac{2}{p-1}e^{-\frac{2(p+1)}{p-1}s}Xf(X)$.
\\We apply the mean value theorem to say that there exists $\theta \in (0,1)$ such that:
\begin{equation}\label{70}
g_{p}(A+X)-g_{p}(A)-Xg_{p}'(A )-\frac{X^2}{2}g''_{p}(A)=\frac{X^3}{6}g_{p}^{(3)}(A+\theta X ).
\end{equation}
Remarking that 
\begin{equation}\label{gnew1}
g_{p}^{(3)}(A+\theta X )=\frac{2}{p-1}e^{-\frac{2(p+1)}{p-1}s}(3f''(A+\theta X )+(A+\theta X)f^{(3)}(A+\theta X )).
\end{equation}
After straightforward computations of the expressions of $f'$, $f''$ and $f^{(3)}$, we write the following
\begin{equation}\label{s2o321}
\frac{|X|^3}{6}|g_{p}^{(3)}(A+\theta X )|\leq C|q_{1}|^{3}|\overline{w}_{1}+\theta q_{1}|^{p-2}.
\end{equation}
Hence, if $\frac{|q_{1}|}{|\overline{w}_{1}|}\geq \frac{1}{2}$, we obtain
\begin{equation}\label{s22o321}
\frac{|X|^3}{6}|g_{p}^{(3)}(A+\theta X )|\leq C|q_{1}|^{p+1}.
\end{equation}
For the case when $\frac{|q_{1}|}{|\overline{w}_{1}|}< \frac{1}{2}$, similarly to the proof of $\eqref{so311}$ and $\eqref{so360}$, we derive if $p\geq 2$
\begin{equation}\label{s21o321}
\frac{|X|^3}{6}|g_{p}^{(3)}(A+\theta X )|\leq C|q_{1}|^{3}|\overline{w}_{1}|^{p-2}
\end{equation}
and if $1<p<2$, we obtain
\begin{equation}\label{s23o321}
\frac{|X|^3}{6}|g_{p}^{(3)}(A+\theta X )|\leq C|q_{1}|^{p+1}.
\end{equation}
Adding estimates $\eqref{s22o321}$, $\eqref{s21o321}$ and $\eqref{s23o321}$, we write
\begin{equation}\label{s24o321}
\frac{|X|^3}{6}|g_{p}^{(3)}(A+\theta X )|\leq C\delta_{p\geq 2}|\overline{w}_{1}|^{p-2}|q_{1}|^{3}+C|q_{1}|^{p+1},
\end{equation}
where $\delta_{p\geq 2}$ is $0$ if $1<p<2$ and $1$ otherwise.
Collecting $\eqref{hs101}$, $\eqref{ena50}$ in Lemma $\ref{T10}$, inequality $\eqref{s24o321}$ and the a priori estimate $\eqref{200}$ together combined with the Hardy-Sobolev inequality, we deduce that 
\begin{equation}\label{hs111}
 |\int_{-1}^{1}\frac{\partial \widehat{F}(\overline{w}_{1},q_{1},s)}{\partial s }\rho dy|\leq C \Vert q\Vert^{1+\overline{p}}_{\H},
 \end{equation}
where $\overline{p}$ is defined in $\eqref{307}$. Now we deal with the term written in $\eqref{hs48}$. From $\eqref{ap1}$, $\eqref{305}$ and the a priori estimate $\eqref{200}$ combined with the Hardy-Sobolev inequality, we obtain
 \begin{equation}\label{hs51}
|d'\int_{-1}^{1} \partial_{d}\overline{w}_{1}\Big(\widehat{f}(\overline{w}_{1},q_{1},s)-\frac{\partial \widehat{F}(\overline{w}_{1},q_{1},s)}{\partial\overline{w}_{1} }\Big)\rho dy|\leq C\Vert q\Vert^4_{\H}+\frac{C}{s^{a}}\Vert q\Vert^3_{\H}.
\end{equation}
We apply the mean value theorem to say that there exists $\theta \in (0,1)$ such that:
\begin{equation}\label{hs0151}
|\frac{\partial \widehat{F}(\overline{w}_{1},q_{1},s)}{\partial\overline{w}_{1} }|\leq C|q_{1}|^{2}|\overline{w}_{1}+\theta q_{1}|^{p-2}+C|q_{1}|^{2}|\overline{w}_{1}|^{p-2}.
\end{equation}
According to inequality $\eqref{hs0151}$ and $\eqref{ap2}$ combined with the Hardy-Sobolev inequality, we can see that 
 \begin{equation}\label{hs112}
 |\int_{-1}^{1}\frac{\partial \widehat{F}(\overline{w}_{1},q_{1},s)}{\partial\overline{w}_{1} }\frac{\partial\overline{w}_{1}}{\partial s }\rho dy|\leq \frac{C}{s^{a}}\Vert q\Vert^{2}_{\H}.
\end{equation}
Gathering the bounds $\eqref{hs46}$, $\eqref{hs111}$, $\eqref{hs51}$, $\eqref{hs112}$ and the a priori estimate $\eqref{200}$, it follows that 
 \begin{equation}\label{hs52}
 \Big |\int_{-1}^{1}q_{2}\widehat{f}(\overline{w}_{1},q_{1},s)\rho dy -\frac{d}{ds}\int_{-1}^{1}\widehat{F}(\overline{w}_{1},q_{1},s)\rho dy \Big | \leq C\Vert q\Vert^{1+\overline{p}}_{\H}+ \frac{C}{s^{a}}\Vert q\Vert^{2}_{\H}.
 \end{equation}
 Finally, $\eqref{hs42}$ and $\eqref{hs52}$ ends the proof of Claim $\ref{T11}$ as well as $(i)$ and $(ii)$ of Proposition $\ref{T8}$.
\\$(iii)$ This inequality is a consequence of the coercivity of the quadratic form $\varphi_{d}$ on the space $\H_{-}^{d}$ stated in $\eqref{raw}$ and $\eqref{raw1}$.
\\From equation $\eqref{21}$, the fact that $q_{2}=\partial_{s}q_{1}+d'\partial_{d}\overline{w}_{1}$ and the definition $\eqref{60}$ of $L_{d}$, we write 
 \begin{eqnarray}\label{e550}
 \frac{d}{ds}\int_{-1}^{1}q_{1}q_{2}\rho dy&=&\int_{-1}^{1}q_{2}\partial_{s}q_{1}\rho dy+\int_{-1}^{1}q_{1}\partial_{s}q_{2}\rho dy \nonumber \\
 &=&\int_{-1}^{1}q_{2}(q_{2}-d'\partial_{d}\overline{w}_{1})\rho dy \\
 && +\int_{-1}^{1}q_{1}\Big( \pounds q_{1}+\psi(d,y)q_{1}-\frac{p+3}{p-1}q_{2}-2y\partial_{y}q_{2}+h(d,y,q_{1})\Big)\rho dy \nonumber \\
&& +\int_{-1}^{1}\overline{V}q^2_{1}\rho dy-d'\int_{-1}^{1}q_{1} \partial_{d}\overline{w}_{2}\rho dy+\int_{-1}^{1}q_{1}\widehat{f}(\overline{w}_{1},q_{1},s) \rho dy.\nonumber
\end{eqnarray}
 According to $\eqref{e550}$ and the proof of item $(iii)$ of Proposition $5.2$ page $110$ in Merle and Zaag \cite{fh}, we write directly that 
 \begin{eqnarray}\label{pins}
 \frac{d}{ds}\int_{-1}^{1}q_{1}q_{2}\rho dy&\leq & -\frac{4}{5}\alpha_{-}^2+C\int_{-1}^{1} q_{-,2}^2\frac{\rho}{1-y^2}dy+C\alpha_{1}^2\nonumber\\
 &&+\int_{-1}^{1}\overline{V}q^2_{1}\rho dy-d'\int_{-1}^{1}q_{1} \partial_{d}\overline{w}_{2}\rho dy+\int_{-1}^{1}q_{1}\widehat{f}(\overline{w}_{1},q_{1},s) \rho dy.\nonumber
\end{eqnarray}
We see from $\eqref{100000}$, $\eqref{na3}$ and the fact that $s^*$ is large enough such that for all $s\geq s^*$, the following estimate holds
\begin{equation}\label{e561}
\Big |\int_{-1}^{1}\overline{V}q^2_{1}\rho dy\Big |\leq Ce^{-s}\int_{-1}^{1}q^2_{1}\rho dy+\frac{C}{s^{a}}\int_{-1}^{1}|\kappa(d,y)|^{p-1}q^2_{1}\rho dy\leq \frac{C}{s^{a}}\Vert q\Vert^2_{\H}\leq \frac{1}{5}\Vert q\Vert^2_{\H}.
\end{equation}
We can see from $\eqref{e11}$, the smallness condition $\eqref{200}$ and Lemma $\ref{T21}$, we have
\begin{eqnarray}\label{e562}
\Big |\int_{-1}^{1}q_{1}\widehat{f}(\overline{w}_{1},q_{1},s) \rho dy\Big |&\leq &C\delta_{p\geq 2}\int_{-1}^{1}\kappa^{p-2}(d,y)|q_{1}|^3\rho dy+C\int_{-1}^{1}|q_{1}|^{p+1}\rho dy\nonumber\\
&&\leq C\Vert \kappa \Vert^{p-2}_{L^{p+1}_{\rho}}\Vert q_{1}\Vert^{3}_{L^{p+1}_{\rho}}+C\delta_{p\geq 2}\Vert q_{1}\Vert^{p+1}_{L^{p+1}_{\rho}}\nonumber\\
&&\leq C\Vert q\Vert^{\overline{p}+1}_{\H}\leq \frac{1}{5}\Vert q\Vert^{2}_{\H}.
\end{eqnarray}
Using $\eqref{e51}$ to write 
\begin{equation}\label{e564}
\Big |\int_{-1}^{1}q_{1} \partial_{d}\overline{w}_{2}\rho dy\Big |\leq \frac{C}{s^{a}(1-d^2)}\int_{-1}^{1}|q_{1}|\kappa(d,y)\rho dy.
\end{equation}
According to $\eqref{e564}$, Lemma $\ref{T21}$, the smallness condition $\eqref{200}$, $(i)$ of Proposition $\ref{T8}$ and the fact that $s^*$ is large enough, we deduce that for all $s\geq  s^*$, we have
\begin{equation}\label{e566}
\Big |d'\int_{-1}^{1}q_{1} \partial_{d}\overline{w}_{2}\rho dy\Big |\leq \frac{C|d'|}{s^{a}(1-d^2)}\Vert q\Vert_{\H}\leq \frac{C}{s^{a}}\Vert q\Vert^3_{\H}+\frac{C}{s^{2a}}\Vert q\Vert^2_{\H}\leq \frac{1}{5}\Vert q\Vert^2_{\H}.
\end{equation}
Combining $\eqref{e550}$ with $\eqref{e561}$, $\eqref{e562}$ and $\eqref{e566}$ to deduce the proof of $(iii)$ of Proposition $\ref{T8}$.
\\$(iv)$ We use first the hypothesis $\eqref{111}$ of Theorem $\ref{T4}$ and then we search to an upper bound to the functional $E(w(s),s)$. Recalling that 
 $$E(w(s),s)=E_{0}(w(s))+I(w(s),s)+J(w(s),s),$$
with 
\begin{equation}\label{hs3}
I(w(s),s)\hspace{-0.1cm}=\hspace{-0.1cm}-e^{\frac{-2(p+1)s}{p-1}} \int_{-1}^{1} F(e^{\frac{2s}{p-1}}w)\rho {\mathrm{d}}y\,{\rm and}\,J(w(s),s)\hspace{-0.1cm}=\hspace{-0.1cm}-\frac{1}{s^{\frac{a+1}{2}}}\int_{-1}^{1} w\partial_{s}w \rho dy.
\end{equation} 
Using the definition of $q(y,s)$ written in $\eqref{301}$, we can make an expansion of $E_{0}(w(s))$ defined in $\eqref{cp}$ for $q\longrightarrow 0$ in $\H$ and get after some straightforward computations
\begin{eqnarray}\label{br4}
E_{0}(w(s)) &= & E_{0}(\overline{w}_{1}(d,y,s))+\frac{1}{2}\varphi_{d}(q,q)-\int_{-1}^{1}H(d,y,q_{1})\rho {\mathrm{d}}y \nonumber\\
&&+R_{1}(s)+R_{2}(s)+R_{3}(s),
\end{eqnarray}
with 
\begin{itemize}
\item $R_{1}(s)=\int_{-1}^{1}\Big(-\pounds\overline{w_{1}}+2\frac{p+1}{(p-1)^2}\overline{w_{1}}-|\overline{w_{1}}|^{p}\Big)q_{1}\rho  {\mathrm{d}}y,$
\item $R_{2}(s)=\frac{p}{2}\int_{-1}^{1}\Big(\kappa^{p-1}-\overline{w}_{1}^{p-1}\Big )q_{1}^2\rho  {\mathrm{d}}y,$
\item $R_{3}(s)=\frac{1}{2}d'^2\int_{-1}^{1}(\partial_{d}\overline{w}_{1})^2\rho  {\mathrm{d}}y+d'\int_{-1}^{1}\partial_{s}\overline{w}_{1}\partial_{d}\overline{w}_{1}\rho  {\mathrm{d}}y-d'\int_{-1}^{1}\partial_{d}\overline{w}_{1}q_{2}\rho  {\mathrm{d}}y+\int_{-1}^{1}\overline{w}_{2}q_{2}\rho  {\mathrm{d}}y.$
\end{itemize}
Recalling that $\varphi_{d}$ (the same bilinear form used by Merle and Zaag \cite{fh}) satisfies the following
 \begin{equation}\label{br8}
 \varphi_{d}(q,q)\leq c_{2} \alpha_{-}^2-c_{3}\alpha_{1}^2 ,
\end{equation}
for some $c_{2} > 0$ and $c_{3} > 0$. Using $\eqref{200}$, $\eqref{raw1}$, $\eqref{307}$ and $\eqref{e207}$, to obtain 
  \begin{equation}\label{br9}
- \int_{-1}^{1}H(d,y,q_{1})\rho {\mathrm{d}}y\leq C\Vert q\Vert^{\overline{p}+1}_{\H}\leq C\epsilon^{\overline{p}-1}(\alpha^2_{1}(s)+\alpha^2_{-}(s)),
 \end{equation}
 where $\overline{p}=\min(p,2)$.
Recalling that $\kappa (d,y)$ defined in $\eqref{400}$ is a stationary solution of the unperturped case, so we can write the following
\begin{equation}\label{br20}
\pounds \kappa (d,y)-2\frac{p+1}{(p-1)^2}\kappa (d,y)+|\kappa (d,y)|^{p}=0,
\end{equation}
fore more detail of equality $\eqref{br20}$ we can see the proof of equality $(47)$ page 58 in Merle and Zaag \cite{fh}. It is easy to check
\begin{equation}\label{br25} 
R_{1}(s)\hspace{-0.1cm}=\hspace{-0.1cm}\int_{-1}^{1}\Big(\pounds(\kappa (d,y)-\overline{w_{1}})+2\frac{p+1}{(p-1)^2}(\overline{w_{1}}-\kappa (d,y))+(|\kappa (d,y)|^{p}-|\overline{w_{1}}|^{p})\Big)q_{1}\rho  {\mathrm{d}}y.
\end{equation}
After an integration by parts, using the fact that $|a^p-1|\leq C|a-1|$ where $a$ is bounded, together with Claim $\ref{T14}$, the invariance of equation $\eqref{1}$ and the norm in $\H_{0}$ under the Lorentz transform written in $\ref{Tekk}$, we get 
\begin{equation}\label{br21} 
R_{1}(s)\leq \frac{C}{s^a}.
\end{equation}
We use again Claim $\ref{T14}$ and the fact that $|a^p-1|\leq C|a-1|$ where $a$ is bounded to write 
\begin{equation}\label{br10}
R_{2}(s) \leq \frac{C}{s^a}\Vert q\Vert^{2}_{\H}.
 \end{equation}
According to the classical inequality $ab\leq \frac{1}{2}(a^2+b^2)$, inequalities $\eqref{ap1}$, $\eqref{ap2}$ and $\eqref{305}$ on $d'$, we can derive
\begin{equation}\label{br11}
R_{3}(s)\leq \frac{C}{s^a}.
 \end{equation}
We use the expression $\eqref{301}$ of $w$ and exactly the same techniques used in the proof of Lemma 2.1 page 1121 in our paper \cite{X}, we prove
\begin{equation}\label{br5}
I(w(s),s)\leq \frac{C}{s^a}.
\end{equation}
From the expression $\eqref{301}$ of $w$ and $\eqref{hs3}$ of $J(w(s),s)$, we can write
  \begin{eqnarray}\label{br6}
J(w(s),s)&=&-\frac{1}{s^{\frac{a+1}{2}}}\int_{-1}^{1} \overline{w}_{1}\partial_{s}\overline{w}_{1} \rho dy-\frac{1}{s^{\frac{a+1}{2}}}\int_{-1}^{1} \overline{w}_{1}\partial_{s}q_{1} \rho dy\nonumber\\
&&-\frac{1}{s^{\frac{a+1}{2}}} \int_{-1}^{1} q_{1}\partial_{s}\overline{w}_{1} \rho dy -\frac{1}{s^{\frac{a+1}{2}}} \int_{-1}^{1} q_{1}\partial_{s}q_{1} \rho dy.
 \end{eqnarray}
Using the fact that $q_{2}=\partial_{s}q_{1}+d'\partial_{d}\overline{w}_{1}$, the classical inequality $ab\leq \frac{1}{2}(a^2+b^2)$, inequalities $\eqref{ap1}$, $\eqref{ap2}$, $\eqref{br6}$ and $\eqref{305}$ on $d'$, we can deduce
  \begin{equation}\label{br7}
  J(w(s),s)\leq  \frac{C}{s^{a}}.
  \end{equation}
To conclude we need to combine the expression of the functional $E(w(s),s)$ with condition $\eqref{111}$, identity $\eqref{br4}$, inequalities $\eqref{br8}$, $\eqref{br9}$, $\eqref{br21}$, $\eqref{br10}$, $\eqref{br11}$, $\eqref{br5}$, $\eqref{br7}$ together with the smallness condition $\eqref{200}$ and taking $\epsilon $ small enough so that $C\epsilon^{\overline{p}-1}\leq \frac{c_{3}}{4}$ to deduce that 
$$-\frac{C}{s^{\frac{a+1}{2}}}\leq E(w(s))- E_{0}(\overline{w}_{1}(d,y,s))\leq \Big(\frac{c_{2}}{2}+\frac{c_{3}}{4}\Big )\alpha_{-}^2-\frac{c_{3}}{4}\alpha_{1}^2+\frac{C}{s^{a}},$$
which yields $\eqref{309}$ and concludes the proof of Proposition $\ref{T8}$.

\Box

\begin{section}{Polynomial decay of the different components}\label{TS2222}
\end{section}
We are here mainly interested in proving the polynomial decay which may appear as a rough estimate compared to $\eqref{115}$. However it is a key result to get Theorem $\ref{T4}$. This estimate guarantee that $\Vert q(s)\Vert_{\H}\longrightarrow 0$ as $s \longrightarrow \infty$ where $q$ is defined in $\eqref{301}$. More precisely we prove the following proposition:
\begin{prop}\label{Tds}{\bf (Polynomial decay).} Assume that for $w\in \C([s^*,\infty ),\H)$ a solution of equation $\eqref{1}$ the following conditions hold:
$$\,\,\forall\,\,s\geq s^*\,\,\,E(w(s),s)\geq  E_{0}(\overline{w}_{1}(d,y,s) )-\frac{C}{s^{\frac{a+1}{2}}}$$
and 
$$\Big\Vert \left(
  \begin{array}{ccc}
    w(s^*) \\
    \partial_{s} w(s^*)\\
  \end{array}
\right)-\omega^*\left(
  \begin{array}{ccc}
    \overline{w}_{1}(d^*,.,s^*) \\
    \overline{w}_{2}(d^*,.,s^*)\\
  \end{array}
\right)\Big\Vert_{\H}\leq \epsilon^*,$$
for some $d^*\in (-1,1)$, $\omega^*\in \lbrace -1,1\rbrace$ and $\epsilon^*\in (0,\epsilon_{0}]$, where $\H$ and its norm are defined in $\eqref{110}$, $\overline{w}_{1}$ and $\overline{w}_{2}$ are defined respectively in $\eqref{87}$ and $\eqref{e505}$, there exists $d_{\infty}\in (-1,1)$ such that 
$$|d_{\infty}-d^*|\leq C\epsilon^*(1-d^{*2}).$$
Then there exists positive constant $C$ such that we have for all $s\geq s^*$,
 \begin{equation}\label{ds115}
\Big\Vert \left(
  \begin{array}{ccc}
    w(s) \\
    \partial_{s} w(s)\\
  \end{array}
\right)-\omega^*\left(
  \begin{array}{ccc}
    \overline{w}_{1}(d_{\infty},.,s) \\
   \overline{w}_{2}(d_{\infty},.,s)\\
  \end{array}
\right)\Big\Vert_{\H}\leq  \frac{C}{s^{\frac{a+1}{4}}},
\end{equation}
\end{prop}

\begin{nb}
Admitting Proposition $\ref{Tds}$, we can prove that Proposition $\ref{tawtaw}$ follows trivially from Proposition $\ref{Tds}$ by the classical triangular inequality and we use the fact that the solution $\overline{w}_{1}$ of equation $\eqref{1}$ approaches the stationary solutions of equation $\eqref{1}$ when ($f\equiv 0$) namely $\kappa (d,y)$ defined in $\eqref{400}$
\end{nb}
A good understanding of Proposition $\ref{Tds}$ gives us the permission to introduce a parameter $d(s)\in (-1,1)$ such that $\Vert q(s)\Vert_{\H}\longrightarrow 0$ as $s\,\,\longrightarrow \infty$ with $q$ is defined in $\eqref{301}$. Following this crucial information, we obtain the exponential decay and conclude our main Theorem $\ref{T4}$. In the following, we showed that if $w(s^*)$ is close enough to some class of solution $\overline{w}_{1}(d^*,.,s)$ and satisfies an energy barrier, then $w(s)$ converges to a neighboring class of solution as $s\rightarrow \infty$. Our aim is to show the convergence of $ \left(
  \begin{array}{ccc}
    w(s) \\
    \partial_{s} w(s)\\
  \end{array}
\right)$ as $s\rightarrow \infty$ to some $\left(
  \begin{array}{ccc}
    \overline{w}_{1}(d_{\infty},y,s) \\
   \overline{w}_{2}(d_{\infty},y,s)\\
  \end{array}
\right)$ for some $d_{\infty}$ close to $d^*\in (-1,1)$ in order to obtain Proposition $\ref{Tds}$.

 \bigskip

After the proof of the Proposition $\ref{T8}$, we can adapt with no difficulty the proof given in the case when $f\equiv 0$ treated by Merle and Zaag \cite{fh}. It happens that the same adaptation pattern works in the present case to obtain the polynomial decay. To be accurate and concise in our result, we are going to give the detail of the proof. Let us first introduce a more adapted notation and rewrite Proposition $\ref{T8}$.
\\If we introduce 
\begin{equation}\label{oaM}
\theta (s)=\frac{1}{2}\log \Big (\frac{1+d(s)}{1-d(s)}\Big ),\,\,\,A(s)=\alpha^2_{1}(s)\,\,\,{\rm and} \,\,\,B(s)=\alpha^2_{-}(s)+2R_{-}(s),
\end{equation}
(note that $d(s)=\tanh (\theta (s))$), then we see from $\eqref{raw1}$ and $\eqref{307}$ that if $\eqref{200}$ holds, then $|B(s)-\alpha^2_{-}|\leq C \epsilon^{\overline{p}-1}(\alpha^2_{1}+\alpha^2_{-})$, hence 
\begin{equation}\label{oaMM}
\frac{99}{100}\alpha^2_{-}-\frac{1}{100}A(s) \leq B(s) \leq \frac{101}{100}\alpha^2_{-}+\frac{1}{100}A(s)
\end{equation}
for $\epsilon$ small enough. Therefore, using Proposition $\ref{T8}$, estimates $\eqref{200}$, $\eqref{raw1}$ and the fact that $\theta' (s)=\frac{d'(s)}{1-d^2(s)}$, we derive the following:
\begin{cor}\label{Tnouv}{\bf (Equations in the new framework)}
There exist positive $\epsilon_{4}$, $K_{0}$, $K_{1}$ and $C_{i}$ for $i\in \lbrace 0,\,1,\,2,\,3\rbrace $ such that if $w$ is a solution to equation $\eqref{1}$ such that $\eqref{101}$ and $\eqref{200}$ hold for some $\epsilon\leq \epsilon_{4}$, where $q$ is defined in $\eqref{301}$, then using the notation $\eqref{oaM}$, we have for all $s\geq s^*$,
\begin{enumerate}[{\rm(i)}]
\item {\bf (Size of the solution)}
\begin{equation}\label{oa300}
\frac{1}{K_{0}}(A(s)+B(s))\leq \Vert q\Vert^2_{\H} \leq K_{0}(A(s)+B(s))\leq K^2_{0} \epsilon^2,
\end{equation}
\begin{equation}\label{oaest}
|\theta' (s)|\leq K_{0}(A(s)+B(s))+\frac{C_{0}}{s^a} \leq K_{0}^2\Vert q\Vert^2_{\H} +\frac{C_{0}}{s^a},
\end{equation}
\begin{equation}\label{oa301}
\Big|\int_{-1}^{1}q_{1}q_{2}\rho dy\Big|\leq  K_{0}(A(s)+B(s)) .
\end{equation}

\item {\bf (Equations)}
\begin{equation}\label{oa200}
\frac{3}{2}A(s)-K_{0} \epsilon B(s) -\frac{C_{1}}{s^a}\leq A'(s)\leq \frac{5}{2}A(s) +K_{0} \epsilon B(s) +\frac{C_{1}}{s^a},
\end{equation}
\begin{equation}\label{oa201}
B'(s)\leq -\frac{8}{p-1}\int_{-1}^{1}q_{-,2}^2\frac{\rho}{1-y^2}dy+K_{0} \epsilon (A(s)+B(s))+\frac{C_{2}}{s^a} ,
\end{equation}
\begin{equation}\label{oa305}
\frac{d}{ds}\int_{-1}^{1}q_{1}q_{2}\rho dy\leq -\frac{1}{10}B(s) +  K_{0}\int_{-1}^{1}q_{-,2}^2\frac{\rho}{1-y^2}dy+K_{0} A(s).
\end{equation}
\item {\bf (Energy barrier)} If $\eqref{111}$ holds, then
\begin{equation}\label{oaTEX}
A(s)\leq K_{1} B(s) +\frac{C_{3}}{s^{\frac{a+1}{2}}}.
\end{equation}
\end{enumerate}
\end{cor}

Coming at this level, we are able to establish the polynomial decay and conclude the proof of Proposition $\ref{Tds}$.
\\{\it {\bf Proof of Proposition $\ref{Tds}$}}: 
Consider $w\in \C([s^*,+\infty ),\H)$ for some $s^*$ large enough a solution of equation $\eqref{1}$ such that $\eqref{111}$ and $\eqref{112}$ hold for some $d^*\in (-1,1)$. Up to replacing $w(y,s)$ by $-w(y,s)$ we may assume that $\omega^*=1$ in $\eqref{ds115}$. Consider then $\epsilon =2KK_{2}\epsilon^*$ where $K$ is given in Proposition $\ref{T7}$ and $K_{2}$ will be fixed later. If
\begin{equation}\label{nouv1}
\epsilon^* \leq \epsilon_{1}\,\,\,{\rm and }\,\,\,\epsilon \leq \epsilon_{4},
\end{equation}
then we see that Proposition $\ref{T7}$, Corollary $\ref{Tnouv}$ and $\eqref{oaMM}$ apply respectively with $\epsilon^*$ and $\epsilon $. In particular, there is a maximal solution $d(s)\in \C^1([s^*,+\infty ),(-1,1))$ such that $\eqref{101}$ holds for all $s\in [s^*,+\infty )$ where $q(y,s)$ is defined in $\eqref{301}$ and 
\begin{equation}\label{nouv2}
|\theta (s^*)-\theta^*|+\Vert q(s^*)\Vert_{\H}\leq K\epsilon^*\,\,\,{\rm with }\,\,\,\theta^*=\frac{1}{2}\log \Big (\frac{1+d^*}{1-d^*}\Big ).
\end{equation}
If in addition we have
\begin{equation}\label{nouv3}
K_{2}\geq 1\,\,\,{\rm hence, }\,\,\,\epsilon \geq 2K\epsilon^*,
\end{equation}
then, we can give two definitions:
\begin{itemize}
\item We define first from $\eqref{nouv2}$ and $\eqref{nouv3}$ $s^*_{1}\in (s^*,+\infty )$ such that for all $s\in [s^*,s_{1}]$,
\begin{equation}\label{nouv4}
\Vert q (s)\Vert_{\H}< \epsilon
\end{equation}
and if $s^*_{1}<+\infty $, then $\Vert q(s^*_{1}) \Vert_{\H}=\epsilon $.
\item Then, we define $s^*_{2}\in [s^*,s^*_{1}]$ as the first $s\in [s^*,s^*_{1}]$ such that
\begin{equation}\label{nouv5}
A(s)\geq \frac{B(s)}{20K_{0}}+\frac{C_{4}}{s^{\frac{a+1}{2}}},
\end{equation}
for some positive $C_{4}$ and $K_{0}$ is introduced in Corollary $\ref{Tnouv}$, or $s^*_{2}=s^*_{1}$ if $\eqref{nouv5}$ is never satisfied on $[s^*,s^*_{1}]$.
  \end{itemize}

Here we outline our formal approach into three steps:
\begin{itemize}
\item In Step 1, using $\eqref{nouv5}$, we integrate equations $\eqref{oa201}$ and $\eqref{oa305}$ on the time interval $[s^*,s^*_{2}]$ and obtain for some positive $K_{3}$, $\eta_{1}$, $C_{5}$, $C_{6}$ and some $f_{0}(s)$
$$\forall \,\,\,s\in [s^*,s^*_{1}]\,\,\frac{1}{K_{3}}\Vert q \Vert^2_{\H}-\frac{C_{5}}{s^{\frac{a+1}{2}}}\leq f_{0}(s)\leq K_{3}\Vert q \Vert^2_{\H}+\frac{C_{5}}{s^{\frac{a+1}{2}}}\,\,\,{\rm and }\,\,\,f'_{0}(s)\leq -\frac{\eta_{1}}{80} f_{0}(s)+\frac{C_{6}}{s^{\frac{a+1}{2}}}.$$
\item In Step 2, integrating equation $\eqref{oa200}$ satisfied by $A(s)$ on the time interval $[s^*,s^*_{2}]$, we obtain some polynomial estimate.
\item In Step 3, we conclude the proof by showing first that $s^*_{1}-s^*_{2}\leq \sigma_{0}$ for some $\sigma_{0}$, then $s^*_{1}=+\infty $. Finally, integrating the equation obtained in Step 1, we conclude the proof of Proposition $\ref{Tds}$.
\end{itemize}

{\bf Step 1:} Integration of the equations on $[s^*,s^*_{2}]$. We claim the following.
\begin{cl}\label{NOUV10}
There exist positive $\epsilon_{5}$, $K_{3}$, $\eta_{1}$, $C_{5}$, $C_{6}$, $C_{7}$ and $f_{0} (s)\in \C^1([s^*,s^*_{2}],\R^+)$ such that if $\epsilon \leq \epsilon_{5}$, then for all $s\in [s^*,s^*_{2}]$,
\begin{enumerate}[{\rm(i)}]
\item  
\begin{equation}\label{mos}
\frac{1}{2}f_{0}(s)-\frac{C_{5}}{s^{\frac{a+1}{2}}}\leq B(s)\leq 2f_{0}(s)+\frac{C_{5}}{s^{\frac{a+1}{2}}}
\end{equation}
and
\begin{equation}\label{mos1}
f'_{0}(s)\leq -\frac{\eta_{1}}{80} f_{0}(s)+\frac{C_{6}}{s^{\frac{a+1}{2}}}.
\end{equation}
\item $$\Vert q (s)\Vert_{\H}\leq K_{3}\Vert q (s^*)\Vert_{\H}e^{-\frac{\eta_{1}}{160}(s-s^*)}+\frac{C_{7}}{s^{\frac{a+1}{4}}}.$$
\end{enumerate}
\end{cl}
{\it Proof}: (i) By definition of $s^*_{2}$, we see that 
 \begin{equation}\label{nouv11}
\forall \,\,\,s\in [s^*,s^*_{2}],\,\,\,A(s)\leq \frac{B(s)}{20K_{0}}+\frac{C_{4}}{s^{\frac{a+1}{2}}},
\end{equation}
where $A(s)$ and $B(s)$ are defined in $\eqref{oaM}$. Since $[s^*,s^*_{2}]\subset [s^*,s^*_{1}]$, the interval where $\eqref{nouv4}$ is satisfied, we can apply Corollary $\ref{Tnouv}$. Therefore, using equations $\eqref{oa201}$ and $\eqref{oa305}$, we write for all $s\in [s^*,s^*_{2}]$
\begin{equation}\label{nouv12}
B'(s)\leq -\frac{8}{p-1}\int_{-1}^{1}q_{-,2}^2\frac{\rho}{1-y^2}dy+C \epsilon B(s)+\frac{C_{2}}{s^{\frac{a+1}{2}}},
\end{equation}
\begin{equation}\label{nouv13}
\frac{d}{ds}\int_{-1}^{1}q_{1}q_{2}\rho dy\leq -\frac{1}{20}B(s) +  K_{0}\int_{-1}^{1}q_{-,2}^2\frac{\rho}{1-y^2}dy+\frac{C_{4}K_{0}}{s^{\frac{a+1}{2}}},
\end{equation}
for some $C>0$ depends on $K_{0}$ and $\epsilon $ small enough. We claim that
\begin{equation}\label{nouv14}
f_{0}(s)=B(s)+\eta_{1} \int_{-1}^{1}q_{1}q_{2}\rho dy
\end{equation}
satisfies all the desired property, where $\eta_{1}>0$ will be fixed small independant of $\epsilon $. Using $\eqref{oa301}$ and $\eqref{nouv11}$, we see that if $\eta_{1}$ is small enough, then we get for all $s\in [s^*,s^*_{2}]$,
\begin{equation}\label{nouv15}
\frac{1}{2}B(s)-\frac{C_{4}\eta_{1}K_{0}}{s^{\frac{a+1}{2}}}\leq f_{0}(s)\leq 2B(s)+\frac{C_{4}\eta_{1}K_{0}}{s^{\frac{a+1}{2}}}.
\end{equation}
Using $\eqref{nouv11}$ and the equivalence of norms $\eqref{oa300}$, we obtain for some $C>0$
\begin{equation}\label{nouv16}
\frac{1}{C}\Vert q(s)\Vert^2_{\H}-\frac{C_{4}\eta_{1}K_{0}}{s^{\frac{a+1}{2}}}\leq f_{0}(s)\leq C\Vert q(s)\Vert^2_{\H}+\frac{C_{4}\eta_{1}K_{0}}{s^{\frac{a+1}{2}}}.
\end{equation}
We combine $\eqref{nouv15}$ and $\eqref{nouv16}$ to conclude the proof of $\eqref{mos}$. Then using $\eqref{nouv12}$, $\eqref{nouv13}$ and $\eqref{nouv14}$, we have for all $s\in [s^*,s^*_{2}]$,
\begin{eqnarray} \label{nouv17}
f'_{0}(s) &\leq &-\Big(\frac{1}{20}\eta_{1}-C \epsilon \Big )B(s) -\Big(\frac{8}{p-1}-K_{0}\eta_{1}\Big )     \int_{-1}^{1}q_{-,2}^2\frac{\rho}{1-y^2}dy +\frac{C_{4}K_{0}}{s^{\frac{a+1}{2}}} \nonumber\\
&\leq &-\frac{\eta_{1}}{40}B(s)+\frac{C_{4}K_{0}}{s^{\frac{a+1}{2}}},
\end{eqnarray}
where $\eta_{1}$ is small enough independant of $\epsilon$ and $\epsilon$ is choosen small enough. Using $\eqref{nouv2}$, $\eqref{nouv15}$ and $\eqref{nouv17}$, to write 
\begin{equation}\label{bint70}
f'_{0}(s)\leq -\frac{\eta_{1}}{80}f_{0}(s)+\frac{C_{4}K_{0}}{s^{\frac{a+1}{2}}}
\end{equation}
and deduce the proof of $\eqref{mos1}$ which concludes the proof of item (i).
\\(ii) Multyplying inequality $\eqref{bint70}$ by $e^{\frac{\eta_{1}}{80}(s-s^*)}$, we obtain the following:
\begin{equation}\label{bint}
\Big(e^{\frac{\eta_{1}}{80}(s-s^*)}f_{0}(s)\Big)'\leq \frac{C_{4}K_{0}e^{\frac{\eta_{1}}{80}(s-s^*)}}{s^{\frac{a+1}{2}}}.
\end{equation}
We integrate $\eqref{bint}$ between $s^*$ and $s$ to obtain
\begin{equation}\label{bint1}
f_{0}(s)\leq e^{-\frac{\eta_{1}}{80}(s-s^*)}f_{0}(s^*) +C_{4}K_{0}e^{-\frac{\eta_{1}}{80}s} \int_{s^*}^{s}\frac{e^{\frac{\eta_{1}}{80}t}}{t^{\frac{a+1}{2}}}dt.
\end{equation}
 That's implie that there exists $C_{7}>0$ depends on $C_{4}$, $K_{0}$ and $\eta_{1}$ such that
 \begin{equation}\label{bint2}
f_{0}(s)\leq e^{-\frac{\eta_{1}}{80}(s-s^*)}f_{0}(s^*) +\frac{C_{7}}{s^{\frac{a+1}{2}}}.
\end{equation}
 Using $\eqref{nouv16}$ to conclude the proof of Claim $\ref{NOUV10}$.

\Box

{\bf Step 2:} Integration of the equations on $[s^*_{2},s^*_{1}]$. We claim the following.
\begin{cl}\label{NOUV18}
The following items hold,
\begin{enumerate}[{\rm(i)}]
\item There exists $\epsilon_{6}>0$ such that for all $\sigma >0$, there exist $K_{4}(\sigma)>0$ and $C_{8}(\sigma)>0$ such that if $\epsilon\leq \epsilon_{6}$, then
$$\forall\,\,\,s\in [s^*_{2},\min(s^*_{2}+\sigma ,s^*_{1})],\,\,\,\Vert q (s)\Vert_{\H}\leq K_{4} e^{-\frac{\eta_{1}}{160}(s-s^*)}\Vert q(s^*)\Vert_{\H}+\frac{C_{8}}{(s^*_{2})^{\frac{a+1}{4}}}.$$
\item There exist positive $\epsilon_{7}$, $C_{12}$, $C_{13}$ and $C$ such that if $\epsilon\leq \epsilon_{7}$ and $s\geq s^*$, then
\begin{equation}\label{nouv19}
\forall\,\,\,s\in [s^*_{2},s^*_{1}],\,\,\,B(s)\leq (A(s)+\frac{C_{12}}{s^{a}})\Big(20K_{0}e^{-\frac{(s-s^*_{2})}{2}}+2C\epsilon+\frac{C_{13}}{s^{a}}\Big),
\end{equation}
where $K_{0}$ is introduced in Corollary $\ref{Tnouv}$.
\end{enumerate}
\end{cl}
{\it Proof}: (i) Using equations $\eqref{oa200}$ and $\eqref{oa201}$, we see that for all $s\in [s^*_{2},\min(s^*_{2}+\sigma ,s^*_{1})]$,
\begin{equation}\label{bint3}
(A(s)+B(s))'\leq 3(A(s)+B(s)) +\frac{C}{s^a},
\end{equation}
 with $C=\max (C_{1},C_{2})$. Multyplying inequality $\eqref{bint3}$ by $e^{-3(s-s^*_{2})}$, we obtain
\begin{equation}\label{bint4}
\Big(e^{-3(s-s^*_{2})}(A(s)+B(s))\Big)'\leq \frac{Ce^{-3(s-s^*_{2})}}{s^a}.
\end{equation}
We integrate $\eqref{bint4}$ between $s^*_{2}$ and $s$, we write
 \begin{equation}\label{bint5}
A(s)+B(s)\leq e^{3(s-s^*_{2})}(A(s^*_{2})+B(s^*_{2}))+Ce^{3(s-s^*_{2})}\int_{s^*_{2}}^{s}\frac{e^{-3(t-s^*_{2})}}{t^a}dt.
\end{equation}
Thanks to $\eqref{oa300}$ and we use the fact that $\frac{1}{2}\leq \frac{s}{s^*_{2}}\leq 2$,
we can see that there exist $C>0$ such that
\begin{equation}\label{bint6}
\Vert q(s)\Vert_{\H}\leq K_{0}e^{\frac{3}{2}\sigma}\Vert q(s^*_{2})\Vert_{\H}+\frac{Ce^{\frac{3}{2}\sigma}}{(s^*_{2})^{a}}.
\end{equation}
Using $(ii)$ in Claim $\ref{NOUV10}$ with $s=s^*_{2}$ gives the conclusion of $(i)$.
\\(ii) By definition of $s^*_{1}$, $\eqref{nouv4}$ is satisfied for all $s\in [s^*_{2} ,s^*_{1}]$, hence, Corollary $\ref{Tnouv}$ applies and equations $\eqref{oa200}$ and $\eqref{oa201}$ holds.
\\Let us first prove that
\begin{equation}\label{nouv20}
\forall\,\,\,s\in (s^*_{2},s^*_{1}],\,\,\,A(s)\geq \frac{B(s)}{20K_{0}}+\frac{C_{4}}{s^{\frac{a+1}{2}}},
\end{equation}
where $K_{0}$ is introduced in Corollary $\ref{Tnouv}$. We need to assume that $s^*_{2}<s^*_{1}$, otherwise the set $(s^*_{2},s^*_{1}]$ is empty. Let $g_{0}(s)=A(s)- \frac{B(s)}{20K_{0}}-\frac{C_{4}}{s^{\frac{a+1}{2}}}$, where $A(s)$ and $B(s)$ are defined in $\eqref{oaM}$. From equations $\eqref{oa200}$ and $\eqref{oa201}$, we write for some $C>0$ and for all $s\in [s^*_{2},s^*_{1}]$
\begin{equation}\label{nouv21}
A'(s)\geq \frac{3}{2}A(s)-C\epsilon B(s)-\frac{C_{1}}{s^{a}} ,\,\,\,B'(s)\leq C\epsilon (A(s)+B(s))+\frac{C_{2}}{s^{a}},
\end{equation}
then, we derive the function $g_{0}$, we obtain
\begin{eqnarray}\label{bint61}
g'_{0}(s)&=&\Big(A(s)- \frac{B(s)}{20K_{0}}-\frac{C_{4}}{s^{\frac{a+1}{2}}}\Big)'\nonumber\\
&\geq &\frac{3}{2}A(s)-C\epsilon B(s)-\frac{C_{9}}{s^{a}}-\frac{C\epsilon}{20K_{0}} (A(s)+B(s))+\frac{C_{4}(a+1)}{2s^{\frac{a+3}{2}}}\nonumber\\
&\geq & A(s)- \frac{B(s)}{20K_{0}}-\frac{C_{4}}{s^{\frac{a+1}{2}}}+\frac{C_{4}(a+1)}{2s^{\frac{a+3}{2}}}+\frac{C_{4}}{s^{\frac{a+1}{2}}}-\frac{C_{9}}{s^{a}},
\end{eqnarray}
for $\epsilon $ small enough. Note that $s^*$ is large enough, so for all $s\geq s^*$, we have:
\begin{equation}\label{bint60}
\frac{C_{4}}{s^{\frac{a+1}{2}}}-\frac{C_{9}}{s^{a}}=\frac{C_{4}}{s^{\frac{a+1}{2}}}\Big(1-\frac{C_{9}}{C_{4}s^{\frac{a-1}{2}}}\Big )\geq \frac{C_{4}}{2s^{\frac{a+1}{2}}}\geq \frac{C_{4}}{2s^{\frac{a+3}{2}}}.
\end{equation}
In view of $\eqref{bint61}$ and $\eqref{bint60}$, we write 
$$g'_{0}(s)\geq A(s)- \frac{B(s)}{20K_{0}}-\frac{C_{4}}{s^{\frac{a+1}{2}}}+\frac{C_{4}(a+2)}{2s^{\frac{a+3}{2}}}=g_{0}(s)+\frac{C_{4}(a+2)}{2s^{\frac{a+3}{2}}}.$$
Since, by definition of $s^*_{2}$, we have $g_{0}(s^*_{2})\geq 0$ and $\eqref{nouv20}$ follows. 
\\Using $\eqref{nouv20}$ and $\eqref{nouv21}$, we obtain for $\epsilon $ small enough and for all $s\in (s^*_{2},s^*_{1}]$ the following
\begin{equation}\label{nouv22}
A'(s)\geq \frac{3}{2}A(s)-20K_{0}C\epsilon ( A(s)-\frac{C_{4}}{s^{\frac{a+1}{2}}})-\frac{C_{1}}{s^{a}}\geq  A(s)-\frac{C_{1}}{s^{a}},
\end{equation}
The same reasoning in $\eqref{bint6}$ can be applied to write the following
\begin{equation}\label{bint8}
A(s)\geq e^{s-s^*_{2}}A(s^*_{2})-\frac{C_{10}}{s^{a}}.
\end{equation}
If $q(s^*_{2})\equiv 0$, then $w(y,s^*_{2})\equiv \overline{w}_{1}(d(s^*_{2}),y,s^*_{2})$ by $\eqref{301}$ and from the uniquness of solutions to equation $\eqref{1}$, we have $w(y,s)\equiv\overline{w}_{1}(d(s^*_{2}),y,s^*_{2})$ and $q(y,s)\equiv 0$ for all $s\geq s^*_{2}$, hence $A(s)=B(s)=0$ by $\eqref{oa300}$ and $\eqref{nouv19}$ follows trivially.
\\Now, if $q(s^*_{2})\not\equiv 0$, we can define $h(s)=\frac{B(s)}{A(s)+\frac{C_{11}}{s^{a}}}$ for all $s\in (s^*_{2},s^*_{1}]$ and derive from $\eqref{nouv21}$ and $\eqref{nouv22}$ for all $s\in (s^*_{2},s^*_{1}]$ 
\begin{eqnarray*}
h'(s)&=&\frac{B'(s)(A(s)+\frac{C_{11}}{s^{a}})-B(s)(A'(s)-\frac{a C_{11}}{s^{a+1}})}{(A(s)+\frac{C_{11}}{s^{a}})^2}\\
&\leq &\frac{ \Big(C\epsilon (A(s)+B(s))+\frac{C_{2}}{s^a}\Big)(A(s)+\frac{C_{11}}{s^{a}})-B(s)A'(s)+\frac{a C_{11}B(s)}{s^{a+1}}}{(A(s)+\frac{C_{11}}{s^{a}})^2}\\
&\leq & \Big(C\epsilon +\frac{aC}{s^{a+1}e^{s-s^*_{2}}A(s^*_{2})}+\frac{2C}{s^{a}e^{s-s^*_{2}}A(s^*_{2})}-1\Big)h(s)+C\epsilon +\frac{C-C\epsilon}{s^{a}e^{s-s^*_{2}}A(s^*_{2})},
\end{eqnarray*}
for $\epsilon $ small enough. We know that $s^*$ is large enough and $\epsilon $ is small enough, one can check that for all $s\geq s^*$, we have:
\begin{equation}\label{bint50}
h'(s)\leq -\frac{h(s)}{2}+C\epsilon +\frac{C}{2s^{a}e^{s-s^*_{2}}A(s^*_{2})}.
\end{equation}
Integrating $\eqref{bint50}$ between $s^*_{2}$ and $s$ to obtain
\begin{equation}\label{bint51}
h(s)\leq e^{-\frac{s-s^*_{2}}{2}}h(s^*_{2})+2C\epsilon +\frac{C}{s^{a}},
\end{equation}
where $C$ depends on $a$, $s^*_{2} $ and $A(s^*_{2})$, we can write now
$$B(s)\leq (A(s)+\frac{ C_{11}}{s^{a}})\Big( e^{-\frac{s-s^*_{2}}{2}}\frac{B(s^*_{2})}{A(s^*_{2})-\frac{C}{(s^*_{2})^{a}}}+2C\epsilon+\frac{C }{s^{a}}\Big ).$$
Using $\eqref{nouv20}$ and taking $\epsilon $ small enough gives $\eqref{nouv19}$ and concludes the proof of Claim $\ref{NOUV18}$.

\Box

{\bf Step 3:} Conclusion of the proof. We use Step 1 and Step 2 to conclude the proof of Proposition $\ref{Tds}$ here. Let us first fix $\sigma_{0}$ such that
\begin{equation}\label{nouv23}
 20K_{0}e^{-\sigma_{0}}+2C\epsilon+\frac{C_{13} }{\sigma_{0}^{a}}\leq \frac{1}{2K_{1}},
\end{equation}
where $K_{0}$, $K_{1}$ are introduced in Corollary $\ref{Tnouv}$. Then, we impose the condition
\begin{equation}\label{bint80}
\epsilon =2K_{2}K\epsilon^*,\,\,\,\,{\rm where}\,\,\,\,K_{2}=\max(2,K_{3},K_{4}).
\end{equation}
Finally, we fix
$$\epsilon_{0} =\min \Big(1,\epsilon_{i},\,\,\,\,{\rm for}\,\,\,\,i\in \lbrace 1,\,4,\,5,\,6,\,7\rbrace \Big ),\,\,\,C_{14}=\max(C_{0},C_{7},C_{8})$$
and the constants are defined in Proposition $\ref{T7}$, Corollary $\ref{Tnouv}$ and Claims $\ref{NOUV10}$ and $\ref{NOUV18}$.
\\Now if $\epsilon^*\leq \epsilon_{0}$, then Corollary $\ref{Tnouv}$ and Steps 1 and 2 apply. We claim that for all $s\in [s^*,\,s^*_{1}]$,
\begin{eqnarray}\label{nouv26}
\Vert q (s)\Vert_{\H} &\leq & K_{2}\Vert q (s^*)\Vert_{\H}e^{-\frac{\eta_{1}}{160}(s-s^*)}+\frac{C_{14}}{s^{\frac{a+1}{4}}}\nonumber\\
&\leq & K_{2}K\epsilon^*e^{-\frac{\eta_{1}}{160}(s-s^*)}+\frac{C_{14}}{s^{\frac{a+1}{4}}}\leq \frac{\epsilon}{2}e^{-\frac{\eta_{1}}{160}(s-s^*)}+\frac{C_{14}}{s^{\frac{a+1}{4}}}.
\end{eqnarray}
Indeed, if $s\in [s^* ,\min(s^*_{2}+\sigma_{0} ,s^*_{1})]$, then this comes from (ii) of Claim $\ref{NOUV10}$ or (i) of Claim $\ref{NOUV18}$ and the definition $\eqref{bint80}$ of $K_{2}$.
\\Now, if $s^*_{2}+\sigma_{0} <s^*_{1}$ and $s\in [s^*_{2}+\sigma_{0} ,s^*_{1}]$, then we have from $\eqref{nouv19}$ and the definition of $\sigma_{0}$ that $B(s)\leq \frac{1}{2K_{1}} \Big(A(s)+\frac{ C_{11}}{s^{a}}\Big )$ on the one hand. On the other hand, from (iii) of Corollary $\ref{Tnouv}$, we have $A(s)\leq K_{1}B(s)+\frac{C_{3}}{s^{\frac{a+1}{2}}}$, we use $\eqref{oa300}$ to deduce that $\eqref{nouv26}$ is satisfied.
\\In particular, we have for all $s\in [s^* ,s^*_{1}]$, $\Vert q (s)\Vert_{\H}\leq \frac{\epsilon}{2}+\frac{C_{14}}{s^{\frac{a+1}{4}}}$. Note that $s^*$ is large enough, so we can say that 
 for all $s\in [s^*,s^*_{1}] $, we have:
$$\Vert q (s)\Vert_{\H}\leq \frac{3\epsilon}{4},$$
hence, by definition of $s^*_{1}$, this means that $s^*_{1}=+\infty$. Therefore, from $\eqref{nouv26}$ and $\eqref{oaest}$, we have  
\begin{equation}\label{nouv27}
\forall\,\,\,\,s\geq s^*,\,\,\,\,\Vert q (s)\Vert_{\H}\leq \frac{C_{15}}{s^{\frac{a+1}{4}}}\,\,\,\,{\rm and}\,\,\,\,|\theta' (s)| \leq \frac{C_{16}}{s^{\frac{a+1}{2}}} ,
\end{equation}
 Hence, there is $\theta_{\infty}\in \R$ such that $\theta (s)\longrightarrow \theta_{\infty}$ as $s\longrightarrow \infty$ and 
\begin{equation}\label{nouv28}
\forall \,\,\,\,s\geq s^*,\,\,\,\,|\theta (s)-\theta_{\infty}|\leq \frac{C_{17}}{s^{\frac{a+1}{2}}} .
\end{equation}
 Taking $s=s^*$ here and using $\eqref{nouv2}$, we see that $|\theta_{\infty}-\theta^*|\leq C\epsilon^*$.
\\If $d_{\infty}=\tanh \theta_{\infty}$, then we see that $|d_{\infty}-d^*|\leq C(1-d^{*2})\epsilon^*$.
\\Using the definition of $q$ given in $\eqref{301}$, inequalities $\eqref{e10t}$, $\eqref{nouv27}$ and $\eqref{nouv28}$, we write
$$\Big\Vert \left(
  \begin{array}{ccc}
    w(s) \\
    \partial_{s} w(s)\\
  \end{array}
\right)-\omega^*\left(
  \begin{array}{ccc}
    \overline{w}_{1}(d_{\infty},y,s) \\
   \overline{w}_{2}(d_{\infty},y,s)\\
  \end{array}
\right)\Big\Vert_{\H}$$
$$\leq \Big\Vert \left(
  \begin{array}{ccc}
    w(s) \\
    \partial_{s} w(s)\\
  \end{array}
\right)-\omega^*\left(
  \begin{array}{ccc}
    \overline{w}_{1}(d(s),y,s) \\
   \overline{w}_{2}(d(s),y,s)\\
  \end{array}
\right)\Big\Vert_{\H}+\Big\Vert \left(
  \begin{array}{ccc}
    \overline{w}_{1}(d(s),y,s) \\
    \overline{w}_{2}(d(s),y,s)\\
  \end{array}
\right)-\omega^*\left(
  \begin{array}{ccc}
    \overline{w}_{1}(d_{\infty},y,s) \\
   \overline{w}_{2}(d_{\infty},y,s)\\
  \end{array}
\right)\Big\Vert_{\H}$$
$$\leq \Vert q(s)\Vert_{\H} +C|\theta (s)-\theta_{\infty}|\leq  \frac{C_{18}}{s^{\frac{a+1}{4}}}.$$
 This concludes the proof of Proposition $\ref{Tds}$.

\Box

\bigskip

As a consequence of our polynomial decay obtained in Proposition $\ref{Tds}$, we can deduce that for $d(s)\in \C^1([s^*,+\infty ),(-1,1))$ introduced in $\eqref{101}$, we have 
$$\Vert q(s)\Vert_{\H}\longrightarrow 0\,\,\,{\rm as}\,\,\, s\longrightarrow +\infty .$$
 This crucial information helps us to obtain in the next section the desired exponential decay.
\subsection{Proof of Theorem $\ref{T4}$}
In this subsection we prove Theorem $\ref{T4}$. We should keep in mind that from Proposition $\ref{Tds}$, we have $\Vert q(s)\Vert_{\H}\longrightarrow 0$ as $s\longrightarrow +\infty$ for the same $d(s)$ introduced in $\eqref{101}$. Following this information, we can see that 
\begin{equation}\label{oa400}
A(s)+B(s)\longrightarrow 0 \,\,\,\,\,\,{\rm as}\,\,\,\,\,\,s\longrightarrow \infty ,
\end{equation}
where $A(s)$ and $B(s)$ are defined in $\eqref{oaM}$. We first show that $A(s)$ is controlled by $ B(s)$, which is not our goal but this control and suitable refinements of some results obtained in the previous part allow us to find the exponential decay and conclude the proof of our Theorem $\ref{T4}$. We start by the following lemma.
\begin{lem}\label{Tzin}
For all $s\geq s^*$, we have:
$$A(s)\leq \frac{1}{4} B(s).$$
\end{lem}

{\it Proof}: We define for all $s\geq s^*$ the function:
$$\gamma (s)=A(s)-\frac{1}{4} B(s).$$
A direct consequence of $\eqref{306}$ and $\eqref{310}$, we can choose $\epsilon $ small enough in Corollary $\ref{Tnouv}$ and $s^*$ large such that for all $s\geq s^*$, the following estimates hold:
$$A'(s)\geq \frac{1}{2}A(s) -\frac{3}{64} B(s) ,$$
$$ B'(s)\leq  A(s)+\frac{1}{16} B(s).$$
We can see that
$$\gamma'(s)=A'(s)-\frac{1}{4} B'(s)\geq \frac{1}{4}\gamma (s),$$
since $\gamma (s)\longrightarrow 0$ as $s\longrightarrow \infty $ (see $\eqref{oa400}$), implies $\gamma (s)\leq 0$, which conclude the proof of Lemma $\ref{Tzin}$.

\Box

This way, we are in a position to perform some refinements of the estimates on the function $f_{0}$ defined in $\eqref{nouv14}$ and some other estimates obtained in the previous subsection. More precisely, we have the following:
\begin{cl}\label{T13}
For all $s\geq s^*$, we have 
\begin{equation}\label{oa410}
\frac{1}{2}f_{0}(s)\leq B(s)\leq 2f_{0}(s),
\end{equation}
\begin{equation}\label{oa600}
f'_{0}(s)\leq -\frac{\eta_{1}}{8}f_{0}(s),
\end{equation}
\begin{equation}\label{oa440}
\Vert q(s)\Vert_{\H} \leq C\Vert q(s^*)\Vert_{\H}e^{-\frac{\eta_{1}}{16}(s-s^*)} \leq \tilde{K}  e^{-\frac{\eta_{1}}{16}(s-s^*)},
\end{equation}
\begin{equation}\label{oaams}
|\theta' (s)|\leq \tilde{K}^2 e^{-\frac{\eta_{1}}{8}(s-s^*)},
\end{equation}
such that $\tilde{K}$ depends on $\epsilon $, the constant $K_{0}$ and the function $f_{0}$ defined respectively in $\eqref{nouv14}$ and Corollary $\ref{Tnouv}$.
\end{cl}

\bigskip

{\it Proof}:
Using $\eqref{oa301}$, $\eqref{nouv14}$ and Lemma $\ref{Tzin}$, we see that if $\eta_{1}$ is small enough, then we get 
\begin{equation}\label{oa421}
\frac{1}{2}B(s)\leq  f_{0}(s)   \leq 2B(s).
\end{equation}
We use again Lemma $\ref{Tzin}$, estimate $\eqref{oa421}$ and the equivalence of norms $\eqref{oa300}$, we obtain for some $C>0$ 
\begin{equation}\label{oa420}
\frac{1}{C}\Vert q(s)\Vert^2_{\H}\leq  f_{0}(s)   \leq C\Vert q(s)\Vert^2_{\H}.
\end{equation}
Using $\eqref{310}$, Corollary $\ref{Tnouv}$ and Lemma $\ref{Tzin}$, we get for some $C>0$ and for all $s\geq s^*$
\begin{equation}\label{oa402}
B'(s)\leq -\frac{8}{p-1}\int_{-1}^{1}q_{-,2}^2\frac{\rho}{1-y^2}dy+\frac{5K_{0}\epsilon}{4} B(s)+\frac{C}{s^{a}}B(s),
\end{equation}
\begin{equation}\label{oa403}
\frac{d}{ds}\int_{-1}^{1}q_{1}q_{2}\rho dy\leq -\frac{1}{20}B(s) +  K_{0}\int_{-1}^{1}q_{-,2}^2\frac{\rho}{1-y^2}dy.
\end{equation}
Then using $\eqref{oa402}$, $\eqref{oa403}$, taking $\epsilon $ small enough and $s^*$ large enough, we obtain for all $s\geq s^*$
\begin{eqnarray}\label{oa404}
f'_{0}(s) &\leq &-(\frac{1}{20}\eta_{1} -\frac{5K_{0}\epsilon}{4}-\frac{C}{s^{a}} )B (s)-(\frac{8}{p-1}-  K_{0}\eta_{1} )\int_{-1}^{1}q_{-,2}^2\frac{\rho}{1-y^2}dy \leq -\frac{\eta_{1}}{4}B(s)\nonumber\\
&&\leq -\frac{\eta_{1}}{8}f_{0}(s),
\end{eqnarray}
which conclude the proof of $\eqref{oa600}$.
\\Integrating $\eqref{oa404}$, we get for all $s\geq s^*$, $f_{0}(s)\leq f_{0}(s^*)e^{-\frac{\eta_{1}}{8}(s-s^*)}$. Using $\eqref{oa420}$ and the a priori estimate $\eqref{200}$, we can deduce that $\eqref{oa440}$ follows.
\\The proof of $\eqref{oaams}$ follows directly from $\eqref{oa440}$. Which conclude the proof of Claim $\ref{T13}$.

\Box

At this level, we are ready to adapt the proof of the case when $f\equiv 0$ and deduce our main result in Theorem $\ref{T4}$. We are going to give the deduction of the proof:
\\From $\eqref{oaams}$, we can see that there exists $\theta_{\infty}\in \R$ such that $\theta (s)\longrightarrow \theta_{\infty}$ as $s\longrightarrow \infty$ and 
$$\forall \,\,\,\,s\geq s^*,\,\,\,\,|\theta (s)-\theta_{\infty}|\leq \tilde{K}^2 e^{-\frac{\eta_{1}}{8}(s-s^*)}.$$
Using the definition of $q$ given in $\eqref{301}$, inequalities $\eqref{e10t}$, $\eqref{e15t}$, $\eqref{oa440}$ and $\eqref{oaams}$, we write:
$$\Big\Vert \left(
  \begin{array}{ccc}
    w(s) \\
    \partial_{s} w(s)\\
  \end{array}
\right)-\omega^*\left(
  \begin{array}{ccc}
    \overline{w}_{1}(d_{\infty},y,s) \\
   \overline{w}_{2}(d_{\infty},y,s)\\
  \end{array}
\right)\Big\Vert_{\H}$$
$$\leq \Big\Vert \left(
  \begin{array}{ccc}
    w(s) \\
    \partial_{s} w(s)\\
  \end{array}
\right)-\omega^*\left(
  \begin{array}{ccc}
    \overline{w}_{1}(d(s),y,s) \\
   \overline{w}_{2}(d(s),y,s)\\
  \end{array}
\right)\Big\Vert_{\H}+\Big\Vert \left(
  \begin{array}{ccc}
    \overline{w}_{1}(d(s),y,s) \\
    \overline{w}_{2}(d(s),y,s)\\
  \end{array}
\right)-\omega^*\left(
  \begin{array}{ccc}
    \overline{w}_{1}(d_{\infty},y,s) \\
   \overline{w}_{2}(d_{\infty},y,s)\\
  \end{array}
\right)\Big\Vert_{\H}$$
$$\leq \Vert q(s)\Vert_{\H} +C|\theta (s)-\theta_{\infty}|\leq \tilde{K_{0}} e^{-\frac{\eta_{1}}{16}(s-s^*)},$$
where $\tilde{K_{0}}$ depends on $\tilde{K}$. This concludes the proof of Theorem $\ref{T4}$.

\bigskip

{\bf Acknowledgment:} The authors wish to thank Professor Hatem ZAAG for many fruitful discussions, valuable suggestions and guidance in this work. Part of this work was done when the second author was visiting the Laboratoire Analyse Géométrie et Applications (LAGA) of university Sorbonne Paris Nord. He is grateful to LAGA for the hospitality and the stimulating atmosphere.

\appendix{}
\begin{section}{Construction of a particular solution of equation $\eqref{1}$ in similarity variables}\label{TS2}
\end{section}
As already written in the introduction, our aim in this section is to show the detail of the construction of the solution defined in $\eqref{87}$ and $\eqref{e505}$. We start by the case where $d=0$ and the deduction of the case where $d\neq 0$ will follow from the fact that $\eqref{y}$ is invariant under the Lorentz transform.
\\{\bf Case} $d=0$. Note that $\kappa (d,y)=\kappa_{0}$ is not a solution of equation $\eqref{1}$. In order to find $\phi (s)$ a solution of equation $\eqref{1}$ such that $\phi (s)\longrightarrow \kappa_{0}$ as $s\longrightarrow +\infty$, it is equivalent to show a solution $\varphi_{T} (t)$ of equation $\eqref{y}$ such that $\varphi_{T} (t)$ blows-up at time $T$ with $\varphi_{T} (t)\sim \kappa_{0} (T-t)^{-\frac{2}{p-1}}$ as $t\longrightarrow T$. 
\\{\bf Case} $d\neq 0$. We know that equation $\eqref{y}$ is invariant under the Lorentz transform. So, for any $d \in (-1,1)$, if we define $\overline{\varphi}(d,x,t)=\varphi_{0}(\frac{t-dx}{\sqrt{1-d^2}})$, then $\overline{\varphi}$ is a solution of equation $\eqref{y}$ which blows-up at the line $\lbrace t=dx \rbrace $. Let us transform it with the definition of $w$:
\begin{equation} \label{3}
y=-\frac{x}{t},\,\,\,\,\,\,s=-\log(-t)\,\,{\rm and }\,\, \overline{w}_{1}(d,y,s)=(-t)^\frac{2}{p-1}\overline{\varphi}(d,x,t).\,\,\,\,\,\,\,\,\,\,\,\\
\end{equation}
Some simple computation gives:
$$\overline{w}_{1}(d,y,s)=e^{\frac{-2s}{p-1}}\varphi_{0}\Big(-e^{-s}\frac{1+dy}{\sqrt{1-d^2}}\Big).$$
In other words, we look in the following lemma for $\varphi_{T} (t)\in \C^2([-t_{0},T),\R)$ a solution of $\varphi"=\varphi^{p}+f(\varphi )$ with $\varphi_{T} (t)\sim \kappa_{0} (T-t)^{-\frac{2}{p-1}}$ as $t\rightarrow T$ and $f$ defined in $\eqref{fy}$.
\begin{lem}\label{T5}
The following items hold,
\begin{enumerate}[{\rm(i)}]
\item There exists $\varphi_{0} (t)\in \C^2([-t_{0},0),\R)$ solution of
$\varphi"(t)=\varphi^{p}(t)+f(\varphi (t) ),$ with $\varphi_{0} (t)\sim \kappa_{0} (-t)^{-\frac{2}{p-1}}$ as $t\longrightarrow 0$.
\item If $\varphi_{T} (t)=\varphi_{0} (t-T)$, then $\varphi_{T}"(t)=\varphi_{T}^{p}(t)+f(\varphi_{T} (t))$, with $\varphi_{T} (t)\sim \kappa_{0} (T-t)^{-\frac{2}{p-1}}$ as $t\longrightarrow T$.
\item $\phi (s)=e^{\frac{-2s}{p-1}}\varphi_{T}(T-e^{-s})=e^{\frac{-2s}{p-1}}\varphi_{0}(-e^{-s}).$
\end{enumerate}
\end{lem}
{\it Proof}: We have the following associated {\bf ODE} to the {\bf PDE} $\eqref{y}$.
\begin{equation}\label{ode}
\left\{
  \begin{array}{ll}
  \varphi"(t)=\varphi^{p}(t)+f(\varphi(t) ),\\
  (\varphi(0),\varphi'(0))=(A,B),\,\,\,\,\,\,\,\,\,\,\,\,\,\,\,\,\,\,\,\,\,\,\,\\
  \end{array}
\right.
\end{equation}
with $A>0$ and $B=B(A)>0$ such that $B^2-\frac{A^{p+1}}{p+1}\geq 0.$ By the Cauchy theory, there exists a maximal solution $\varphi_{A} (t)$ defined in $[0,T_{A})$, with $T_{A}\leq +\infty$. In order to conclude the proof of Lemma $\ref{T5}$, we proceed in 3 steps.
\\{\bf Step 1:} In this step we prove that for $A$ large enough and for all $t\in [0,T_{A})$, we have $\varphi'_{A}(t)\geq 0.$ Let $A_{0}$ large enough, such that 
$$\forall \,\,\,\xi \geq A_{0},\,\,\,|f(\xi)|\leq \frac{\xi^{p}}{2},$$
we consider $A\geq A_{0}$. Since, we have $\varphi_{A}"(t)\geq \frac{\varphi_{A}^{p}(t)}{2},$ whenever $\varphi_{A}(t)\geq A$, using a contradiction argument, the result follows.
\\{\bf Step 2:} The solution $\varphi_{A} (t)$ is an increasing function, so for all $t\in [0,T_{A})$, we have $\varphi_{A} (t)\geq \varphi_{A} (0)= A> 0$. 
\\Recalling that $\varphi'_{A}(t)\geq 0$, so, by multiplying inequality $\varphi_{A}"(t)\geq \frac{\varphi_{A}^{p}(t)}{2}$ by $\varphi'(t)$, we obtain 
$$\frac{d}{dt}\Xi (\varphi_{A} (t))\geq 0,$$
with $\Xi (\varphi_{A}(t))=(\varphi_{A}'(t))^2-\frac{\varphi_{A}^{p+1}(t)}{p+1}$. From the fact that $B^2\geq \frac{A^{p+1}}{p+1} $, we obtain
$$\Xi (\varphi_{A}(t))\geq \Xi (\varphi_{A} (0))\geq 0.$$
As mentioned above, $\varphi_{A}(t)> 0$ for all $t\in [0,T_{A})$, so 
$$\frac{\varphi_{A}'(t)}{\varphi_{A}^{\frac{p+1}{2}}(t)}\geq \frac{1}{\sqrt{p+1}}.$$
After integration between $0$ and $T_{A}$, we deduce that $\varphi_{A} (t)$ blow's up in finite time $T_{A}$.
\\{\bf Step 3:} If $\varphi(t)$ is a solution of equation $\eqref{ode}$ which blow's up in finite time $T$ (for example $\varphi(t) =\varphi_{A}(t)$ constructed in Step 2).
If $\epsilon > 0$, then there exists $t_{0}(\epsilon)$ such that $\forall \,\,\,t\in [t_{0}(\epsilon), T)$, we have:
$$|f(\varphi_{A} (t))|\leq \epsilon \varphi_{A}^{p}(t).$$
Next, we multiplie $\eqref{ode}$ by $\varphi'$ and we integrate between $0$ and $t$ in order to get the following:
\begin{equation} \label{odee}
 2\frac{\varphi_{A}^{p+1}(t)(1-\epsilon )}{p+1}-C \leq\varphi_{A}'^{2}(t)\leq 2\frac{\varphi_{A}^{p+1}(t)(1+\epsilon )}{p+1}+C,
\end{equation}
with $C\leq \epsilon \frac{\varphi_{A}^{p+1}}{p+1}.$ Solving $\eqref{odee}$, yields to  
$$\kappa_{0}(T_{A}-t)^{-\frac{2}{p-1}}(1-2\epsilon )^{-\frac{1}{p-1}}\leq\varphi_{A}(t)\leq \kappa_{0}(T_{A}-t)^{-\frac{2}{p-1}}(1+2\epsilon )^{-\frac{1}{p-1}},$$
where $\kappa_{0}$ is introduced in $\eqref{400}$. Since $\epsilon>0$ is arbitrary, we can find that 
$$\varphi_{A}(t)\sim \kappa_{0}(T-t)^{-\frac{2}{p-1}},$$ 
which concludes the proof of Lemma $\ref{T5}$.

\Box

\begin{nb}
Combining $(i)$ and $(iii)$ of Lemma $\ref{T5}$, we see that $\phi (s)$ is a solution of equation $\eqref{1}$ with 
\begin{equation}\label{fek}
\phi (s)\longrightarrow \kappa_{0}\,\,\,\,{\rm as}\,\,\,\,  s\longrightarrow +\infty.
\end{equation}

\end{nb}

\bigskip

\begin{nb}
In particular, when $d=0$, we show the following equivalence 
$$\phi (s)=e^{\frac{-2s}{p-1}}\varphi_{0}(-e^{-s})\Leftrightarrow \varphi_{0}(t)= \phi (s)(-t)^{-\frac{2}{p-1}}.$$
We now focus on the case when $d\neq 0$, we exploite the expression $\eqref{400}$ of $\kappa (d,y)$ to write the following
\begin{eqnarray*}
\overline{w}_{1}(d,y,s)&=& e^{\frac{-2s}{p-1}}\varphi_{0}\Big(-e^{-s}\frac{1+dy}{\sqrt{1-d^2}}\Big)=e^{\frac{-2s}{p-1}}\Big(e^{-s}\frac{1+dy}{\sqrt{1-d^2}}\Big)^{\frac{-2}{p-1}}\phi (-\log(e^{-s}\frac{1+dy}{\sqrt{1-d^2}}))\\
\overline{w}_{1}(d,y,s)&=&\frac{(1-d^2)^{\frac{1}{p-1}}}{(1+dy)^{\frac{2}{p-1}}}\phi (s-\log(\frac{1+dy}{\sqrt{1-d^2}}))=\kappa (d,y)\frac{\phi (s-\log(\frac{1+dy}{\sqrt{1-d^2}}))}{\kappa_{0}}.
\end{eqnarray*}
Finally, we find our solution defined in $\eqref{87}$ and $\eqref{e505}$.
\end{nb}
\begin{nb}
In the case when $f\equiv 0$, $\varphi_{0}=\kappa_{0}(-t)^{-\frac{2}{p-1}}$ and $\overline{\varphi}(d,x,t)=\kappa_{0}\Big(\frac{-t+dx}{\sqrt{1-d^2}}\Big)^{-\frac{2}{p-1}}$, we can see that 
$$\overline{w}_{1}(d,y,s)=\kappa (d,y).$$
In this case, $\varphi_{0}\sim \kappa_{0}(-t)^{-\frac{2}{p-1}}$, therefore 
\begin{equation}\label{fer}
\overline{w}_{1}(d,y,s)\sim\kappa (d,y)\,\,\,\,{\ as}\,\,\,\, s\longrightarrow \infty .
\end{equation}
\end{nb}
\section{Property of the particular solution $\phi(s)$ }\label{T160}
The following claim shows the asymptotic behavior of the particular solution $\phi(s)$ of equation $\eqref{1}$, which is crucial in many steps in this paper.
\begin{cl}\label{T14}{\bf (Equivalent to $\phi(s)-\kappa_{0}$)} For all $p>1$ and $a>1$, we have, 
\begin{enumerate}[{\rm(i)}]
\item $\phi(s)-\kappa_{0} \sim  -\frac{\kappa_{0}}{p-1}\Big(\frac{p-1}{4s}\Big)^{a}\,\,\,{\rm as}\,\,s \rightarrow +\infty ,$
\item $|\phi'(s)|\leq \frac{C}{s^a},$
\end{enumerate}
where $\phi(s)$ satisfies equation $\eqref{1}$, the constant $\kappa_{0}$ is defined in $\eqref{400}$ and $C>0$. 
\end{cl}
{\it Proof}: The proof of item $(ii)$ is a direct consequence from the proof of $(i)$ and we use the fact that $\phi(s)$ satisfies equation $\eqref{1}$. We are going now to give the proof of item $(i)$.
\\We proceed in three steps:
\\{\bf Step 1:}
Let us recall the following equation:
\begin{equation}\label{73}
\varphi_{0}''(t)=\varphi_{0}^{p}(t)+f(\varphi_{0}(t) ),
\end{equation}
and also from the previous part $\varphi_{0}(t)\sim \kappa_{0}(-t)^{-\frac{2}{p-1}}$, $\varphi'_{0}(t)>0$ and $\varphi_{0}''(t)>0$.
 \\Multiplying equation $\eqref{73}$ by $\varphi'_{0}(t)$ and integrating between $0$ and $t$, we can see that 
 \begin{equation}\label{74}
 \varphi'_{0}(t)=\sqrt{\frac{2\varphi_{0}^{p+1}(t)}{p+1}+2F(\varphi_{0}(t))+M_{0}},
\end{equation}
where $M_{0}=M_{0}\Big(\varphi'_{0}(0),\varphi_{0}(0)\Big)$ and $F$ is defined in $\eqref{90}$.
\\{\bf Step 2:}
We use now the following self-similar change of variables:
\begin{equation} \label{75}
\phi(s)=(-t)^\frac{2}{p-1}\varphi_{0}(t)\,\,\,{\rm and }\,\,\,\,\,s=-\log(-t).
\end{equation}
Some computation gives 
\begin{equation} \label{76}
\phi'(s)=-\frac{2}{p-1}\phi(s)+\sqrt{\frac{2\phi^{p+1}(s)}{p+1}+2e^{\frac{-2(p+1)s}{p-1}}F(e^{\frac{2s}{p-1}}\phi(s))+e^{\frac{-2(p+1)s}{p-1}}M_{0}}.
\end{equation}
{\bf Step 3:}
According to Step 1 and Step 2, we may try to find an equivalent to $\phi(s)-\kappa_{0}$.
\\If we note by $q=\phi(s)-\kappa_{0}$, according to equation $\eqref{76}$, we can see that $q(s)$ satisfies the following equation 
\begin{equation} \label{78}
\hspace{-0.1cm} q'(s)\hspace{-0.1cm}=\hspace{-0.1cm}-\frac{2(q(s)+\kappa_{0})}{p-1}+\sqrt{\frac{2(q(s)+\kappa_{0})^{p+1}}{p+1}+e^{\frac{-2(p+1)s}{p-1}}(2F(e^{\frac{2s}{p-1}}(q(s)+\kappa_{0}))+M_{0})}.
\end{equation}
Using Taylor expansion to derive formally from $\eqref{fek}$ the following
\begin{equation} \label{77}
 (q(s)+\kappa_{0})^{p+1}=\kappa_{0}^{p+1}+(p+1)\kappa_{0}^{p}q(s)+ \O(q^2(s))    .
\end{equation}
An integration by part gives when $\theta\rightarrow +\infty $
\begin{equation} \label{79}
 F(\theta )= \frac{|\theta|^{p+1}}{(p+1)\log^{a}(2+\theta^2)}+\O\Big(\frac{|\theta|^{p+1}}{\log^{a+1}(2+\theta^2)}\Big). 
\end{equation}
According to $\eqref{77}$ and $\eqref{79}$, we can write 
\begin{equation}\label{80}
F(e^{\frac{2s}{p-1}}(q(s)+\kappa_{0}) )=\frac{e^{\frac{2(p+1)s}{p-1}}\kappa_{0}^{p+1}}{(p+1)\Big(\frac{4s}{p-1}\Big)^{a}}+\mathit{o} (\frac{e^{\frac{2(p+1)s}{p-1}}}{s^{a}}).
\end{equation}
Combining equation $\eqref{78}$, $\eqref{77}$ and $\eqref{80}$, we write 
\begin{equation} \label{it}
\hspace{-0.1cm}q'(s)\hspace{-0.1cm} = -\frac{2\kappa_{0}}{p-1}-\frac{2q(s)}{p-1}+\hspace{-0.2cm}\sqrt{\frac{2}{p+1}[\kappa_{0}^{p+1}\hspace{-0.1cm}+\hspace{-0.1cm}(p+1)\kappa_{0}^{p}q(s)\hspace{-0.1cm}+\hspace{-0.1cm}\frac{\kappa_{0}^{p+1}}{2\Big(\frac{4s}{p-1}\Big)^{a}}]\hspace{-0.1cm}+ \hspace{-0.1cm}\O(q^2(s))\hspace{-0.1cm}+\hspace{-0.1cm}\mathit{o} (\frac{1}{s^{a}})}.
\end{equation}
From $\eqref{it}$ and the fact that $\sqrt{\frac{2\kappa_{0}^{p+1}}{p+1}}=\frac{2\kappa_{0}}{p-1}$, it is simple to write the following
$$q'(s) = -\frac{2\kappa_{0}}{p-1}-\frac{2q(s)}{p-1}+\frac{2\kappa_{0}}{p-1}\sqrt{1+\frac{p+1}{\kappa_{0}}q(s)+\Big (\frac{p-1}{4s}\Big )^a+\O(q^2(s))+\mathit{o} (\frac{1}{s^{a}})}.$$
Using Taylor expansion of the function $Z\longrightarrow \sqrt{1+Z}$ as $Z\longrightarrow 0$, applying at $Z=\frac{p+1}{\kappa_{0}}q(s)+\Big (\frac{p-1}{4s}\Big )^a+\O(q^2(s))+\mathit{o} (\frac{1}{s^{a}})$, we derive
\begin{equation} \label{81}
\hspace{-0.1cm}q'(s) = -\frac{2\kappa_{0}}{p-1}-\frac{2q(s)}{p-1}+\frac{2\kappa_{0}}{p-1}\Big(1\hspace{-0.1cm}+\frac{p+1}{2\kappa_{0}}q(s)+\frac{1}{2}\Big(\frac{p-1}{4s}\Big)^{a}+\O (q^2(s))+\hspace{-0.1cm}\mathit{o} (\frac{1}{s^{a}})\Big).
\end{equation}
Finally, by $\eqref{81}$ it is obvious to write
\begin{equation} \label{83}
q'(s)=q(s)+\frac{\kappa_{0}}{p-1}\Big(\frac{p-1}{4s}\Big)^{a}+\O (q^2(s))+\mathit{o} (\frac{1}{s^{a}}) .
\end{equation}
We would like now to prove that $q(s)=\phi(s)-\kappa_{0} \sim  -\frac{c_{p}^{a}}{s^{a}},$ where $c_{p}^{a}=\frac{\kappa_{0}(p-1)^{a-1}}{4^{a}}$. To do that let $g_{-}(s)=-\frac{1}{s^{a-\frac{1}{2}}},$ $g_{+}(s)=\frac{1}{s^{\frac{1}{2}}}$ and $g(s)=-\frac{c_{p}^{a}}{s^{a}}$. We start by remarking that the flow is transverse outgoing in the curve of $g_{+}$ and $g_{-}$. We explain in the following this fact.
Firstly, we start by the case when $q(s)=g_{-}(s)$. A simple derivation gives $g'_{-} (s)=\frac{a-\frac{1}{2}}{s^{a+\frac{1}{2}}}>0.$ We exploite $\eqref{83}$ to write when $s$ is large enough
\begin{equation} \label{sai5}
q'(s)=-\frac{1}{s^{a-\frac{1}{2}}}+\frac{c_{p}^{a}}{s^{a}}+\O (\frac{1}{s^{2a-1}})+\mathit{o} (\frac{1}{s^{a}})<0 .
\end{equation}
We conclude that $q'(s)<g'_{-} (s)$ for $s$ large enough and the flow is transverse outgoing in the curve of $g_{-}$.
Now, we treat the case when $q(s)=g_{+}(s)$. A simple derivation gives $g'_{+} (s)=-\frac{1}{2s^{\frac{3}{2}}}<0.$ We exploite $\eqref{83}$ to write when $s$ is large enough
\begin{equation} \label{sai6}
q'(s)=\frac{1}{s^{\frac{1}{2}}}+\frac{c_{p}^{a}}{s^{a}}+\O (\frac{1}{s})+\mathit{o} (\frac{1}{s^{a}})>0 .
\end{equation}
We conclude that $q'(s)>g'_{+} (s)$ for $s$ large enough and the flow is transverse outgoing in the curve of $g_{+}$. The conclusion of this part is that the flow is transverse outgoing in the curve of $g_{+}$ and $g_{-}$. Three cases then arise:
\\{\bf The first case:}\underline{ If for all $s$ large enough, we have $g_{-}(s)\leq q(s) \leq g_{+}(s).$} 
\\Since, $q(s)\longrightarrow 0$ as $s\longrightarrow +\infty$, we rewrite $\eqref{83}$ as follow 
\begin{equation} \label{sai}
q'(s)=q(s)+\frac{c_{p}^{a}}{s^{a}}+\O (q^2(s))+\mathit{o} (\frac{1}{s^{a}}) .
\end{equation}
In this case we have $\O (q^2(s))+\mathit{o} (\frac{1}{s^{a}})=\mathit{o} (\frac{1}{s^{a}})$, which implie that we can write 
$$\O (q^2(s))+\mathit{o} (\frac{1}{s^{a}})=\frac{\varepsilon (s)}{s^{a}}\,\,\,{\rm with }\,\,\,\varepsilon (s)\longrightarrow 0\,{\rm as }\,s\longrightarrow +\infty.$$
Multyplying equation $\eqref{sai}$ by $e^{-\tau}$ and we integrate between $s$ and $+\infty$, we obtain the following
\begin{equation} \label{sai1}
\hspace{-0.1cm}-e^{-s}q(s)=\int_{s}^{+\infty}  \frac{c_{p}^{a}e^{-\tau}}{\tau^{a}}d\tau+\int_{s}^{+\infty}\frac{\varepsilon (\tau)e^{-\tau}}{\tau^{a}} d\tau,\,\,\,{\rm with }\,\,\,\varepsilon (s)\longrightarrow 0 \,{\rm as }\,s\longrightarrow +\infty,
\end{equation}
with equation $\eqref{sai1}$, we can deduce that
$$q(s)\sim  -\frac{c_{p}^{a}}{s^{a}},\,\,\,{\rm as }\,\,\,s\longrightarrow +\infty.$$
\\{\bf The second case:} \underline{If there exist $\tilde{s}$ large enough such that $ q(\tilde{s}) < g_{-}(\tilde{s})$.}
\\In this case, for all $s\geq \tilde{s}$, $q(s) < g_{-}(s)=-\frac{1}{s^{a-\frac{1}{2}}}.$  
From this fact, we obtain for some $C>0$ the following
$$\frac{c_{p}^{a}}{s^{a}}+\mathit{o} (\frac{1}{s^{a}})\leq \frac{C}{s^{a}}< -\frac{Cq(s)}{s^{\frac{1}{2}}}$$
and deduce that 
\begin{equation} \label{e81}
\frac{c_{p}^{a}}{s^{a}}+\mathit{o} (\frac{1}{s^{a}})=\O(-\frac{q(s)}{s^{\frac{1}{2}}}).
\end{equation} 
Combining $\eqref{sai}$ and $\eqref{e81}$, we get
\begin{equation} \label{sai2}
q'(s)=q(s)+\O(-\frac{q(s)}{s^{\frac{1}{2}}})+\O (q^2(s))=q(s)\Big(1-\O(\frac{1}{s^{\frac{1}{2}}})+\O (q(s))\Big) .
\end{equation}
We use the fact that in this case $q(s)<0$ and $1-\O(\frac{1}{s^{\frac{1}{2}}})+\O (q(s))>\frac{1}{2}$, therefore
\begin{equation} \label{sai11}
q'(s)\leq \frac{q(s)}{2}.
\end{equation} 
Solving $\eqref{sai11}$ and using the fact that $q(\tilde{s})<0$, we derive the following 
$$q(s)\leq e^{\frac{s}{2}-\frac{\tilde{s}}{2}}q(\tilde{s})\rightarrow -\infty,\,\,\,{\rm as }\,\,\,s\longrightarrow +\infty.$$
 Recalling that $q$ is small, we get a contradiction.
\\{\bf The third case:} \underline{ If there exist $\tilde{s}$ large enough such that $ q(\tilde{s}) > g_{+}(\tilde{s})$.}
\\In this case, for all $s\geq \tilde{s}$, $g_{+}(s)=\frac{1}{s^{\frac{1}{2}}}<q(s)$, which implie that
$$ \frac{c_{p}^{a}}{s^{a}}+\mathit{o} (\frac{1}{s^{a}})=\O (q^2(s)).$$
We use again $\eqref{sai}$, to write
\begin{equation} \label{sai3}
q'(s)=q(s)+\O (q^2(s))= q(s)(1+\O (q(s)))\geq q(s).
\end{equation}
We solve $\eqref{sai3}$ and we use the fact that $q(\tilde{s})>0$, we derive the following
$$q(s)\geq e^{s-\tilde{s}}q(\tilde{s})\rightarrow +\infty,\,\,\,{\rm as }\,\,\,s\longrightarrow +\infty,$$
which is a contradiction because $q$ is small. Consequently, we get
\begin{equation} \label{82}
q(s)=\phi(s)-\kappa_{0} \sim  -\frac{\kappa_{0}}{p-1}\Big(\frac{p-1}{4s}\Big)^{a},\,\,\,{\rm as }\,\,\,s\longrightarrow +\infty.
\end{equation} 
 This conclude the proof of Claim $\ref{T14}$.
 
 \Box

\bigskip
 
 In the rest of this section, we are going to give the proof of some identities used in many steps of this paper. We claim the following:

\begin{cor}\label{T15}
For all $(d,y)\in (-1,1)^2$ and $s\in [s^*,+\infty )$, we have the following:
\begin{equation}\label{e51}
|\partial_{d} \overline{w}_{2}(d,y,s)| \leq \frac{C\kappa(d,y)}{s^{a}(1-d^2)},
\end{equation}
\begin{equation}\label{ap1}
|\partial_{d}\overline{w}_{1}(d,y,s)|+|\partial_{d}\overline{w}_{2}(d,y,s)|+\Vert \partial_{d}\overline{\psi}(d,y,s)\Vert_{L_{\rho}^{\frac{p+1}{p-1}}}\leq \frac{C}{1-d^2},
\end{equation}
\begin{equation}\label{ap2}
|\partial_{s}\overline{w}_{1}(d,y,s)|\leq \frac{C}{s^a}\kappa(d,y),
\end{equation}
\begin{equation}\label{e54}
|\partial_{d}(\tilde{\phi}(d,y,s))|+\Vert \partial_{s}\overline{\psi}(d,y,s)\Vert_{L_{\rho}^{\frac{p+1}{p-1}}}\leq \frac{C}{s^{a}(1-d^2)}.
 \end{equation} 
\end{cor}
{\it Proof}: 
\\{\bf Proof of $\eqref{e51}$:} Some simple calculations give the following
\begin{eqnarray}\label{hsat}
 \partial_{d} \overline{w}_{2}(d,y,s)&=&\frac{1}{\kappa_{0}}\partial_{d}(\kappa(d,y))\phi'(s-\log\Big(\frac{1+dy}{\sqrt{1-d^2}}\Big))\nonumber\\
 &&-\kappa(d,y)\phi''(s-\log\Big(\frac{1+dy}{\sqrt{1-d^2}}\Big))\frac{y+d}{(1+dy)(1-d^2)}.
 \end{eqnarray}
From $(ii)$ of Claim $\ref{T14}$, the expression $\eqref{hsat}$ of $\partial_{d} \overline{w}_{2}(d,y,s)$ and inequality $\eqref{zet}$, we conclude the proof of $\eqref{e51}$.
\\{\bf Proof of $\eqref{ap1}$:}
By a careful calculation, we write
\begin{eqnarray}\label{hs40}
 \partial_{d}\overline{w}_{1}(d,y,s)&=&-\frac{\kappa(d,y)}{\kappa_{0}(1+dy)}\Big[ \frac{2d}{(p-1)}\phi(s-\log \Big(\frac{1+dy}{\sqrt{1-d^2}}\Big))\nonumber\\
 &&-\frac{y+d}{1-d^2}\phi'(s-\log \Big(\frac{1+dy}{\sqrt{1-d^2}}\Big))\Big ].
 \end{eqnarray}
From the expression $\eqref{hs40}$ of $\partial_{d}\overline{w}_{1}(d,y,s)$ and Claim $\ref{T14}$, we can see that:
\begin{equation}\label{ap20}
|\partial_{d}\overline{w}_{1}(d,y,s)|\leq \frac{C}{s^a(1-d^2)}\Big [\frac{1-d^2}{1+dy}+\frac{|y+d|}{1+dy}\Big ].
\end{equation}
Using $\eqref{ap20}$, $\eqref{e104}$, the expression $\eqref{hsat}$ of $\partial_{d}\overline{w}_{2}(d,y,s)$, the fact that $\partial_{d}\kappa(d,y)\leq \frac{C}{1-d^2}$, Claim $\ref{T14}$ and the fact that $s^*$ is large enough such that for all $s\geq s^*$, we have:
\begin{equation}\label{ap21}
|\partial_{d}\overline{w}_{1}(d,y,s)|+|\partial_{d}\overline{w}_{2}(d,y,s)| \leq \frac{C}{1-d^2}.
\end{equation}
According to the expression $\eqref{e700}$ of $\overline{\psi}(d,y,s)$, we write for all $(d,y,s)\in (-1,1)^2\times [s^*,+\infty )$, we obtain 
\begin{equation}\label{ap60}
|\partial_{d}\overline{\psi}(d,y,s)|\leq \frac{C}{(1+dy)^2}.
\end{equation}
 We apply $(iii)$ of Claim $\ref{T20}$ and $\eqref{ap60}$ to conclude that 
 \begin{equation}\label{50000}
\Vert \partial_{d}\overline{\psi}(d,y,s)\Vert_{L_{\rho}^{\frac{p+1}{p-1}}}\leq \frac{C}{(1+dy)^2}.
\end{equation}
 Combining $\eqref{ap21}$ and $\eqref{50000}$ to conclude the proof of $\eqref{ap1}$.
\\{\bf Proof of $\eqref{ap2}$:} Inequality $\eqref{ap2}$ is a direct consequence of item $(ii)$ of Claim $\ref{T14}$.
\\{\bf Proof of $\eqref{e54}$:}
We write in the following the expression of $\partial_{d}(\tilde{\phi}(d,y,s))$
\begin{equation}\label{ap80}
\partial_{d}(\tilde{\phi}(d,y,s))=-\frac{1}{\kappa_{0}}\phi'\Big(s-\log\Big(\frac{1+dy}{\sqrt{1-d^2}}\Big)\Big)\frac{y+d}{(1+dy)(1-d^2)}.
\end{equation}
From $\eqref{ap80}$, we have:
\begin{equation}\label{ap50}
|\partial_{d}(\tilde{\phi}(d,y,s))|\leq \frac{1}{\kappa_{0}}|\phi'(s-\log\Big(\frac{1+dy}{\sqrt{1-d^2}}\Big))|\frac{y+d}{(1+dy)(1-d^2)^{ \frac{3}{2}}}\leq \frac{C}{s^{a}(1-d^2)}.
\end{equation}
From the expression $\eqref{e700}$ of $\overline{\psi}(d,y,s)$ and $(ii)$ of Claim $\ref{T14}$, we can see that 
\begin{equation}\label{ap30}
|\partial_{s}\overline{\psi}(d,y,s)|\leq \frac{C}{s^a(1+dy)^2}.
\end{equation} 
Then, we apply $(iii)$ of Claim $\ref{T20}$ combined with $\eqref{ap50}$ and $\eqref{ap30}$ to deduce $\eqref{e54}$.

\Box

\medskip

We recall in the following claim some basic bounds and properties used in many steps in our paper and exceptionally in the proof of the Proposition $\ref{T7}$.
\begin{cl}\label{T17}{\bf (Useful properties):}
For all $(d,d_{1},d_{2},y)\in (-1,1)^4$, we have the following:
\begin{equation}\label{e104}
|y+d|+(1-d^2)+(1-y^2)\leq C|1+dy|,
\end{equation}
\begin{equation}\label{zet}
|\partial_{d}(\kappa(d,y))| \leq \frac{C\kappa(d,y)}{1-d^2},
\end{equation}
\begin{equation}\label{47}
|W_{\lambda ,2}^{d}(y)|+(1-y^2)|\partial_{y}W_{\lambda ,2}^{d}(y)|\leq C\kappa(d,y),
\end{equation}
\begin{equation}\label{e10t}
 \Vert \overline{w}_{1}(d_{1},y,s)-\overline{w}_{1}(d_{2},y,s)\Vert_{\H_{0}}\leq C|\theta_{1}-\theta_{2}|,
 \end{equation}
 \begin{equation}\label{e15t}
 \Vert \overline{w}_{2}(d_{1},y,s)-\overline{w}_{2}(d_{2},y,s)\Vert_{\L^2_{\rho}}\leq C|\theta_{1}-\theta_{2}|,
 \end{equation}
where $W_{\lambda ,2}^{d}$ is defined in $\eqref{e4}$ and $\theta_{i}=\frac{1}{2}\log \Big(\frac{1+d_{i}}{1-d_{i}}\Big )$, for $i\in \lbrace 1,2\rbrace$.
\end{cl}
{\it Proof}: 
\\{\bf Proof of $\eqref{zet}$:} The proof is ommited since it is classical and known from \cite{fh}.
\\{\bf Proof of $\eqref{47}$:}
The proof is ommited since it is the same as the proof of $(iii)$ of Lemma 4.4 page 85 in \cite{fh}.
\\{\bf Proof of $\eqref{e10t}$:}
We exploite here the proof of inequality $(174)$ in Merle and Zaag \cite{fh} page $102$. Since, we have 
from item $(i)$ of Claim $\ref{T14}$ to deduce:
\begin{equation}\label{ap25}
 \overline{w}_{1}(d_{1},y,s)-\kappa(d_{1},y)\leq -\frac{\kappa_{0}}{p-1}\Big(\frac{p-1}{4s}\Big)^{a} \kappa(d_{1},y)
\end{equation}
and 
\begin{equation}\label{ap26}
\kappa(d_{2},y)-\overline{w}_{1}(d_{2},y,s)\leq \frac{\kappa_{0}}{p-1}\Big(\frac{p-1}{4s}\Big)^{a} \kappa(d_{2},y).
\end{equation}
Combining $\eqref{ap25}$, $\eqref{ap26}$ and the triangular inequality to conclude that 
\begin{equation}\label{ap100}
\Vert \overline{w}_{1}(d_{1},y,s)-\overline{w}_{1}(d_{2},y,s)\Vert_{\H_{0}}\leq \frac{C}{s^{a}}\Vert \kappa(d_{1},y)-\kappa(d_{2},y)\Vert_{\H_{0}}.
\end{equation}
Coming at this level, we apply the proof of inequality $(174)$ in Merle and Zaag \cite{fh} page $102$ to obtain $\eqref{e10t}$.
\\{\bf Proof of $\eqref{e15t}$:} The proof is similar as the proof of $\eqref{e10t}$. According to the expression of $\overline{w}_{2}(d,y,s)$ and item $(ii)$ of Claim $\ref{T14}$, we obtain:
 \begin{equation}\label{ap28}
\Vert \overline{w}_{2}(d_{1},y,s)-\overline{w}_{2}(d_{2},y,s)\Vert_{\L^2_{\rho}}\leq \frac{C}{s^{a}}\Vert \kappa(d_{1},y)-\kappa(d_{2},y)\Vert_{\L^2_{\rho}}.
\end{equation}
It is clear that we have
\begin{equation}\label{ap29}
\Vert \kappa(d_{1},y)-\kappa(d_{2},y)\Vert_{\L^2_{\rho}}\leq \Vert \kappa(d_{1},y)-\kappa(d_{2},y)\Vert_{\H_{0}}.
\end{equation}
Now, we are in position to apply the proof of inequality $(174)$ in Merle and Zaag \cite{fh} page $102$ to obtain $\eqref{e15t}$.

\Box
\bigskip

\medskip
 We recall now from \cite {fh} the following estimate:

\begin{cl}{\rm (Integral computation table )}.\label{T20}
Consider for some $\alpha >-1$ and $\beta \in \R$ the following integral,
$$I(d)=\int_{-1}^{1}\frac{(1-y^2)^{\alpha}}{(1+dy)^{\beta}}dy,$$
then there exists $K(\alpha ,\beta)$ such that the following holds for all $d\in (-1,1)$,
\\ $(i)$ if $\alpha +1-\beta >0$, then $\frac{1}{K}\leq I(d)\leq K$,
\\ $(ii)$ if $\alpha +1-\beta =0$, then $\frac{1}{K}\leq I(d)/\log (1-d^2) \leq K$,
\\ $(iii)$ if $\alpha +1-\beta <0$, then $\frac{1}{K}\leq I(d) (1-d^2)^{-(\alpha +1)+\beta} \leq K$.
\end{cl}
{\it Proof}: 
See the proof of Claim 4.3 page $84$ in Merle and Zaag \cite {fh}.

\Box

\bigskip

We would like now to recall that the Lorentz transform keeps equation $\eqref{1}$ and the norms in $\H_{0}$ invariant. More precisely, we write from Merle and Zaag \cite {fh} the following lemma:
\begin{lem}{\rm (The invariance of the Lorentz transform in similarity variables)}\label{Tekk}Consider $w(y,s)$ a solution of equation $\eqref{1}$ defined for all $|y|<1$ and $s\in (s_{0},s_{1})$ for some $s_{0}$ and $s_{1}$ in $\R$ and introduce for any $d\in (-1,1)$, the function $W\equiv \tau_{d}(w)$ defined by,
$$W(Y,S)=\frac{(1-d^2)^{\frac{1}{p-1}}}{(1+dY)^{\frac{2}{p-1}}}w(y,s),\,\,\,{\rm where}$$
$$y=\frac{Y+d}{1+dY}\,\,\,{\rm and}\,\,\,s=S-\log \Big (\frac{1+dY}{\sqrt{1-d^2}}\Big ).$$
Then $W(Y,S)=\tau_{d}(w)$ is also a solution of $\eqref{1}$ defined for all $|y|<1$ and 
$$S\in \Big( s_{0}+\frac{1}{2}\log \frac{1+|d|}{1-|d|},s_{1}-\frac{1}{2}\log \frac{1+|d|}{1-|d|} \Big ).$$
Also, we have the continuity of $\tau_{d}$ in $\H_{0}$, 
$$\Vert \tau_{d}(w)\Vert_{\H_{0}}\leq C\Vert w\Vert_{\H_{0}}.$$
\end{lem}
{\it Proof}: See the proof of Lemma 2.6 page 54 and Lemma 2.8 page 57 in Merle and Zaag \cite {fh}.

\Box

We end this section by the Hardy-Sobolev identity in the space $\H_{0}$ defined in $\eqref{e91}$ combined with some basic bounds on the stationary solution of $\eqref{1}$ when $f\equiv 0$, namely $\kappa (d,y)$ defined in $\eqref{400}$: 
\begin{lem}\label{T21}We have the following identities,
\begin{enumerate}[{\rm(i)}]
\item {\bf (A Hardy-Sobolev type identity)}For all $h\in \H_{0}$, we have
$$\|h\|_{L_{\frac{\rho}{1-y^2}}^{2}(-1,1)}+\|h\|_{L_{\rho}^{p+1}(-1,1)}+\|h(1-y^2)^{\frac{1}{p-1}}\|_{L^{\infty}(-1,1)}\leq C \Vert h\Vert_{\H_{0}}.$$
\item{\bf (Boundedness of $\kappa (d,y)$ in several norms)}For all $d\in (-1,\,1)$, we have
$$\|\kappa (d,y)\|_{L_{\frac{\rho}{1-y^2}}^{2}(-1,1)}+\|\kappa (d,y)\|_{L_{\rho}^{p+1}(-1,1)}+\|\kappa (d,y)(1-y^2)^{\frac{1}{p-1}}\|_{L^{\infty}(-1,1)}+ \Vert \kappa (d,y)\Vert_{\H_{0}}\leq C.$$
\end{enumerate}
\end{lem}
{\it Proof}: For the proof of $(i)$, we can see the proof of Lemma 2.2 page $51$ in Merle and Zaag \cite {fh}. For the proof of $(ii)$, we need to use $(i)$ and identity $(49)$ page $(59)$ in \cite {fh}.
 
\Box

\noindent{\bf Address}:\\
 Imam Abdulrahman Bin Faisal University P.O. Box 1982 Dammam, Saudi Arabia.
\vspace{-4mm}
\begin{verbatim}
e-mail: mahamza@iau.edu.sa
\end{verbatim}
Facult\'e Des Sciences de Gafsa. Laboratoire \'{E}quations Aux Dérivées Partielles LR03ES04, 2092 Tunis, Tunisie.
\vspace{-4mm}
\begin{verbatim}
e-mail: saidiomar1985@gmail.com
\end{verbatim}


\begin{thebibliography}{10}
 
\bibitem{sa} 
S. Alinhac.
\newblock Blow up for nonlinear hyperbolic equations,.
\newblock In {\em volume~17 of Progress in Nonlinear Differential Equations and
  their Applications}, pages Birkh\"auser Boston Inc., Boston, MA. 1995.

\bibitem{sa1}
S. Alinhac.
\newblock A minicourse on global existence and blowup of classical solutions to
  multidimensional quasilinear wave equations.
\newblock In {\em Journ\'ees ``\'{E}quations aux {D}\'eriv\'ees {P}artielles''
  ({F}orges-les-{E}aux, 2002)}, pages Exp. No. I, 33. Univ. Nantes, Nantes,
  2002.

\bibitem{cf}
C.~Antonini and F.~Merle.
\newblock Optimal bounds on positive blow-up solutions for a semilinear wave
  equation.
\newblock {\em Internat. Math. Res. Notices}, (21):1141--1167, 2001.

\bibitem{BR}
A.~Bressan.
\newblock On the asymptotic shape of blow-up.
\newblock {\em Indiana Univ. Math. J.}, 39(4):947--960, 1990.

\bibitem{BR1}
A.~Bressan.
\newblock Stable blow-up patterns.
\newblock {\em J. Differential Equations.}, 98(1):57--75, 1992.

\bibitem{LA}
L.~A. Caffarelli and A. Friedman.
\newblock Differentiability of the blow-up curve for one dimensionel nonlinear wave equations.
\newblock {\em Arch. Rational Mech. Anal.}, 91(1):83--98, 1985.

\bibitem{LA1}
L.~A. Caffarelli and A. Friedman.
\newblock The blow-up boundary for nonlinear wave equations.
\newblock {\em Trans. Amer. Math. Soc.}, 297(1):223--241, 1986.

\bibitem{C}
R.~C\^{o}te.
\newblock Construction of solutions to the subcritical gKdV equations with a given asymptotical behavior.
\newblock {\em J. Funct. Anal.}, 241(1):143--211, 2006.

\bibitem{C1}
R.~C\^{o}te.
\newblock Construction of solutions to $L^2$-critical KdV equations with a given asymptotical behavior.
\newblock {\em Duke Math. J.}, 138(3):487--531, 2007.

\bibitem{C2}
R.~C\^{o}te and H. Zaag.
\newblock Construction of a multisoliton blowup solution to the semilinear wave equation in one space dimension.
\newblock {\em Comm. Pure Appl. Math.}, 66(10):1541--1581, 2013.

\bibitem{DSS}
R. Donninger, W. Schlag and A. Soffer.
\newblock On pointwise decay of linear waves on a Schwarzschild black hole background
\newblock {\em Comm. Math. Phys.}, 309(1):51--86, 2012.

\bibitem{DS}
R. Donninger and B. Sch\"{o}rkhuber.
\newblock Stable self-similar blow-up for energy subcritical wave equation.
\newblock {\em Dyn. Partial Differ. Equ.}, (9):63--87, 2012.

\bibitem{DS1}
R. Donninger and B. Sch\"{o}rkhuber.
\newblock Stable blow-up dynamics for energy supercritical wave equations.
\newblock {\em Trans. Amer. Soc.}, 366(4):2167--2189, 2014.

\bibitem{GIN}
T. Ghoul, S. Ibrahim and V. T. Nguyen.
\newblock Construction of type II blow up solutions for the 1-corotational energy supercritical wave maps.
\newblock {\em J. Differential Equations.}, 265(7):2968--3047, 2018.


\bibitem{GHM}
T. Ghoul and N. Masmoudi.
\newblock Stability of infinite time blow up for the Patlak-Keller-Segel system. [To appear in Communication of Pure and Applied Mathematics].


\bibitem{TVH}
T. Ghoul, V. T. Nguyen and H. Zaag.
\newblock Construction and stability of type I blow up solutions for non-variational semilinear parabolic systems.
\newblock {\em Adv. Pure Appl. Math.}, 10(4):299--312, 2019.

\bibitem{TVH1}
T. Ghoul, V. T. Nguyen and H. Zaag.
\newblock Construction and stability of type I blow up solutions for a higher order semilinear parabolic systems. 
\newblock {\em Adv. Nonlinear Anal.}, 9(1):388--412, 2020.



\bibitem{GMS}
Y. Giga, S. Matsui and S. Sasayama.
\newblock Blow-up rate for semilinear heat equations with subcritical nonlinearity .
\newblock {\em Indiana Univ. Math.J}, 53(2):483--514, 2004.
\bibitem{XXX}
 M. A. Hamza.
 \newblock The blow-up rate for strongly perturbed semilinear wave equations in the conformal regime without a radial assumption.
 \newblock {\em  Asymptotic Analysis.},  97, no. 3-4, 351-378, 2016.

\bibitem{X}
 M. A. Hamza and O. Saidi.
 \newblock The blow-up rate for strongly perturbed semilinear wave equations.
 \newblock {\em  J. Dyn. Diff. Equat}, 26(1): 2014.

\bibitem{XX}
 M. A. Hamza and O. Saidi.
 \newblock The blow-up rate for strongly perturbed semilinear wave equations in the conformal case.
 \newblock {\em  Math. Phys. Anal. Geom.}, 18(1), Art. 15, 2015.



\bibitem{MH}
M.~A. Hamza and H.~Zaag.
\newblock A {L}yapunov functional and blow-up results for a class of perturbed
  semilinear wave equations.
\newblock {\em Nonlinearity}, 25(9):2759--2773, 2012.


\bibitem{MH1}
M.~A. Hamza and H.~Zaag.
\newblock Lyapunov functional and blow-up results for a class of
              perturbations of semilinear wave equations in the critical
              case.
\newblock {\em J. Hyperbolic Differ. Equ}, 9(2):195--221, 2012.

\bibitem{MH2}
M.~A. Hamza and H.~Zaag.
\newblock  Blow-up behavior for the Klein-Gordon and other perturbed semilinear wave equations.
\newblock {\em  Bull. Sci. math}, 137(8):1087--1109, 2013.

\bibitem{MH3}
M.~A. Hamza and H.~Zaag.
\newblock  Prescribing the center of mass of a multi soliton solution for a perturbed semilinear wave equation.
\newblock {\em  J. Differ. Equ}, 267(6):3524--3560, 2019.

\bibitem{HZN}
M.~A. Hamza and H.~Zaag.
\newblock  The blow-up rate for non-scaling invariant semilinear wave equations.
\newblock {\em  J. Math. Annal. Appl}, 483(2), 2020.

\bibitem{HZNn}
M.~A. Hamza and H.~Zaag.
\newblock  The blow-up rate for non-scaling invariant semilinear wave equations in higher dimensions.
\newblock {\em  Nonlinear Analysis} Volume 212, 2021.

\bibitem{HZNnn}
M.~A. Hamza and H.~Zaag.
\newblock  The blow-up rate for non-scaling invariant semilinear heat equation.
\newblock {\em  arXiv:2012.00768v1} (Submitted), 2021.




\bibitem{SW}
S. Kichenassamy and W. Littman.
\newblock Blow-up surfaces for nonlinear wave equations.
\newblock {\em I. Comm. Partial Differential Equations.}, 18(3-4):431--452, 1993.

\bibitem{SW1}
S. Kichenassamy and W. Littman.
\newblock Blow-up surfaces for nonlinear wave equations.
\newblock {\em II. Comm. Partial Differential Equations,} 18(11):1869--1899, 1993.


\bibitem{RM}
R. Killip, B. Stovall and M. Visan.
\newblock Blow-up behaviour for the nonlinear Klein-Gordon equation.
\newblock {\em Math. Ann,} 358:289--350, 2014.


\bibitem{ha}
H. A. Levine.
\newblock Instability and non-existence of global solutions to nonlinear wave equations of the form {$Pu\sb{tt}=-Au+{\cal F}(u)$}.
\newblock {\em Trans. Amer. Math. Soc.}, 192:1--21, 1974.


\bibitem{HAG}
H. A. Levine and G. Todorova.
\newblock Blow up of solutions of the {C}auchy problem for a wave equation with
  nonlinear damping and source terms and positive initial energy.
\newblock {\em SIAM J. Math. Anal.}, 5(3):793--805, 2001.

\bibitem{MNZ}
F.~Mahmoudi, N. Nouaili and H.~Zaag.
\newblock Construction of a stable periodic solution to a semilinear heat equation with a prescribed profile.
\newblock {\em Nonlinear Anal.}, 131:300--324, 2016.


\bibitem{mm}
Y.~Martel and F.~Merle.
\newblock Stability of blow-up profile and lower bounds for blow-up rate for the critical generalized KDV equation.
\newblock {\em Ann. of Math. (2)}, 155(1):235--280, 2002.

\bibitem{MAZ}
N.~Masmoudi and H.~Zaag.
\newblock Blow-up profile for the complex Ginzerburg-Landau equation.
\newblock {\em J. Funct. Anal.}, 255(7):1613--1666, 2008.

\bibitem{mr}
F.~Merle and P.~Rapha\"{e}l.
\newblock On universality of blow-up profile for $L^2$ critical nonlinear Shr\"{o}dinger equation.
\newblock {\em Invent. Math.}, 156(3):565--672, 2004.

\bibitem{MRR}
F.~Merle, P.~Rapha\"{e}l and I. Rodnianski.
\newblock Blowup dynamics for smooth data equivariant solutions to the critical Shr\"{o}dinger map problem.
\newblock {\em Invent. Math.}, 193(2):249--365, 2013.


\bibitem{fh3}
F.~Merle and H.~Zaag.
\newblock Determination of the blow-up rate for the semilinear wave equation.
\newblock {\em Amer. J. Math.}, 125(5):1147--1164, 2003.

\bibitem{fh4}
F.~Merle and H.~Zaag.
\newblock Determination of the blow-up rate for a critical semilinear wave
  equation.
\newblock {\em Math. Ann.}, 331(2):395--416, 2005.



\bibitem{kl}
F.~Merle and H.~Zaag.
\newblock Blow-up rate near the blow-up surface for semilinear wave equation.
\newblock {\em Internat. Math. Res. Notices}, 19(1):1127--56, 2005.

\bibitem{fh6}
F.~Merle and H.~Zaag.
\newblock Existence and classification of characteristic points at blow-up for
  a semilinear wave equation in one space dimension.
\newblock {\em Amer. J. Math.}, 134(3):581--648, 2012.

\bibitem{fh}
F.~Merle and H.~Zaag.
\newblock Existence and universality of the blow-up profile for the semilinear
  wave equation in one space dimension.
\newblock {\em J. Funct. Anal.}, 253(1):43--121, 2007.

\bibitem{fh1}
F.~Merle and H.~Zaag.
\newblock Openness of the set of non-characteristic points and regularity of
  the blow-up curve for the 1 {D} semilinear wave equation.
\newblock {\em Comm. Math. Phys.}, 282(1):55--86, 2008.

\bibitem{fh7}
F.~Merle and H.~Zaag.
\newblock Isolatedness of characteristic points for a semilinear wave equation in one space dimension. In
\newblock {S\'{e}minaire sur les Equation aux d\'{e}riv\'{e}es partielle 2009-2010,} pages Exp.No. 11, 10p. Ecole Polytech., Palaiseau 2010.

\bibitem{fh15}
F.~Merle and H.~Zaag.
\newblock On the stability of the notion of non-characteristic points and blow-up profile for semilinear wave equations.  
\newblock {\em Comm. Math. Phys.}, pages 1-34, 2015.
\bibitem{fh10}
F.~Merle and H.~Zaag.
\newblock Dynamics near explicit stationary solutions in similarity variables for solutions of a semilinear wave equation in higher dimensions. 
\newblock {\em Trans. Amer. Math. Soc}, 368(1):27--87, 2016.

\bibitem{rhv}
V. T. Nguyen and H.~Zaag.
\newblock Blow-up results for a strongly perturbed semilinear heat equation: Theoretical analysis and numerical method.
\newblock {\em Anal. PDE.}, 9(1):229--257, 2016.


\bibitem{VH}
V. T. Nguyen and H.~Zaag.
\newblock Construction of a stable blow-up solution for a class of strongly perturbed semilinear heat equations.
\newblock {\em Ann. Sc. Norm. Super. Pisa Cl. Sci.}, 16(4):1275--1314, 2016.

\bibitem{NH1}
N.~Nouaili and H.~Zaag.
\newblock Profile for a Simultaneously Blowing up Solution to a Complex Valued Semilinear Heat Equation.
\newblock {\em Comm. Partial Differential Equations.}, 40(7):1197--1217, 2015.

\bibitem{NH}
N.~Nouaili and H.~Zaag.
\newblock Construction of a blow-up solutions for the complex Ginzburg-Landau in a critical case.
\newblock {\em Arch. Ration. Mech. Anal.}, 228(3):995--1058, 2018.


\bibitem{RS}
P.~Rapha\"{e}l and R. Schweyer.
\newblock On the stability of critical chemotactic aggregation.
\newblock {\em Math. Ann.}, 359(12):267--377, 2014.




\bibitem{GT}
G. Todorova.
\newblock Cauchy problem for a non linear wave equation with non linear damping
  and source terms.
\newblock {\em Nonlinear Anal.}, pages 891--905, 2000.

\bibitem{W1}
G. B. Whitham.
\newblock Linear and nonlinear waves.
\newblock {\em Pure and Applied Mathematics (New York), John Wiley and Sons Inc.}, New York, 1999.
\\Reprint of the 1974 original, A Wiley-Interscience Publication.

\bibitem{Z}
H.~Zaag.
\newblock Blow-up results for vector-valued nonlinear heat equations with no gradient structure.
\newblock {\em Ann. Inst. H. Poincarr\'{e}. Anal. Non Lin\'{e}aire.}, 15(5):581--622, 1998.

 
 \end{thebibliography}
\end{document}